\documentclass{article}%
\usepackage{amsmath}
\usepackage{amsfonts}
\usepackage{amssymb}
\usepackage{graphicx}%
\usepackage{mathrsfs}
\usepackage{algorithm}
\usepackage{algorithmic}
\usepackage{subfigure}
\usepackage{epsfig}
\usepackage{booktabs}
\usepackage{multirow}
\usepackage{rotating}

\providecommand{\U}[1]{\protect \rule{.1in}{.1in}}
\newtheorem{theorem}{Theorem}

\newtheorem{lemma}{Lemma}

\newtheorem{remark}{Remark}

\newtheorem{ass}{Assumption}
\newenvironment{proof}[1][Proof]{\noindent \textbf{#1.} }{\  \rule{0.5em}{0.5em}}

\begin{document}

\title{ Local Nonparametric Estimation for Second-Order Jump-Diffusion Model Using Gamma Asymmetric Kernels
\footnote{This research work is supported by the National Natural
Science Foundation of China (No. 11371317, 11526205,  11626247), the
Fundamental Research Fund of Shandong University (No. 2016GN019) and
the General Research Fund of Shanghai Normal University (No.
SK201720).} }

\date{\today}

\author{Yuping Song\\
School of Finance and Business, Shanghai Normal
         University,\\ Shanghai, 200234, P.R.C.\\~\\Hanchao Wang\footnote{Corresponding author, email: wanghanchao@sdu.edu.cn}\\
Zhongtai Securities Institute for Financial Studies, Shandong
University,\\ Jinan, 250100, P.R.C.\\}

\date{}
\maketitle

\begin{abstract}
This paper discusses the local linear smoothing to estimate the
unknown first and second infinitesimal moments in second-order
jump-diffusion model based on Gamma asymmetric kernels. Under the
mild conditions, we obtain the weak consistency and the asymptotic
normality of these estimators for both interior and boundary design
points. Besides the standard properties of the local linear
estimation such as simple bias representation and boundary bias
correction, the local linear smoothing using Gamma asymmetric
kernels possess some extra advantages such as variable bandwidth,
variance reduction and resistance to sparse design, which is
validated through finite sample simulation study. Finally, we employ
the estimators
for the return of some high frequency financial data.\\


\end{abstract}

\section{Introduction}

Continuous-time models are widely used in economics and finance,
such as interest rate and etc., especially the continuous-time
diffusion processes with jumps. Jump-diffusion process $X_{t}$ is
represented by the following stochastic differential equation:
\begin{equation}
\label{1.1} dX_{t} = \mu(X_{t-})dt + \sigma(X_{t-})dW_{t} +
\int_{\mathscr{E}}c(X_{t-} , z) r(\omega, dt, dz), \end{equation}
which can accommodate the impact of sudden and large shocks to
financial markets. Johannes \cite{jo} provided the statistical and
economic role of jumps in continuous-time interest rate models.
However, in empirical finance or physics the current observation
usually behaves as the cumulation of all past perturbation such as
stock prices by means of returns and exchange rates in Nicolau
\cite{ni2} and the velocity of the particle on the surface of a
liquid in Rogers and Williams \cite{rw}. Furthermore, in the
research field involved with model (\ref{1.1}), the scholars mainly
considered it for the price of a asset not for the returns of the
price. As mentioned in Campbell, Lo and MacKinlay \cite{clm}, return
series of an asset are a complete and scale-free summary of the
investment opportunity for average investors, and are easier to
handle than price series due to their more attractive statistical
properties.

For characterizing this integrated economic phenomenon, moreover,
and the return series, Nicolau \cite{ni} considered the promising
continuous second-order diffusion process (\ref{nm}), which is
motivated by unit root processes under the discrete framework of
Park and Phillips \cite{pp},
\begin{equation}
\label{nm} \left\{
\begin{array}{ll}dY_{t} = X_{t}dt,\\
dX_{t} = \mu(X_{t})dt + \sigma(X_{t})dW_{t}, \end{array} \right.
\end{equation}
which may be an alternative model to describe the dynamic of some
financial data. Considering empirical properties of return series
such as higher peak, fatter tails and etc., here we extend the
continuous autoregression of order two (\ref{nm}) to a discontinuous
and realistic case with jumps for return series based on model
(\ref{1.1}), noted as the second-order jump-diffusion process
(\ref{1.2}) to represent integrated and differentiated processes,
\begin{equation}
\label{1.2} \left\{
\begin{array}{ll}
dY_{t} = X_{t-}dt,\\
dX_{t} = \mu(X_{t-})dt + \sigma(X_{t-})dW_{t} +
\int_{\mathscr{E}}c(X_{t-} , z) r(\omega, dt, dz),
\end{array}
\right.
\end{equation}
where $W_{t}$ is a standard Brownian motion, $\mu(\cdot)$ and
$\sigma(\cdot)$ are the infinitesimal conditional drift and variance
respectively, $\mathscr{E} = \mathbb{R}\setminus\{ 0\},$ $r(\omega,
dt, dz) = (p - q)(dt, dz), ~p(dt, dz)$ is a time-homogeneous Poisson
random measure on $\mathbb{R}_{+} \times \mathbb{R}$ independent of
$W_{t}$, and $q(dt, dz)$ is its intensity measure, that is, $E[p(dt,
dz)] = q(dt, dz) = f(z)dzdt$, $f(z)$ is its L\'{e}vy density. For
empirical financial data, $X_{t}$ represents the continuously
compounded return of underlying assets, $Y_{t}$ denotes the asset
price by means of the cumulation of the returns plus initial asset
value. Note that model (\ref{1.2}) is neither a special case of a
two-dimensional stochastic differential equation nor a stochastic
volatility model without noise in the first coordinate. It is
essentially the linear version of the following nonlinear stochastic
jump-diffusion equation of order two considered in \"{O}zden and
\"{U}nal \cite{ou}
\begin{equation}
\label{1.3} d\dot{X}_{t} = f(t, X_{t}, \dot{X}_{t})dt + g(t, X_{t},
\dot{X}_{t})dW_{t} + j(t, X_{t}, \dot{X}_{t})dN_{t},
\end{equation}
where $d N_{t}$ denotes the infinitesimal increment of Poisson
process and $f, g, j$ are the coefficient functions.

Model (\ref{1.2}) can be used commonly in empirical financial for at
least four reasons. Firstly, in contrary to the usual models for the
price of a asset such as model (\ref{1.1}), model (\ref{1.2})
overcomes the difficulties associated with the nondifferentiability
of a Brownian motion, which can model integrated and differentiated
diffusion processes for the cumulation of all past perturbations in
modern econometric phenomena (similarly as unit root processes in a
discrete framework). Furthermore, the model (\ref{1.2}) can
accommodate nonstationary original process and be made stationary by
linear combinations such as differencing, which is a more widely
adopted technique used in empirical financial. We have verified this
property through the Augmented Dickey-Fuller test statistic for the
sampling time series before and after the difference in the
empirical analysis part. Secondly, compared with the model
(\ref{nm}), model (\ref{1.2}) accommodates the impact by
macroeconomic announcements and a dramatic interest rate cut by the
Federal Reserve in combination with jump component. It has been
testified the existence of jumps modeling by (\ref{nm}) or
(\ref{1.2}) for real financial data through the test statistic
proposed in Barndorff-Nielsen and Shephard \cite{bs2} in empirical
analysis part. Thirdly, model (\ref{1.2}) can directly characterize
the returns with heavy tail (or the log return) of a asset to
specify general properties for returns (such as stationarity in the
mean and weakness in the autocorrelation et al.) more easily than a
diffusion univariate process for the price of a asset (such as stock
prices and nominal exchange rates). Fourthly, the
integro-differential jump-diffusion model (\ref{1.2}) can also be
employed in finance for integrated volatility in stochastic
volatility models with jumps.

For model (\ref{1.2}), \"{O}zden and \"{U}nal \cite{ou} gave the
linearization criterion in terms of coefficients for transforming a
nonlinear stochastic differential equations into linear ones via
invertible stochastic mappings in order to solve exact solutions.
However, we do not know beforehand the specific form of the model
coefficients to verify these criteria in practical applications, so
we should statistically estimate the unknown coefficients in model
(\ref{1.2}) for modeling the economic and financial phenomena based
on the observations. For continuous case (\ref{nm}), Gloter
\cite{g1} \cite{g2} and Ditlevsen and S{\o}rensen \cite{ds}
presented the parametric and semiparametric estimation from a
discretely observed samples. Recently, an inordinate amount of
attention has been focused on nonparametric methods in econometrics
due to the flexibility for handling the nonlinear conditional moment
estimation, not assuming its expression in contrary to parametric
estimation. Nicolau \cite{ni} and Wang and Lin \cite{wl} analyzed
the local nonparametric estimations using symmetric kernel and Hanif
\cite{hm2} studied Nadaraya-Watson estimators using Gamma asymmetric
kernel. For model (\ref{1.2}), Song \cite{syp1} considered
Nadaraya-Watson estimators for the unknown coefficients using
Gaussian symmetric kernel in high frequency data.

In the context of nonparametric estimator with finite-dimensional
auxiliary variables, local polynomial smoothing become an effective
smoothing method, which has excellent properties of full asymptotic
minimax efficiency achievement and boundary bias correction
automatically, one can refer to Fan and Gijbels \cite{fg2} for
better review. In general, the popular choices for kernels in
application of local polynomial estimators are symmetric and
compact. However, local nonparametric estimation constructed with
symmetric kernels for the nonnegative variables or nonnegative part
of the underlying variables in economics and finance is not
approximate for the region near the origin without a boundary
correction. Furthermore, for finite sample size, Seifert and Gasser
\cite{sg} found the problem of unbounded variance in sparse regions
for local polynomial smoothing employing a compact kernel and
proposed local increase of bandwidth in sparse regions of the design
to add more information to reduce the corresponding variability.
Fortunately, local nonparametric smoothing using asymmetric kernels
can effectively solve the above two major problems for nonparametric
estimation.

Besides the usual standard properties of the local linear estimation
such as the simple bias representation and boundary bias correction,
the local linear smoothing using Gamma asymmetric kernel with
unbounded support $[0, \infty)$ have some extra benefits as follows
such as variable bandwidth, variance reduction and resistance to
sparse design. Firstly, the Gamma asymmetric kernel is a kind of
adaptive and flexible smoothing, whose curve shapes can vary with
the smoothing parameter and the location of the design point similar
as the variable bandwidth methods, which is illustrated in FIG 1 for
a fixed smoothing parameter. As mentioned in Gospodinov and Hirukawa
\cite{gh}, unlike the variable bandwidth methods, the smoothing
method constructed with Gamma asymmetric kernel is achieved by only
a single smoothing parameter. Secondly, when the support of Gamma
asymmetric kernel matches the support of the curve of the function
to be estimated, and the curve has sparse regions, the finite
variance of the curve estimation decreases due to the fact that
their variances vary along with the location of the design point $x$
and the increase of the effective sample size for estimation.
Thirdly, Gamma asymmetric kernels are free from boundary bias
(allowing a larger bandwidth to pool more data) and achieve the
optimal rate of convergence in mean square error within the class of
nonnegative kernel estimators with order two. In conclusion,
asymmetric kernels is a combination of a boundary correction device
and a ``variable bandwidth'' method. For a review of application
using asymmetric kernels, one should refer to Chen \cite{ch1}
\cite{ch2} to estimate densities, Chen \cite{ch3} to estimate a
regression curve with bounded support, Bouezmarni and Scaillet
\cite{bs} to apply the asymmetric kernel density estimators to
income data, Kristensen \cite{kr} to consider the realized
integrated volatility estimation, Gospodinov and Hirukawa \cite{gh}
to propose an asymmetric kernel-based method for scalar diffusion
models of spot interest rates, compared with Stanton \cite{st}, to
show that asymmetric kernel smoothing is expected to reduce
substantially both the bias near the origin and the bias that occurs
for high values in the estimation of the drift coefficient. Hanif
\cite{hm1} addressed local linear estimation using Gamma asymmetric
kernels for the infinitesimal moments in model (\ref{1.1}).

\begin{figure*}[!htb]
\centering
\includegraphics[width=5in,height=2in]{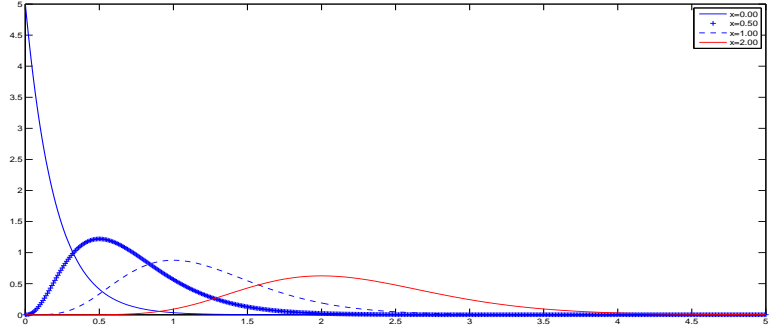}
\caption{Shapes of the Gamma function with a fixed smoothing
parameter (0.2) at various design points ($x = 0.0, 0.5, 1.0, 2.0$)}
\end{figure*}

For model (\ref{1.2}), although Chen and Zhang \cite{cz} discussed
the local linear estimators for the
 unknown functions $\mu(x)$ and $\sigma^{2}(x) +
\int_{\mathscr{E}}c^{2}(x , z)f(z)dz$ based on symmetric kernels.
More previous works focused on the simplification for the bias
representation of the estimator, while less research was considered
deeply for the reduction in variance. Gouri\'{e}roux and Monfort
\cite{gm} and Jones and Henderson \cite{jh} argued that unlike the
symmetric kernels case, the empirical data $X_{t}$ point and the
design point $x$ in asymmetric kernels case are not exchangeable.
Hence, we cannot establish the asymptotic theorems for the Gamma
kernel estimators of unknown coefficients in model (\ref{1.2}) by
the means of the similar approach as that in Bandi and Nguyen
\cite{bn} and Chen and Zhang \cite{cz}. In this paper, we will
propose local linear estimators of $\mu(x)$ and
$\sigma^2(x)+\int_{\mathscr{E}}c^{2}(x,z)f(z)dz$ in model
(\ref{1.2}) using Gamma asymmetric kernels for both bias correction,
especially for boundary point near the origin, and variance
reduction, especially for sparse design point far away from the
origin.

For convenience, we note $X_i = X_{i\Delta_n}$, $\widetilde{X}_i =
\widetilde{X}_{i\Delta_n}$ in the remainder of this paper which is
organized as follows. In Section 2, we propose local linear
estimators with asymmetric kernels and some ordinary assumptions for
model (\ref{1.2}). We present the asymptotic results in Section 3.
Section 4 presents the finite sample performance through Monte Carlo
simulation study. The estimators are illustrated empirically  in
Section 5. Section 6 concludes. Some technical lemmas and the main
proofs are explicitly shown in Section 7.

\section{Local Linear Estimators with Asymmetric Kernels and Assumptions}

According to the asymmetric kernels used in Chen \cite{ch2}, the
Gamma kernel function is defined as
\begin{equation}
\label{2.1} K_{G(x/h + 1, h)} (u) = \frac{u^{x/h} \exp
(-u/h)}{h^{x/h + 1} \Gamma(x/h + 1)}~~0 \leq u \leq \infty,
\end{equation}
where $\Gamma(m) = \int^{\infty}_{0}{y^{m - 1} \exp (- y) dy},~m >
0$ is the Gamma function and $h$ is the smoothing parameter. Note
that we use modified Gamma kernel function $K_{G(x/h + 1, h)} (u)$
instead of $K_{G(x/h, h)} (u)$ due to the fact $K_{G(x/h, h)} (u)$
is unbounded near at $x = 0.$ As is shown in Figure 1, the Gamma
function has shapes varying with the design point $x,$ which changes
the amount of smoothing applied by the asymmetric kernels since the
variance $x h + h^{2}$ of $K_{G(x/h + 1, h)} (u)$ is increasing as
$x$ away from the boundary. Additionally, as discussed in Bouezmarni
and Scaillet \cite{bs} the consistency of gamma kernel estimator
holds even though the true density is unbounded at $x = 0.$
Statistically and theoretically, since the density of Gamma
distribution has support $[0, \infty)$, the Gamma kernel function
does not generate boundary bias for nonnegative variables or
nonnegative part of underlying variables. Furthermore, the
asymptotic variance of the nonparametric Gamma kernel estimation
depends on the design point $x,$ which yields the optimal rate of
convergence in mean integrated squared error of nonnegative kernels.

Different from model (\ref{1.1}), nonparametric estimations
constructed for the coefficients in second-order jump-diffusion
model (\ref{1.2}) give rise to new challenges for two main reasons.

On the one hand, we usually get observations
$\{Y_{i\Delta_n};i=1,2,\cdot \cdot \cdot\}$ rather than
$\{X_{i\Delta_n};i = 1, 2, \cdot \cdot \cdot\}.$ The value of
$X_{t_{i}}$ cannot be obtained from $Y_{t_{i}} = Y_{0} +
\int_{0}^{t_{i}} X_{s} ds$ in a fixed sample intervals.
Additionally, nonparametric estimations of the unknown qualities in
model (\ref{1.2}) cannot in principle be constructed on the
observations $\{Y_{i\Delta_{n}} ; i = 1, 2, \cdot \cdot \cdot\}$ due
to the unknown conditional distribution of $Y$. As Nicolau \cite{ni}
showed, with observations $\{Y_{i\Delta_n};i=1,2,\cdot \cdot
\cdot\}$ and given that
$$Y_{{i\Delta_{n}}} - Y_{{(i-1)\Delta_{n}}} = \int_{(i-1)\Delta_{n}}^{i\Delta_{n}} X_{u} du, $$
we can obtain an approximation value of $X_{i\Delta_n}$ by
\begin{equation}
\label{2.2}
\widetilde{X}_{i\Delta_n}=\frac{Y_{i\Delta_n}-Y_{(i-1)\Delta_n}}{\Delta_n}.
\end{equation}

On the other hand, the Markov properties for statistical inference
of unknown qualities in model (\ref{1.2}) based on the samples
$\{\widetilde{X}_{i\Delta_n}; i = 1, 2, \cdot \cdot \cdot \}$ should
be built, which are infinitesimal conditional expectations
characterized by infinitesimal operators. Fortunately, under Lemma
\ref{l1} we can build the following infinitesimal conditional
expectations for model (\ref{1.2})
\begin{align}
& \label{2.3} E[\frac{\widetilde{X}_{(i+1)\Delta_{n}} -
\widetilde{X}_{i\Delta_{n}}}{\Delta_{n}}|\mathscr{F}_{(i-1)\Delta_{n}}]
 = \mu (X_{(i-1)\Delta_{n}}) + O_{p}(\Delta_{n}),\\
& \label{2.4} E[\frac{(\widetilde{X}_{(i+1)\Delta_{n}} -
\widetilde{X}_{i\Delta_{n}})^{2}}{\Delta_{n}}|\mathscr{F}_{(i-1)\Delta_{n}}]
= \frac{2}{3}\sigma^{2} (X_{(i-1)\Delta_{n}}) +
\frac{2}{3}\int_{\mathscr{E}}c^{2}(X_{(i-1)\Delta_{n}} , z)f(z)dz +
O_{p}(\Delta_{n}).
\end{align}
where $\mathscr{F}_t=\sigma\{X_s,s \leq t\}$. One can refer to
Appendix A in Song, Lin and Wang \cite{slw} for detailed
calculations.

Under $\{\widetilde{X}_{i\Delta_{n}}; i = 0, 1, 2, \cdot \cdot
\cdot\}$, we construct local linear estimators for the unknown
coefficients in model (\ref{1.2}) based on equations (\ref{2.3}) and
(\ref{2.4}). We consider the following weighted local linear
regression to estimate $\mu(x)$ and $M(x)
=\sigma^2(x)+\int_{\mathscr{E}}c(x,z)f(z)dz$ using asymmetric Gamma
kernel, respectively:
\begin{equation}
\label{2.5} \arg \min_{a,
b}\sum_{i=1}^{n}\Big(\frac{\widetilde{X}_{(i+1)\Delta_{n}} -
\widetilde{X}_{i\Delta_{n}}}{\Delta_{n}} - a -
b(\widetilde{X}_{i\Delta_{n}} - x) \Big)^{2}
 K_{G(x/h + 1, h)}\big(\widetilde{X}_{(i-1)\Delta_{n}}\big),
\end{equation}
\begin{equation}
\label{2.6} \arg \min_{a,
b}\sum_{i=1}^{n}\Big(\frac{\frac{3}{2}(\widetilde{X}_{(i+1)\Delta_{n}}
- \widetilde{X}_{i\Delta_{n}})^{2}}{\Delta_{n}} - a -
b(\widetilde{X}_{i\Delta_{n}} - x) \Big)^{2} K_{G(x/h + 1,
h)}\big(\widetilde{X}_{(i-1)\Delta_{n}}\big),
\end{equation} where
$K_{G(x/h + 1, h)}( \cdot )$ is the asymmetric Gamma kernel
function.

The solutions for $a$ to (\ref{2.5}) and (\ref{2.6}) as follows are
respectively the local linear estimators of $\mu(x)$ and $M(x) =
\sigma^2(x)+\int_{\mathscr{E}}c(x,z)f(z)dz,$
\begin{equation}
\label{2.7} \hat{\mu}_n(x) = \frac{\sum_{i=1}^n \omega_{i-1}\Big(
\frac{\widetilde{X}_{i+1}-\widetilde{X}_{i}}
{\Delta_n}\Big)}{\sum_{i=1}^n \omega_{i-1}},
\end{equation}
\begin{equation}
\label{2.8} \hat{M}_{n}(x) = \frac{\sum_{i=1}^n \omega_{i-1}
\frac{3}{2}
\frac{(\widetilde{X}_{i+1}-\widetilde{X}_{i})^2}{\Delta_n}}{\sum_{i=1}^n
\omega_{i-1}}
\end{equation}
where
\begin{eqnarray*}\omega_{i-1} & = & K_{G(x/h + 1, h)}\big(\widetilde{X}_{i-1}\big)(\sum_{j=1}^n K_{G(x/h + 1, h)}\big(\widetilde{X}_{j-1}\big)
(\widetilde{X}_{j} - x)^2\\ & ~ & - (\widetilde{X}_{i} - x)
\sum_{j=1}^{n} K_{G(x/h + 1, h)}\big(\widetilde{X}_{j-1}\big)
(\widetilde{X}_j - x)). \end{eqnarray*}

There are some differences between the estimators given in this
paper and the local linear estimators for the coefficients of model
(\ref{1.1}) in Hanif \cite{hm1}. In (\ref{2.5}) and (\ref{2.6}) the
observations in $(\widetilde{X}_i - x)$ and $K(\widetilde{X}_{i-1})$
are different for the model (\ref{1.2}) which is more complex than
the one in Hanif \cite{hm1}, because we need to calculate some
meaningful conditional expect values of the estimators in the
detailed proof by means of the method introduced by Nicolau
\cite{ni}, but cannot obtain the desired result if they are
identical. Fortunately, $\widetilde{X}_i$ and $\widetilde{X}_{i-1}$
can both approximate the value of $X_{i-1}$, which guarantees the
desired result in the article reasonably from the result in Hanif
\cite{hm1}. However, the local linear smoothing constructed in this
paper cannot be extended to the local polynomial cases, which needs
more than two values to approximate $X_{(i - 1)\Delta_n}.$

We now present some assumptions used in the paper. In what follows,
let $\mathscr{D} = (l , u)$ with $l \geq -\infty$ and $u \leq
\infty$ denote the admissible range of the process $X_{t},$ $K$
denotes $K_{G(x/h + 1, h)}.$

\begin{ass}\label{a1}
(i)~(Local Lipschitz continuity)~ For each $n \in \mathbb{N},$ there
exist a constant $L_{n}$ and a function $\zeta_{n}~:~\mathscr{E}
\rightarrow \mathbb{R}_{+}$ with $\int_{\mathscr{E}}
\zeta_{n}^{2}(z) f(z) dz < \infty$ such that, for any $|x| \leq n,
|y| \leq n,~z \in \mathscr{E}$,
$$|\mu(x) - \mu(y)| + |\sigma(x) - \sigma(y)| \leq L_{n}|x - y|,~~~ |c(x , z) - c(y , z)| \leq \zeta_{n}(z)|x - y|.$$

(ii)~(Linear growthness)~ For each $n \in \mathbb{N}$, there exist
$\zeta_{n}$ as above and $C$, such that for all $x \in \mathbb{R},~z
\in \mathscr{E}$,
$$|\mu(x)| + |\sigma(x)| \leq C (1 + |x|) , ~ |c(x , z)| \leq \zeta_{n}(z)(1 + |x|).$$
\end{ass}

\begin{remark}
Assumption \ref{a1} guarantees the existence and uniqueness of a
solution to $X_{t}$ in Eq. (\ref{1.1}) on the probability space
$(\Omega, \mathscr{F}, P)$, see Jacod and Shiryayev \cite{js}.
\end{remark}

\begin{ass}\label{a2}
The process $X_{t}$ is ergodic and stationary with a finite
invariant measure $\phi(x)$. Furthermore, The process $X_{t}$ is
$\rho-$mixing with $\sum_{i\geq1}\rho(i\Delta_{n}) =
O(\frac{1}{\Delta_{n}^{\alpha}}), ~n \rightarrow \infty,$ where
$\alpha < \frac{1}{2}.$
\end{ass}

\begin{remark}
The hypothesis that $X_{t}$ is a stationary process is obviously a
plausible assumption because for major integrated time series data,
a simple differentiation generally assures stationarity. The same
condition yielding information on the rate of decay of $\rho-$mixing
coefficients for $X_{t}$ was mentioned the Assumption 3 in
Gugushvili and Spereij \cite{gs}.
\end{remark}

\begin{ass}\label{a3}
For any 2 $\leq i$ $\leq n$, $g$ is a differentiable function on
$\mathbb{R}$
 and $\xi_{n,i} = \theta X_{(i-1)\Delta_{n}} + (1 - \theta)
\widetilde{X}_{(i-1)\Delta_{n}}$, $0 \leq \theta \leq 1,$ the
conditions hold:

(i)~~~$\lim_{h \rightarrow 0}E\big[~|h K^{'}(\xi_{n,i})
g(X_{(i-1)\Delta_{n}})|\big] < \infty,$

(ii)~~$\lim_{h \rightarrow 0} h^{1/2} E\big[~|h^{2}
K^{'2}(\xi_{n,i}) g(X_{(i-1)\Delta_{n}})|\big] < \infty$ for
``interior x'',

(iii)~$\lim_{h \rightarrow 0} h E\big[~|h^{2} K^{'2}(\xi_{n,i})
g(X_{(i-1)\Delta_{n}})|\big] < \infty$ for ``boundary x''.
\end{ass}

\begin{remark}
According to the procedure for assumption \ref{a3} in Appendix
(7.1), we can easily deduce the following results:

(i)~~~$\lim_{h \rightarrow 0}E\big[~|K(X_{(i-1)\Delta_{n}})
g(X_{(i-1)\Delta_{n}})|\big] < \infty;$

(ii)~~$\lim_{h \rightarrow 0} h^{1/2} E\big[~|
K^{2}(X_{(i-1)\Delta_{n}}) g(X_{(i-1)\Delta_{n}})|\big] < \infty$
for ``interior x'';

(iii)~$\lim_{h \rightarrow 0} h E\big[~|K^{2}(X_{(i-1)\Delta_{n}})
g(X_{(i-1)\Delta_{n}})|\big] < \infty$ for ``boundary x''.
\end{remark}

\begin{ass}\label{a4}
For all $p \geq 1$, $\sup_{t\geq0} E[|X_{t}|^{p}] < \infty,$ and
$\int_{\mathscr{E}}|z|^{p}f(z)dz < \infty.$
\end{ass}

\begin{remark}
This assumption guarantees that Lemma \ref{l1} can be used properly
throughout the article. If $X_{t}$ is a L\'{e}vy process with
bounded jumps (i.e., $\sup_{t}|\Delta X_{t}| \leq C < \infty$ almost
surely, where C is a nonrandom constant), then $E\{|X_{t}|^{n}\} <
\infty ~\forall n$, that is, $X_{t}$ has bounded moments of all
orders, see Protter \cite{pr}. This condition is widely used in the
estimation of an ergodic diffusion or jump-diffusion from discrete
observations, see Florens-Zmirou \cite{fz}, Kessler \cite{ke},
Shimizu and Yoshida \cite{sy}.
\end{remark}

\begin{ass}\label{a5}
$\Delta_{n} \rightarrow 0,~h \rightarrow 0,~\frac{n
\Delta_{n}}{h}\sqrt{\Delta_{n} \log \Big( \frac{1}{\Delta_{n}}\Big)}
\rightarrow 0,~~h_{n}n\Delta_{n}^{1 + \alpha} \rightarrow \infty,$
as $n \rightarrow \infty.$
\end{ass}

\begin{remark}
The relationship between $h_{n}$ and $\Delta_{n}$ is similar as the
stationary case in Hanif \cite{hm1}, (b1) , (b2) of A8 in Nicolau
\cite{ni} and assumption 7 in Song \cite{syp1}. Wang and Zhou
\cite{wz} presented the optimal bandwidth of symmetric kernel
nonparametric threshold estimator of diffusion function in jump
 - diffusion models. We will select the optimal smoothing parameter
$h_{n}$ for Gamma asymmetric kernel estimation of second-order jump
- diffusion models by means of minimizing the mean square error
(MSE) and $k-$block cross-validation method in Remark \ref{r3.3}.
\end{remark}

\section{ Large sample properties }

\noindent Based on the above assumptions and the lemmas in the
following proof procedure part, we have the following asymptotic
properties. To simplify notations, we define $x \in \mathscr{D}$ to
be a
\begin{equation*}
``interior~x''~~if~~``x/h_{n} \longrightarrow
\infty''~~or~~``boundary~x''~~if~~``x/h_{n} \longrightarrow \kappa''
\end{equation*}

\begin{theorem} \label{thm1}
If Assumptions \ref{a1} - \ref{a5} hold, then
$$\hat{\mu}_n(x) \stackrel{p}{\rightarrow} \mu(x),$$
$$\hat{M}_n(x) \stackrel{p}{\rightarrow} M(x).$$
\end{theorem}

\begin{theorem} \label{thm2}
(i) Under Assumptions \ref{a1} - \ref{a5} and for ``interior x'', if
$h = O((n\Delta_{n})^{-2/5}),$ then
\begin{align*}
& \sqrt{n\Delta_{n} h^{1/2}}\big(\hat{\mu}_{n}(x) - \mu(x) - h
B_{\hat{\mu}_{n}(x)}\big) \stackrel{d}{\rightarrow} N\Big( 0 ,
\frac{M(x)}{2 \sqrt{\pi} x^{1/2} p(x)}\Big),\\
& \sqrt{n\Delta_{n} h^{1/2}}\big(\hat{M}_{n}(x) - M(x) - h
B_{\hat{M}_{n}(x)}\big) \stackrel{d}{\rightarrow} N\Big( 0 ,
\frac{\int_{\mathscr{E}}c^{4}(x , z)f(z)dz}{2 \sqrt{\pi} x^{1/2}
p(x)}\Big),
\end{align*}
where $B_{\hat{\mu}_{n}(x)} =
\frac{x}{2}\mu^{''}(x),~~B_{\hat{M}_{n}(x)} = \frac{x}{2}M^{''}(x),$

(ii) Under Assumptions \ref{a1} - \ref{a5} and for ``boundary x'',
if $h = O((n\Delta_{n})^{-1/5}),$ then
\begin{align*}
& \sqrt{n\Delta_{n} h}\big(\hat{\mu}_{n}(x) - \mu(x) - h^{2}
B^{'}_{\hat{\mu}_{n}(x)}\big) \stackrel{d}{\rightarrow} N\Big( 0 ,
\frac{M(x) \Gamma(2\kappa + 1)}{2^{2\kappa + 1} \Gamma^{2}(\kappa + 1) p(x)}\Big),\\
& \sqrt{n\Delta_{n} h}\big(\hat{M}_{n}(x) - M(x) - h^{2}
B^{'}_{\hat{M}_{n}(x)}\big) \stackrel{d}{\rightarrow} N\Big( 0 ,
\frac{\int_{\mathscr{E}}c^{4}(x , z)f(z)dz \Gamma(2\kappa +
1)}{2^{2\kappa + 1} \Gamma^{2}(\kappa + 1) p(x)}\Big),
\end{align*}
where $B^{'}_{\hat{\mu}_{n}(x)} = \frac{1}{2}(2 +
\kappa)\mu^{''}(x),~~B^{'}_{\hat{M}_{n}(x)} = \frac{1}{2}(2 +
\kappa)M^{''}(x).$
\end{theorem}

\begin{remark}
\label{r3.2} In this article, we considered the asymptotic
consistency and weak convergence of local linear estimators for the
unknown quantities in the second-order jump-diffusion model which
can be directly used to model the returns of a asset, based on Gamma
asymmetric kernel in Theorem \ref{thm1} and Theorem \ref{thm2}. The
main method to obtain the asymptotic properties for the estimators
of model (\ref{1.2}) is to approximate the estimator for model
(\ref{1.2}) by the similar estimator for model (\ref{1.1}) in
probability. One can refer to Nicolau \cite{ni} for the same idea.
Fortunately, lemma \ref{l2} in the proofs section builds the bridge,
which provides us the desired properties of local linear estimator
based on Gamma asymmetric kernel for model (\ref{1.1}). We only
discussed the stationary jump-diffusion in lemma 2. Actually, the
similarly theoretical and numerical results as that in Theorem
\ref{thm1} and Theorem \ref{thm2} also hold for the univariate case:
the jump-diffusion model which can be employed to model the price of
a asset, whether it is stationary or not. With the similar procedure
as Bandi and Nguyen \cite{bn}, lemma \ref{l2} and some conditions on
the local time in Wang and Zhou \cite{wz}, one can easily deduce
their asymptotic consistency and normality of the local linear
estimators for the unknown quantities in the univariate
nonstationary jump-diffusion model based on Gamma asymmetric
kernels. It is not our objective in this paper and thus it is less
of a concern here. We will take it into consideration in the future
work.
\end{remark}

\begin{remark}
\label{r3.3} Theorems \ref{thm1} and \ref{thm2} give the weak
 consistency and the asymptotic normality of the local linear
 estimators using Gamma asymmetric kernels. As discussed in Chapman and Pearson \cite{cp},
 the performance of nonparametric kernel estimator depends crucially on the choice of the smoothing parameter $h_{n}.$ Hence, it is very important to
 consider the choice of the smoothing parameter $h_{n}$
for the nonparametric estimation using asymmetric kernels. Here we
will select the optimal smoothing parameter $h_{n}$ based on the
mean square error (MSE). Take $\mu(x)$ for example,

\noindent for ``interior x'',  the optimal smoothing parameter
$h_{n}$ is
$$h_{n, opt} = \left(\frac{1}{n \Delta_{n}} \cdot \frac{M(x)}{2 \sqrt{\pi} x^{1/2} p(x)}
\cdot \frac{4}{[x \mu^{''}(x)]^{2}}\right)^{\frac{2}{5}} =
O_{p}\left(\frac{1}{n \Delta_{n}}\right)^{\frac{2}{5}}$$ and the
corresponding MSE is $O_{p}\left(\frac{1}{n
\Delta_{n}}\right)^{\frac{4}{5}}.$

\noindent For ``boundary x'', the optimal smoothing parameter
$h_{n}$ is
$$h_{n, opt} = \left(\frac{1}{n \Delta_{n}} \cdot \frac{M(x) \Gamma(2\kappa + 1)}{2^{2\kappa + 1}
\Gamma^{2}(\kappa + 1) p(x)} \cdot \frac{4}{[(2 + \kappa)
\mu^{''}(x)]^{2}}\right)^{\frac{1}{5}} = O_{p}\left(\frac{1}{n
\Delta_{n}}\right)^{\frac{1}{5}}$$ and the corresponding MSE is also
$O_{p}\left(\frac{1}{n \Delta_{n}}\right)^{\frac{4}{5}}.$ When $h
\leq h_{n, opt},$ the bias term contributes to larger part of the
mean square error, while the coverage rate contributes to larger
part if $h > h_{n, opt}.$ One can observe that $h_{n, opt}$ for
``boundary x'' is larger than that for ``interior x'' as $n
\rightarrow \infty$ because more sample points are required for
boundary bias reduction.

In practice, we can take the plug-in method studied in Fan and
Gijbels \cite{fg3} to obtain an optimal smoothing bandwidth $h_{n}$
on behalf of MSE. As mentioned in Gospodinov and Hirukawa \cite{gh},
the bandwidth $h_{n}$ constructed above relies on the consistent
estimators for these unknown quantities and they are difficult to
obtain and may give rise to bias. Moreover, Hagmann and Scaillet
\cite{hs2}, regarding global properties, discussed the choice of
bandwidth to ameliorate the adaptability of Gamma asymmetric kernel
estimators since the bandwidth $h_{n}$ constructed above varies with
the change of the design point $x.$ Hence we mention two rules of
thumb on selecting a global smoothing bandwidth here. For
simplicity, we can use bandwidth selector $h_{n} = c \cdot \hat{S}
\cdot T^{-\frac{2}{5}}$ in  Xu and Phillips \cite{xp}, where
$\hat{S}$, $T$ denote the standard deviation of the data and the
time span and $c$ represents different constants for different
estimators to be estimated by means of minimizing the Mean Square
Errors (MSE). Or, we can employ the $k-$block cross-validation
method proposed by Racine \cite{rac} to assess the performance of an
estimator via estimating its prediction error. The main idea is to
minimize the following expression: $CV(h_{n}) =
n^{-1}\sum_{i=k+1}^{n-k}\{\frac{\widetilde{X}_{(i+1)\Delta_{n}} -
\widetilde{X}_{i\Delta_{n}}}{\Delta_{n}} -
\hat{\mu}_{h_{n},-(i-k):-(i+k)}(\widetilde{X}_{i\Delta_{n}})\}^{2}$
with $k = n^{1/4}$ to eliminate the dependence of data, where
$\hat{\mu}_{h_{n},-(i-k):-(i+k)}(\widetilde{X}_{i\Delta_{n}})$ is
the underlying estimator (\ref{2.7}) as a function of bandwidth
$h_{n}$, but without using the $i - k$th to $i + k$th observations.
In practice, the optimal data-dependent choice of block size $k$
should take into consideration the persistence of the empirical
data, which one can refer to in the selection of smoothing parameter
part in Gospodinov and Hirukawa \cite{gh}. These two rules of thumb
on selecting the bandwidth are displayed numerically in simulation
and empirical analysis part. For the further study of the optimal
value of the bandwidth, one can refer to A\"it-Sahalia and Park
\cite{ap}.
\end{remark}

\begin{remark}
\label{r3.4} In addition, the asymptotic normality of local linear
estimator using Gamma asymmetric kernel for $\mu(x)$ in this paper
is different from that in Chen and Zhang \cite{cz}, where their
asymptotic normality was
$$\sqrt{h_{n} n \Delta_n }(\hat{\mu}_{n} (x) - \mu(x) - h_{n}^{2} \frac{1}{2}\mu^{''}(x))
\stackrel{d}{\rightarrow} N(0,V\frac{M(x)}{p(x)})$$ with $V =
\frac{(K_1^2)^2 K_2^0+ (K_1^1)2 K_2^2
-2(K_1^1)(K_1^2)K_2^1}{[K_1^2-(K_1^1)^2]^2}$ which is equal to
$\frac{1}{2 \sqrt{\pi}}$ if the kernel is Gaussian kernel. There are
two main differences: on one hand, the convergence rate of local
linear estimator using Gamma asymmetric kernel is different for the
location of the design point $x$ such as ``interior $x$'' and
``boundary $x$''; on the other hand, the variance of of local linear
estimator using Gamma asymmetric kernel is inversely proportional to
the design $x,$ which shows that the variance decreases as the
design point $x$ increases. Additionally, the optimal bandwidth
$h_{n, opt}^{Gaussian}$ is $O(\frac{1}{n \Delta_{n}})^{\frac{1}{5}}$
for any design point $x$ and the corresponding MSE is $O(\frac{1}{n
\Delta_{n}})^{\frac{4}{5}}.$ For ``interior $x$'', the optimal
smoothing parameter $h_{n, opt}^{Gamma} = O(\frac{1}{n
\Delta_{n}})^{\frac{2}{5}} = O(h_{n, opt}^{Gaussian})^{2},$ which
means that the asymptotic variance of the estimator constructed with
Gamma asymmetric kernel is $O(\frac{1}{h_{n, opt}^{Gaussian} n
\Delta_n})$ the same as that constructed with Gaussian symmetric
kernel.
\end{remark}

\begin{remark}
\label{r3.5}  For ``interior x'', if the smoothing parameter $h_{n}
= O((n\Delta_{n})^{-2/5}),$ the normal confidence interval for
$\mu(x)$ using Gamma asymmetric kernel and Gaussian symmetric kernel
at the significance level $100(1-\alpha)\%$ are constructed as
follows,
\begin{align*}
I_{\mu,\alpha}^{Asym} = & \Bigg[\hat{\mu}_{n}^{Asym}(x) - h_{n}
\cdot \frac{x}{2} \hat{\mu}_{n}^{''}(x) - z_{1-\alpha/2} \cdot
\frac{1}{\sqrt{n\Delta_{n} h_{n}^{1/2}}} \cdot
\sqrt{\frac{\hat{M}_{n}^{Asym}(x)}{2 \sqrt{\pi} x^{1/2} \hat{p}_{n}^{Asym}(x)}},\\
& \hat{\mu}_{n}^{Asym}(x) - h_{n} \cdot \frac{x}{2}
\hat{\mu}_{n}^{''}(x)
+ z_{1-\alpha/2} \cdot \frac{1}{\sqrt{n\Delta_{n} h_{n}^{1/2}}} \cdot
\sqrt{\frac{\hat{M}_{n}^{Asym}(x)}{2 \sqrt{\pi} x^{1/2} \hat{p}_{n}^{Asym}(x)}} \Bigg],\\
I_{\mu,\alpha}^{Sym} = & \Bigg[\hat{\mu}_{n}^{Sym}(x) - h_{n}^{2}
\cdot \frac{1}{2} \hat{\mu}_{n}^{''}(x) - z_{1-\alpha/2} \cdot
\frac{1}{\sqrt{n\Delta_{n} h_{n}}} \cdot
\sqrt{\frac{\hat{M}_{n}^{Sym}(x)}{2 \sqrt{\pi} \hat{p}_{n}^{Sym}(x)}},\\
& \hat{\mu}_{n}^{Sym}(x) - h_{n}^{2} \cdot \frac{1}{2}
\hat{\mu}_{n}^{''}(x) + z_{1-\alpha/2} \cdot
\frac{1}{\sqrt{n\Delta_{n} h_{n}}} \cdot
\sqrt{\frac{\hat{M}_{n}^{Sym}(x)}{2 \sqrt{\pi}
\hat{p}_{n}^{Sym}(x)}} \Bigg].
\end{align*}
$\hat{\mu}_{n}^{Asym}(x)$,
 $\hat{M}_{n}^{Asym}(x)$ , $\hat{\mu}_{n}^{Sym}(x)$,
 $\hat{M}_{n}^{Sym}(x)$ denote the local
linear estimators of $\mu(x), M(x)$ in (\ref{2.7}) and (\ref{2.8})
using Gamma asymmetric kernel or Gaussian symmetric kernel,
respectively. $z_{1-\alpha/2}$ is the inverse CDF for the standard
normal distribution evaluated at $1 - \alpha/2.$
$\hat{p}_{n}^{Asym}(x) = \frac{1}{n} \sum_{i=1}^{n} K_{Gamma(x/h_{n}
+ 1,
h)}\big(\widetilde{X}_{(i-1)\Delta_{n}}\big),~\hat{p}_{n}^{Sym}(x) =
\frac{1}{nh_{n}} \sum_{i=1}^{n} K_{Gaussian} \Big(\frac{x -
\widetilde{X}_{(i-1)\Delta_{n}}}{h_{n}}\Big).$ As Fan and Gijbels
\cite{fg2} showed, the derivative $\hat{\mu}_{n}^{''}(x)$ in
$I_{\mu,\alpha}^{Asym}$ can be estimated by taking the second
derivative of the local linear estimators of $\mu(x)$ in (\ref{2.7})
using Gamma asymmetric kernel. In the similar manner, one can give
the normal confidence intervals for $\mu(x)$ of ``boundary x'' using
Gamma asymmetric kernel or Gaussian symmetric kernel. As Xu
\cite{xkl} considered, the resultant normal confidence interval
$I_{\mu,\alpha}^{Asym}$ or $I_{\mu,\alpha}^{Sym}$ has more correct
coverage rate asymptotically as long as a smoothing parameter
$h_{n}$ is used and the bias and variance can be consistently
estimated.

Similarly, the normal confidence interval for $M(x)$ at a spatial
point $x$ using Gamma asymmetric kernel and Gaussian symmetric
kernel at the significance level $100(1-\alpha)\%$ can be
constructed. The final interval for $M(x)$ should be taken as the
intersection of $I_{M,\alpha}^{Asym}$ or $I_{M,\alpha}^{Sym}$ with
$[0, +\infty)$ to coincide with the nonnegativity of the conditional
variance $M(x).$
\end{remark}

\begin{remark}
\label{r3.6} Here we briefly commented the theoretical comparison
for lengths of confidence intervals for ``interior x'' and
``boundary x'' using Gamma asymmetric kernel or Gaussian symmetric
kernel. One can refer to the simulation part for the numerical
comparison for lengths of confidence intervals.

Take $\mu(x)$ for example. The dominant factors that affect the
length of the confidence interval are the various coverage rate and
the different coefficient in the variance.

For ``interior $x$'', the coverage rate of the local linear
estimator based on Gamma asymmetric kernel is
$\frac{1}{\sqrt{n\Delta_{n} h_{n}^{1/2}}}$, which is much smaller
than $\frac{1}{\sqrt{n\Delta_{n} h_{n}}}$ of that based on Gaussian
symmetric kernel with a given $h_{n}$. Compared with the ones using
Gaussian symmetric kernel, the variance of local linear estimator
using Gamma asymmetric kernel is inversely proportional to the
design $x$, which shows that the length of the confidence interval
decreases as the design point $x$ increases.

For ``boundary $x$'', although the coverage rate of the local linear
estimator based on Gamma asymmetric kernel is the same as that based
on Gaussian symmetric kernel, the coefficient in their variance
differs a little such as $\frac{\Gamma(2\kappa + 1)}{2^{2\kappa + 1}
\Gamma^{2}(\kappa + 1)}$ for Gamma asymmetric kernel while
$\frac{1}{2 \sqrt{\pi}}$ for Gaussian symmetric kernel. Under
numeral calculations, from Table 1 we can conclude that when $\kappa
\leq 0.7$, the variance based on Gamma asymmetric kernel is larger
while when $\kappa \geq 0.75$, the variance based on Gamma
asymmetric kernel is smaller than that based on Gaussian symmetric
kernel. This reveals that the closer to the boundary point or the
larger bandwidth for fixed ``boundary $x$'', the shorter the length
of confidence interval based on Gaussian symmetric kernel, which is
shown in the simulation result.

\begin{table}[!htb]
  \centering
  \caption{The difference between $\frac{\Gamma(2\kappa + 1)}{2^{2\kappa + 1} \Gamma^{2}(\kappa + 1)}$ and $\frac{1}{2 \sqrt{\pi}}$ in the variance for various $\kappa$.}
  \begin{tabular}{ccccccccc}
    \toprule
    Value of $\kappa$ & 0.25  & 0.3   & 0.35  & 0.4   & 0.45  & 0.5   & 0.55  & 0.6 \\
    \midrule
    Difference & 0.0993 & 0.0838 & 0.0701 & 0.0577 & 0.0464 & 0.0362 & 0.0269 & 0.0183 \\
    \midrule
    Value of $\kappa$ & 0.65  & 0.7   & 0.75  & 0.8   & 0.85  & 0.9   & 0.95  & 1 \\
    \midrule
    Difference & 0.0103 & 0.003 & -0.004 & -0.01 & -0.016 & -0.022 & -0.027 & -0.032 \\
    \midrule
    Value of $\kappa$ & 1.25  & 1.5   & 1.75  & 2     & 2.25  & 2.5   & 2.75  & 3 \\
    \midrule
    Difference & -0.053 & -0.07 & -0.083 & -0.095 & -0.104 & -0.112 & -0.12 & -0.126 \\
    \midrule
    Value of $\kappa$ & 3.25  & 3.5   & 3.75  & 4     & 4.25  & 4.5   & 4.75  & 5 \\
    \midrule
    Difference & -0.132 & -0.137 & -0.141 & -0.145 & -0.149 & -0.153 & -0.156 & -0.159 \\
    \bottomrule
    \end{tabular}%
  \label{tab:addlabel}%
\end{table}
\end{remark}

\begin{remark}
\label{r3.7} In contrary to the second-order diffusion model without
jumps (Nicolau \cite{ni}), the second infinitesimal moment estimator
for second-order diffusion model with jumps has a rate of
convergence that is the same as the rate of convergence of the first
infinitesimal moment estimator. Apparently, this is due to the
presence of discontinuous breaks that have an equal impact on all
the functional estimates. As Johannes \cite{jh} pointed out, for the
conditional variance of interest rate changes, not only diffusion
play a certain role, but also jumps account for more than half at
lower interest level rates, almost two-thirds at higher interest
level rates, which dominate the conditional volatility of interest
rate changes. Thus, it is extremely important to estimate the
conditional variance as $\sigma^{2}(x)$ + $\int_{\mathscr{E}}c^{2}(x
, z)f(z)dz$ which reflects the fluctuation of the underlying asset
or the return of the underlying asset.

Meanwhile, for the special case of model (\ref{1.1}) with compound
Poisson jump components, there are several methodologies for the
nonparametric estimation to identify the diffusion coefficient
$\sigma^{2}(x)$, the jump intensity $\lambda(x)$ and the variance of
the jump sizes $\sigma^{2}_{z}$. For single-factor model, one can
use the fourth and sixth moments to identify the jump components
$\lambda(x)$ and $\sigma^{2}_{z}$, then $\sigma^{2}(x)$ can be
identified through second moment like Theorem 3.3, see Johannes
\cite{jh} as following (\ref{3.7}) $\sim$ (\ref{3.9}), or one can
use the threshold estimation for $\sigma^{2}(x)$ , $\lambda(x)$ and
$\sigma^{2}_{z}$, see Mancini and Ren\`{o} (\cite{mr}, Theorems 3.2,
3.7).

Nonparametric estimation to identify the diffusion coefficient
$\sigma^{2}(x)$, the jump intensity $\lambda(x)$ and the variance of
the jump sizes $\sigma^{2}_{z}$ for model (\ref{1.2}) is not our
objective in this paper and thus it is less of a concern here.
However, we give some procedures here to deal with this
identification for the special case of model (\ref{1.2}) similarly
as Johannes \cite{jh} :
\begin{equation}
\label{3.6}
\left\{
\begin{array}{ll}
dY_{t} = X_{t-}dt,\\
dX_{t} = \mu(X_{t-})dt + \sigma(X_{t-})dW_{t} +
d\Big(\sum_{n=1}^{N_{t}}Z_{t_{n}}\Big),
\end{array}
\right. \end{equation} where $N_{t}$ is a doubly stochastic point
process with jump intensity $\lambda(X_{t-})$, and $Z_{t_{n}} \sim
\mathscr{N}(\mu_{z},\sigma_{z})$ with the variance of the jump size
$\sigma^{2}_{z}$. We have the following infinitesimal conditions
(\ref{3.7}) $\sim$ (\ref{3.9}) under simple but tedious calculations
by virtue of Lemma \ref{l1} with $d = 2$ :
\begin{align}
& \label{3.7}
E\Big[\frac{\frac{3}{2}(\widetilde{X}_{(i+1)\Delta_{n}} -
\widetilde{X}_{i\Delta_{n}})^{2}}{\Delta_{n}}|\mathscr{F}_{(i-1)\Delta_{n}}\Big]
= \sigma^{2} (X_{(i-1)\Delta_{n}}) +
\lambda(X_{(i-1)\Delta_{n}})\sigma_{z}^{2} + O_{p}(\Delta_{n}),  \\
& \label{3.8} E\Big[\frac{3(\widetilde{X}_{(i+1)\Delta_{n}} -
\widetilde{X}_{i\Delta_{n}})^{4}}{\Delta_{n}}|\mathscr{F}_{(i-1)\Delta_{n}}\Big]
= 3\lambda(X_{(i-1)\Delta_{n}})(\sigma_{z}^{2})^{2} +
O_{p}(\Delta_{n}), \\
& \label{3.9} E\Big[\frac{3(\widetilde{X}_{(i+1)\Delta_{n}} -
\widetilde{X}_{i\Delta_{n}})^{6}}{\Delta_{n}}|\mathscr{F}_{(i-1)\Delta_{n}}\Big]
= 15\lambda(X_{(i-1)\Delta_{n}})(\sigma_{z}^{2})^{3} +
O_{p}(\Delta_{n}).
\end{align}
Based on (\ref{3.8}) and (\ref{3.9}), we can deduce the $\lambda(x)$
and $\sigma^{2}_{z}$ using kernel estimations, finally the diffusion
coefficient $\sigma^{2}(x)$ can be easily derived from three other
kernel estimations for the conditional variance, jump intensity and
the jump variance. We will take their asymptotic weak consistency
and normality into consideration in the future work based on the
result of Bandi and Nguyen \cite{bn}.
\end{remark}

\section{Monte Carlo Simulation Study}

In this section, we conduct a simple Monte Carlo simulation
experiment aimed at the finite sample performance of the local
linear estimators for both drift and conditional variance functions
constructed with Gamma asymmetric kernels and those constructed with
Gaussian symmetric kernels. Assessment will be made between them by
comparing their mean square error (MSE), coverage rate and length of
confidence band. Our experiment is based on the following data
generating process:
\begin{equation}
\label{mc1} \left\{
\begin{array}{ll}
dY_{t} = X_{t-}dt,\\
dX_{t} = (1-10X_{t-})dt + \sqrt{0.1+0.1X_{t-}^{2}}dW_{t} + dJ_{t},
\end{array}
\right. \end{equation} where the coefficients of continuous part
similar as the ones used in Nicolau (\cite{ni}) (here we add a
constant of 1 in the drift coefficient , which guarantees the
nonnegativity of $X_{t}$ with the initial value of 0.1) and $J_{t}$
is a compound Poisson jump process, that is, $J_{t} =
\sum_{n=1}^{N_{t}}Z_{t_{n}}$ with arrival intensity $\lambda \cdot T
= 20$ or $\lambda \cdot T = 50$ and jump size $Z_{n} \sim
\mathscr{N}(0,0.036^{2})$ or $Z_{n} \sim \mathscr{N}(0,0.1^{2})$
corresponding to Bandi and Nguyen (\cite{bn}), where $t_{n}$ is the
$n$th jump of the Poisson process $N_{t}.$ The parameters $\lambda,
\sigma^{2}$ for $Z_{n}$ are of particular importance, which
represent the intensity and amplitude for the jump component,
respectively. For the time series generated, the two values of
$\lambda$ chosen indicate the moderate and intense rate of
recurrence for jumps, and the two values of $\sigma^{2}$ for $Z_{n}$
chosen mean the slow and high level of jumps.

By taking the integral from 0 to $t$ in the second expression of
(\ref{mc1}), we obtain
\begin{equation}
\label{mc2} X_{t} = 0.1 + t - 10\int_{0}^{t}{X_{s-}}ds +
\int_{0}^{t}{\sqrt{0.1+0.1X_{s-}^{2}}dW_{t}} +
\sum_{n=1}^{N_{t}}Z_{t_{n}}.
\end{equation}
Then we have
\begin{equation}
\label{mc3} Y_{t} = \int_{0}^{t}{X_{s-}}ds = -\frac{1}{10}\Big(X_{t}
- 0.1 - t - \int_{0}^{t}{\sqrt{0.1+0.1X_{s-}^{2}}dW_{t}} -
\sum_{n=1}^{N_{t}}Z_{t_{n}}\Big).
\end{equation}

\noindent $X_{t}$ can be sampled by the Euler-Maruyama scheme
according to (\ref{mc2}), which will be detailed in the following
Algorithm 1 (one can refer to Cont and Tankov \cite{ct}).

\begin{algorithm}[!htb]
\label{alg1} \caption{ Simulation for trajectories of second-order
jump-diffusion model}
\begin{algorithmic}[]
\REQUIRE ~~\\

{\textbf{Step 1:}} generate a standard normal random variate $V$ and
transform it into $D_{i} =
\sqrt{0.1+0.1X_{t_{i-1}}^{2}}*\sqrt{\Delta t_{i}}*V$, where $\Delta
t_{i} = t_{i} - t_{i-1} = \frac{T}{n}$ is the observation time
frequency;

{\textbf{Step 2:}} generate a Poisson random variate $N$ with
intensity $\lambda = 2$;

{\textbf{Step 3:}} generate $N$ random variables $\tau_{i}$
uniformly distributed in $[0, T]$, which correspond the jump times;

{\textbf{Step 4:}} generate $N$ random variables $Z_{\tau_{i}} \sim
\mathscr{N}(0,0.036^{2})$, which correspond the jump sizes;

\noindent One trajectory for $X_{t}$ is
$$X_{t_{i}} = X_{t_{i-1}}+ (1 - 10X_{t_{i-1}})*\Delta t_{i} + D_{i} +
1_{\{t_{i-1} \leq \tau_{i} < t_{i}\}}*Z_{\tau_{i}}.$$

{\textbf{Step 5:}} By substitution of $X_{t_{i}}$ in (\ref{mc3}),
$Y_{t_{i}}$ can be sampled.

\end{algorithmic}
\end{algorithm}

One sample trajectory of the differentiated process $X_{t}$ and
integrated process $Y_{t}$ with $T = 10,~n = 1000,$ $X_{0} = 0.1$
and $Y_{0} = 100$ using Algorithm 1 is shown in FIG 2. Through
observation on FIG 2(b), we can find the following  features of the
integrated process $Y_{t}$: absent mean-reversion, persistent
shocks, time-dependent mean and variance, nonnormality, etc.
\begin{figure}[!htb]
\centering
  \subfigure[The differentiated process $X_{t}$]{
  \label{fig:subfig} 
    \includegraphics[width=0.479\textwidth]{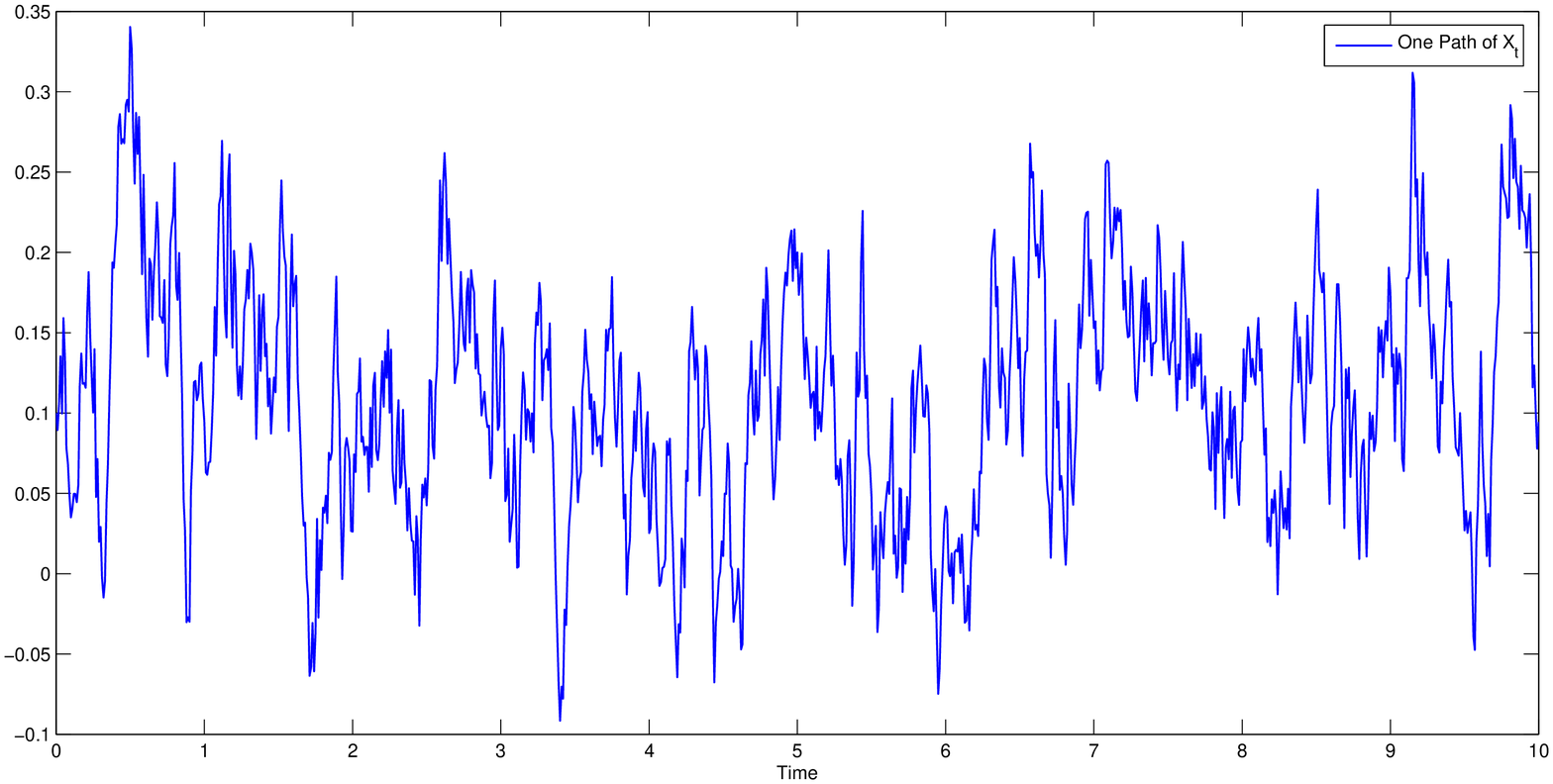}}
 \hspace{0.01in}
  \subfigure[The integrated process $Y_{t}$]{
  \label{fig:subfig} 
    \includegraphics[width=0.479\textwidth]{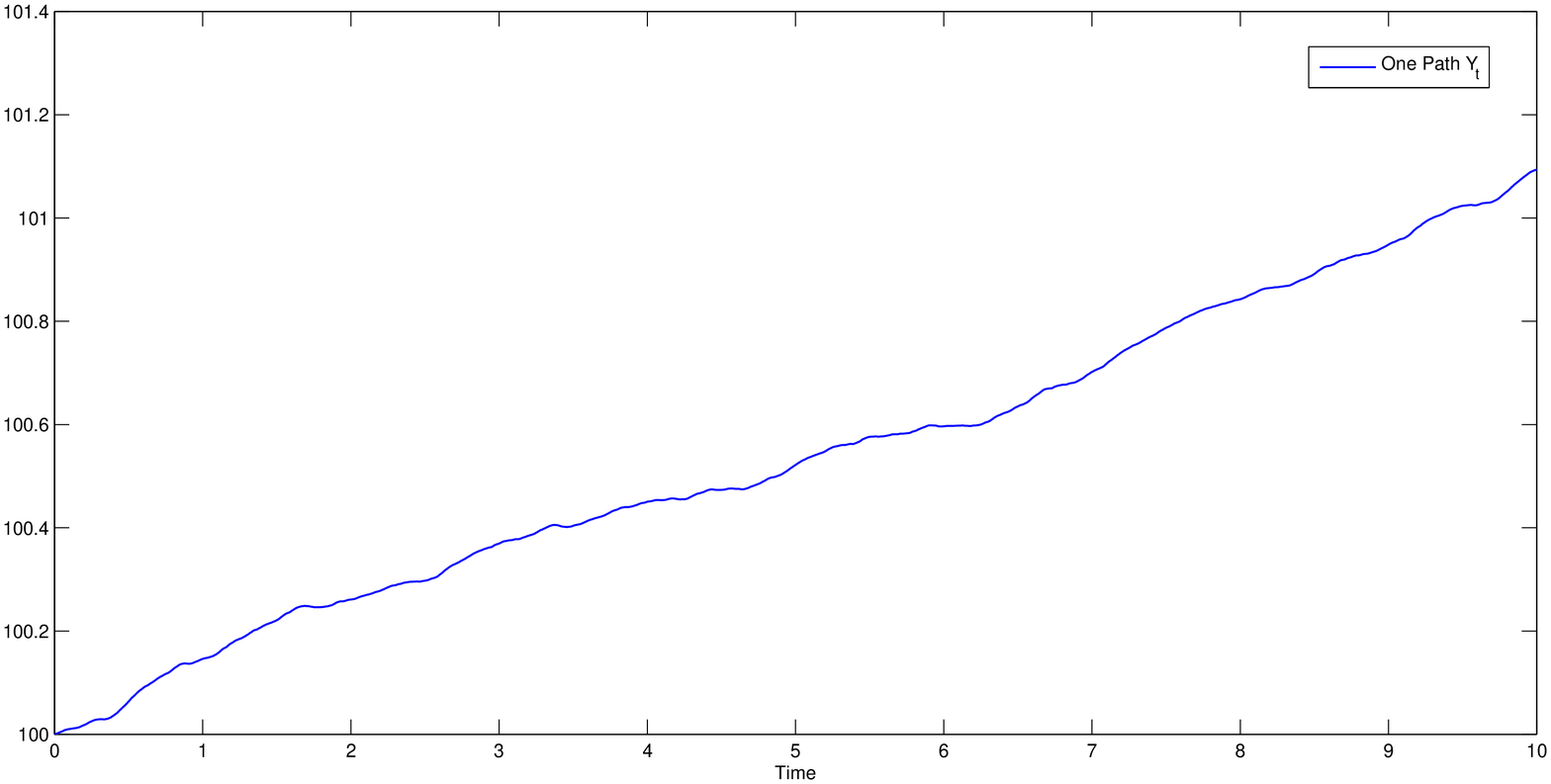}}
  \caption{ the Sample Paths of $X_{t}$ and $Y_{t}$ }
  \label{fig:subfig} 
\end{figure}

Throughout this section, we employ Gamma kernel $K_{G(x/h + 1, h)}
(u) = \frac{u^{x/h} \exp (-u/h)}{h^{x/h + 1} \Gamma(x/h + 1)}$ and
Gaussian kernel $K(x) = \frac{1}{\sqrt{2\pi}}e^{-\frac{x^{2}}{2}}.$
In general, the nonparametric estimation is sensitive to bandwidth
choices, hence we select the data-driven bandwidth via the bandwidth
selector $h = c\hat{S} (n\Delta_{n})^{-\frac{2}{5}} = c\hat{S}
T^{-\frac{2}{5}}$ for ``interior x'' or $h = c\hat{S}
(n\Delta_{n})^{-\frac{1}{5}} = c\hat{S} T^{-\frac{1}{5}}$ for
``boundary x'' similarly as Xu and Phillips \cite{xp}, where
$\hat{S}$ denotes the standard deviation of the data and c
represents different constants for different cases, or
cross-validation method in Racine \cite{rac} according to Remark
\ref{r3.3}. We will compare their mean square error (MSE) under
various lengths of observation time interval $T~(= 50, 100, 500)$
and sample sizes $n~(= 500, 1000, 5000)$ with $\Delta_{n} =
\frac{T}{n}.$ For comparison of coverage rate and length of
confidence band, we will consider three design points such as the
boundary point $x = 0.05$ and interior point $x = 0.15, 0.30,$ which
fall in the range of the simulated data. Meanwhile, we will also
consider five fixed bandwidth $(h_{1}, h_{2}, h_{3}, h_{4}, h_{opt})
= (0.01, 0.02, 0.03, 0.05, 0.0441)$ for estimation of $\mu(x)$ and
$(h_{1}, h_{2}, h_{3}, h_{4}, h_{opt}) = (0.01, 0.02, 0.04, 0.05,
0.0294)$ for estimation of $M(x)$ as that in Xu \cite{xkl}, which
cover the bandwidths used in practice. In this section, the normal
confidence level is assumed to be $95\%.$

\begin{figure}[!htb]
\centering
  \subfigure[ Curve of $MSE(h)$ versus $C$ in $h$ with $C_{opt} = 2.8$ ]{
  \label{fig:subfig} 
    \includegraphics[width=0.479\textwidth]{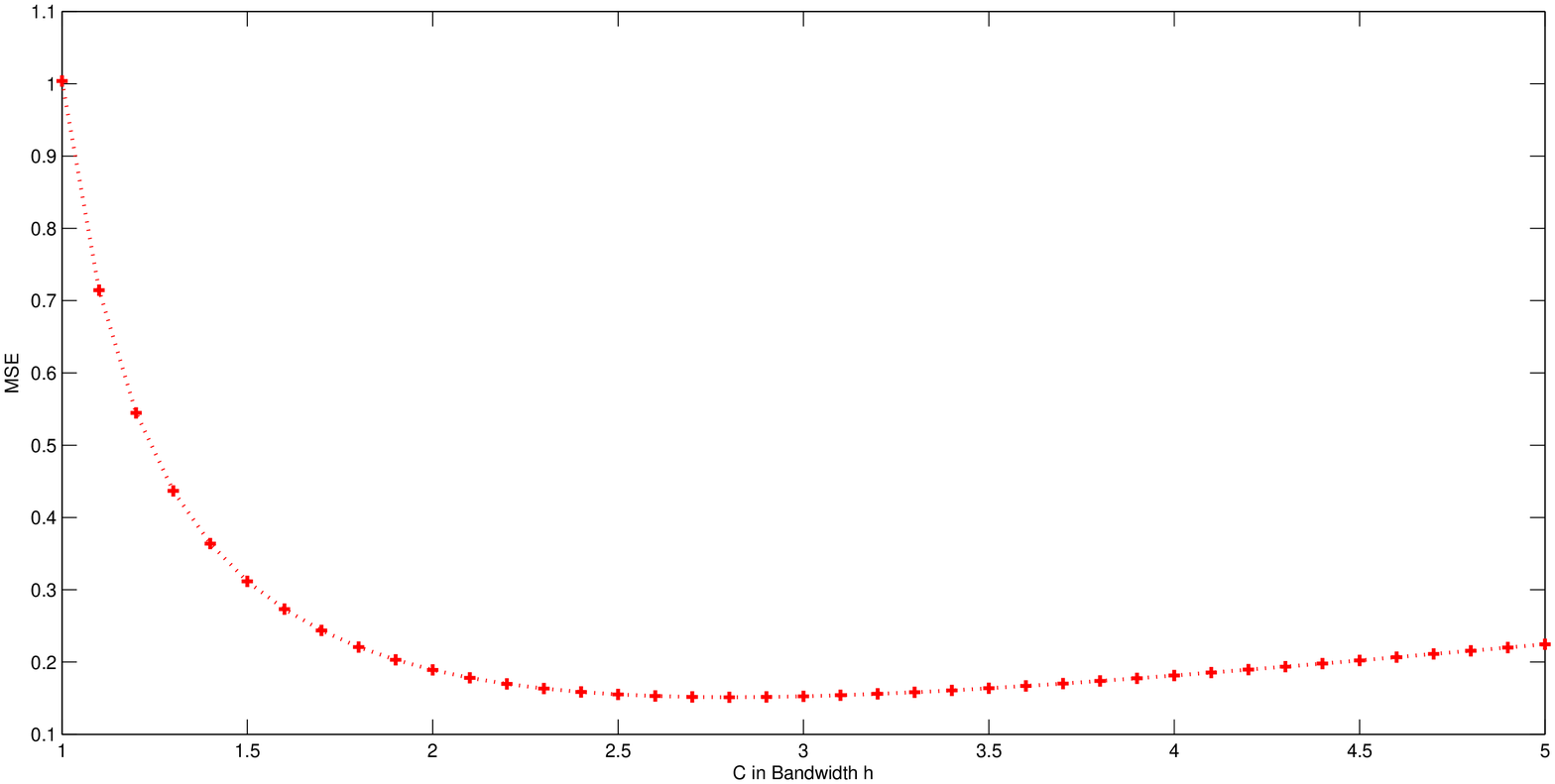}}
 \hspace{0.01in}
  \subfigure[ Curve of $CV(h)$ versus $h$ with $h_{cv} = 0.035$ ]{
  \label{fig:subfig} 
    \includegraphics[width=0.479\textwidth]{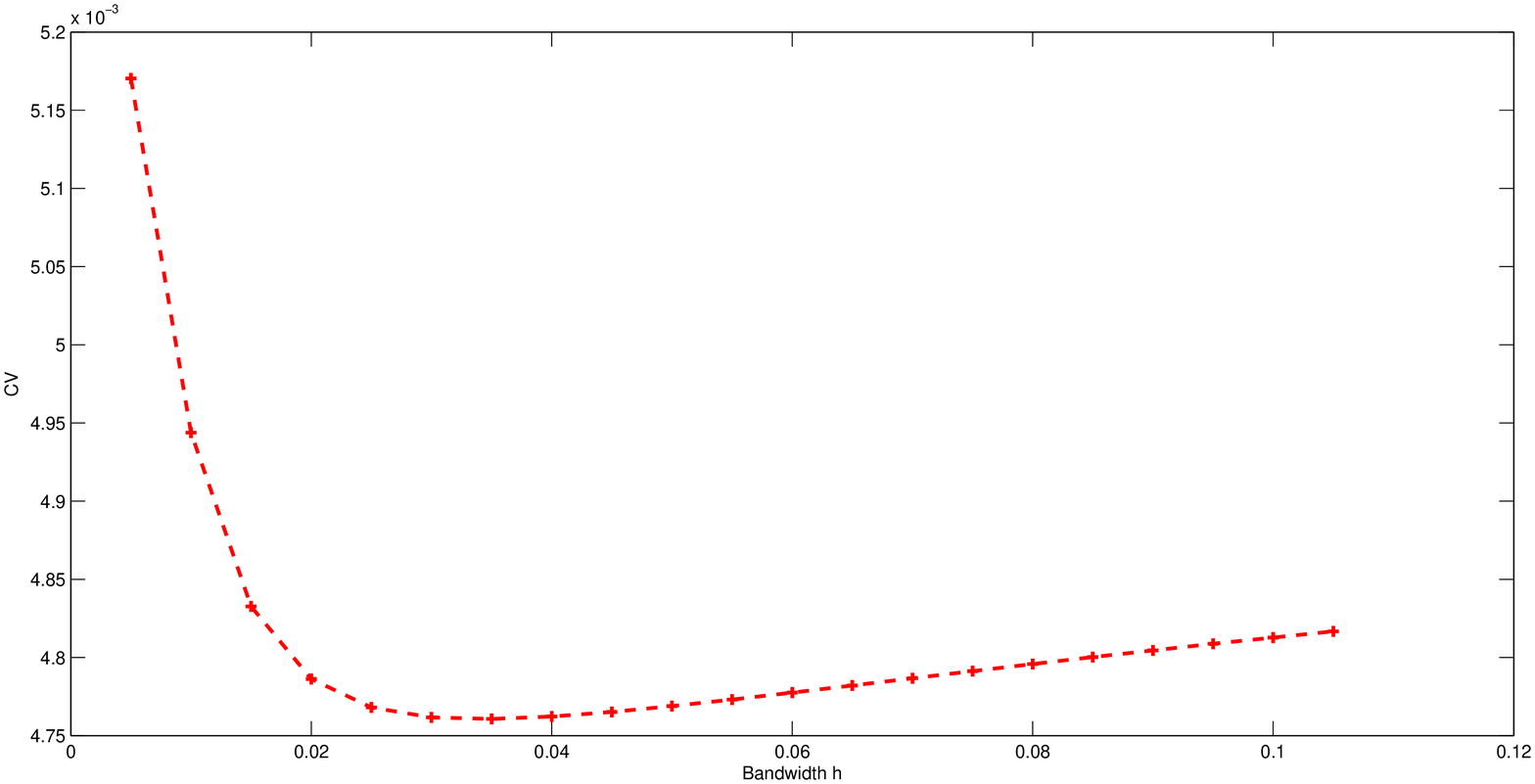}}
  \caption{ Curve for two rules of thumb versus $h$ with $T = 10,~n = 1000$ of $\mu(x) = 1 - 10*x$}
  \label{fig:subfig} 
\end{figure}

\begin{figure}[!htb]
\centering
  \subfigure[ Curve of $MSE(h)$ versus $C$ in $h$ with $C_{opt} = 1.7$  ]{
  \label{fig:subfig} 
    \includegraphics[width=0.479\textwidth]{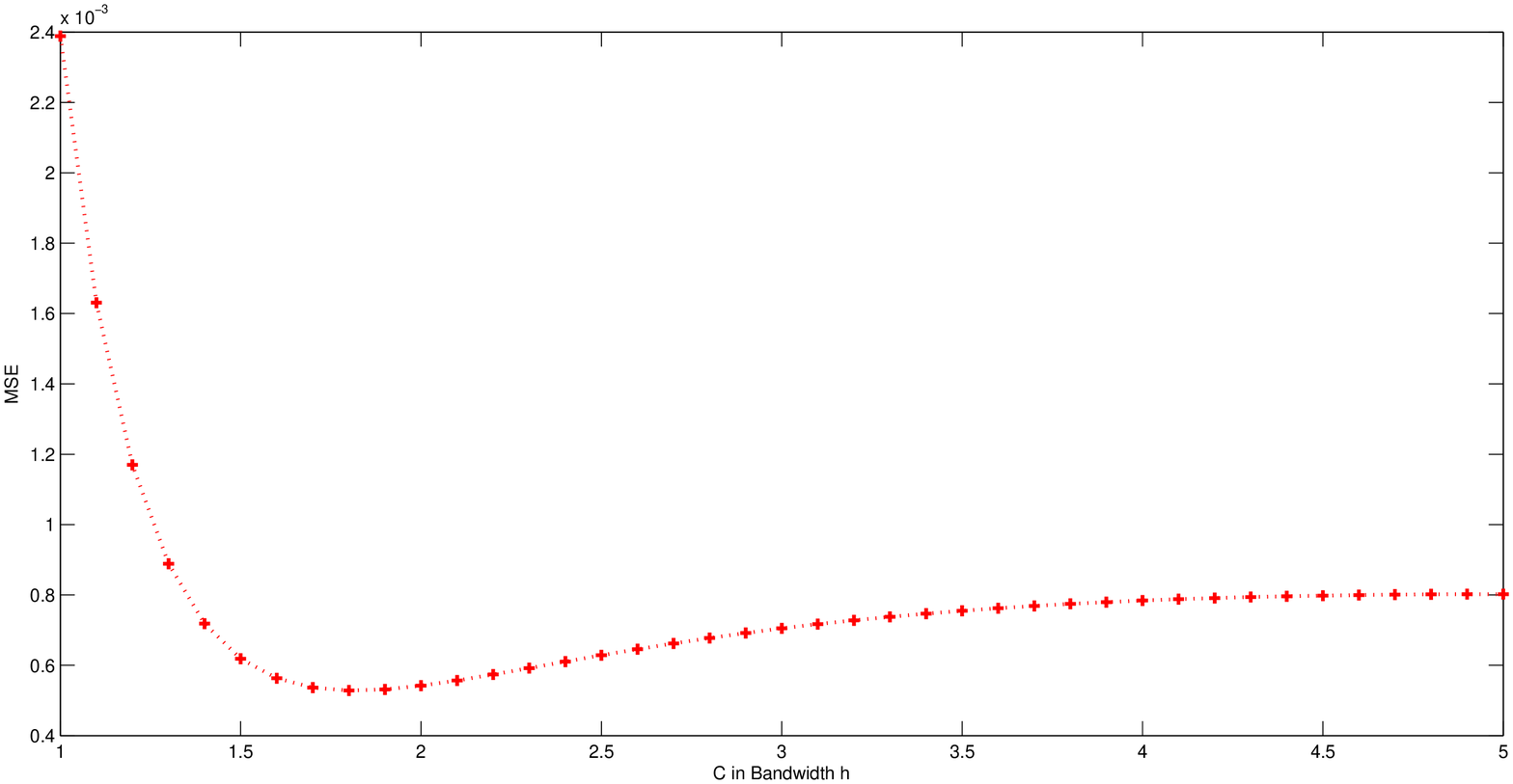}}
 \hspace{0.01in}
  \subfigure[ Curve of $CV(h)$ versus $h$ with $h_{cv} = 0.025$ ]{
  \label{fig:subfig} 
    \includegraphics[width=0.479\textwidth]{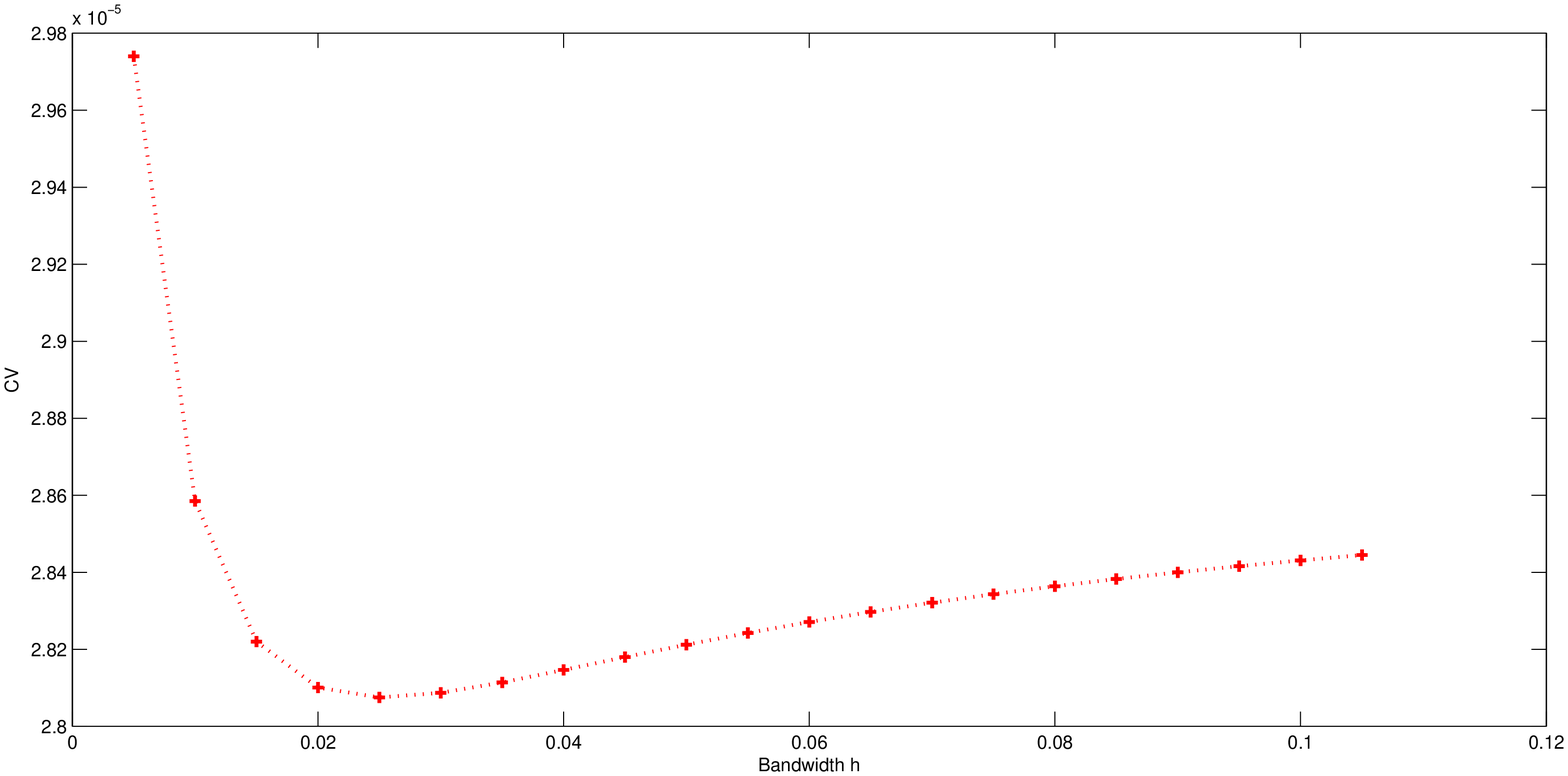}}
  \caption{ Curve for two rules of thumb versus $h$ with $T = 10,~n = 1000$ of $M(x) = 0.1 + 0.1*x^{2} + 2*0.036^{2}$ }
  \label{fig:subfig} 
\end{figure}

Firstly, we select the data-driven bandwidth by calculating MSE or
$k-$block CV as a function of $h_{n}$ from a sample with $T = 10$
and $n = 1000$ under two rules of thumb in Remark \ref{r3.3}~for
estimation of both $\mu(x)$ and $M(x)$, which are shown in Figure 3
and 4. We can get the optimal bandwidths $h_{n}$ for theses two
cases on estimating $\mu(x)$ by means of minimizing MSE or $k-$block
CV, which are $h_{opt} = 0.0683$ with $c_{opt} = 2.8$ which
coincides with that in Xu and Phillips \cite{xp} and $h_{cv} =
0.035.$ Although $h_{cv}$ is smaller than $h_{opt},$ the performance
of the estimator with $h_{cv}$ is a little worse than that with
$h_{opt},$ which can be tested and verified in Figures 5 and 6.

\begin{figure}[!htb]
\centering
\includegraphics[width=1 \textwidth]{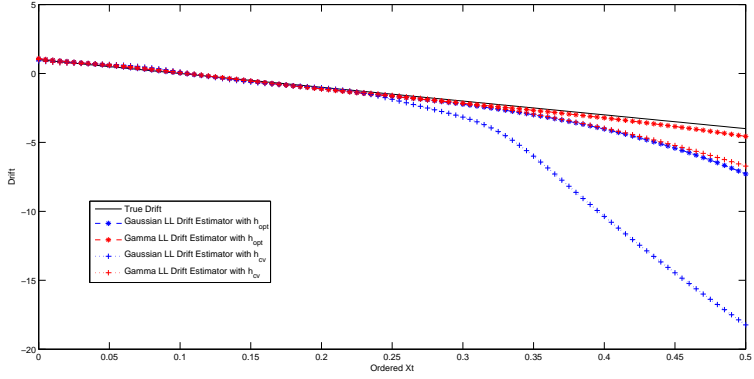}
\caption{ Local Linear Estimators for $\mu(x) = 1 - 10*x$ based on
Gaussian and Gamma kernels with $T = 10,~n = 1000,~h_{opt} =
2.8*\hat{S}*T^{-2/5} = 0.0683,~h_{cv} = 0.035$}
\end{figure}

\begin{figure}[!htb]
\centering
\includegraphics[width=1 \textwidth]{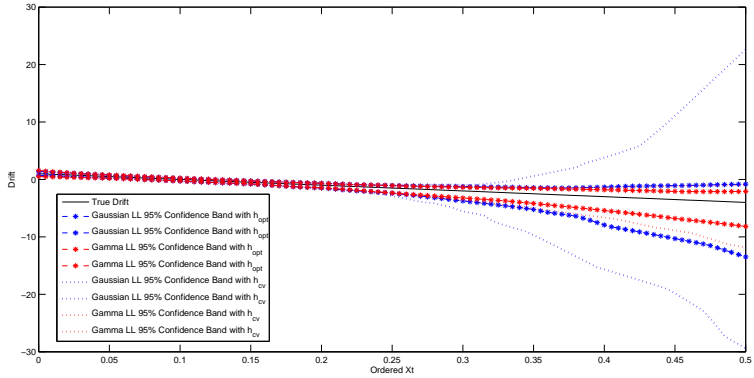}
\caption{ 95\% Monte Carlo confidence intervals for $\mu(x) = 1 -
10*x$ based on Gaussian and Gamma kernels with $T = 10,~n =
1000,~h_{opt} = 2.8*\hat{S}*T^{-2/5} = 0.0683,~h_{cv} = 0.035$}
\end{figure}

\begin{figure}[!htb]
\centering
\includegraphics[width=1 \textwidth]{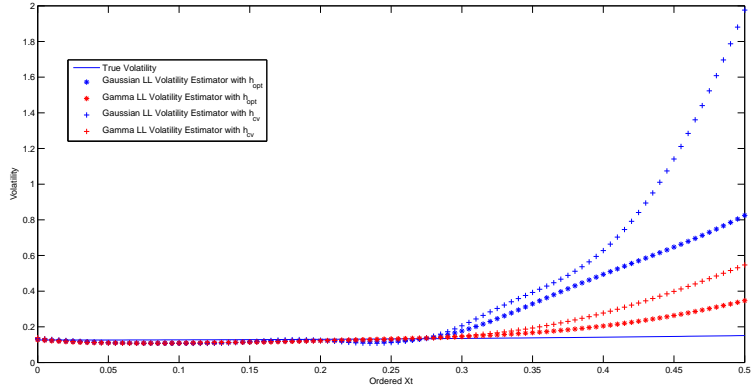}
\caption{ Local Linear Estimators for $M(x) = 0.1 + 0.1*x^{2} +
2*0.036^{2}$ based on Gaussian and Gamma kernels with $T = 10,~n =
1000,~h_{opt} = 1.7*\hat{S}*T^{-2/5} = 0.0441,~h_{cv} = 0.025$}
\end{figure}

\begin{figure}[!htb]
\centering
\includegraphics[width=1 \textwidth]{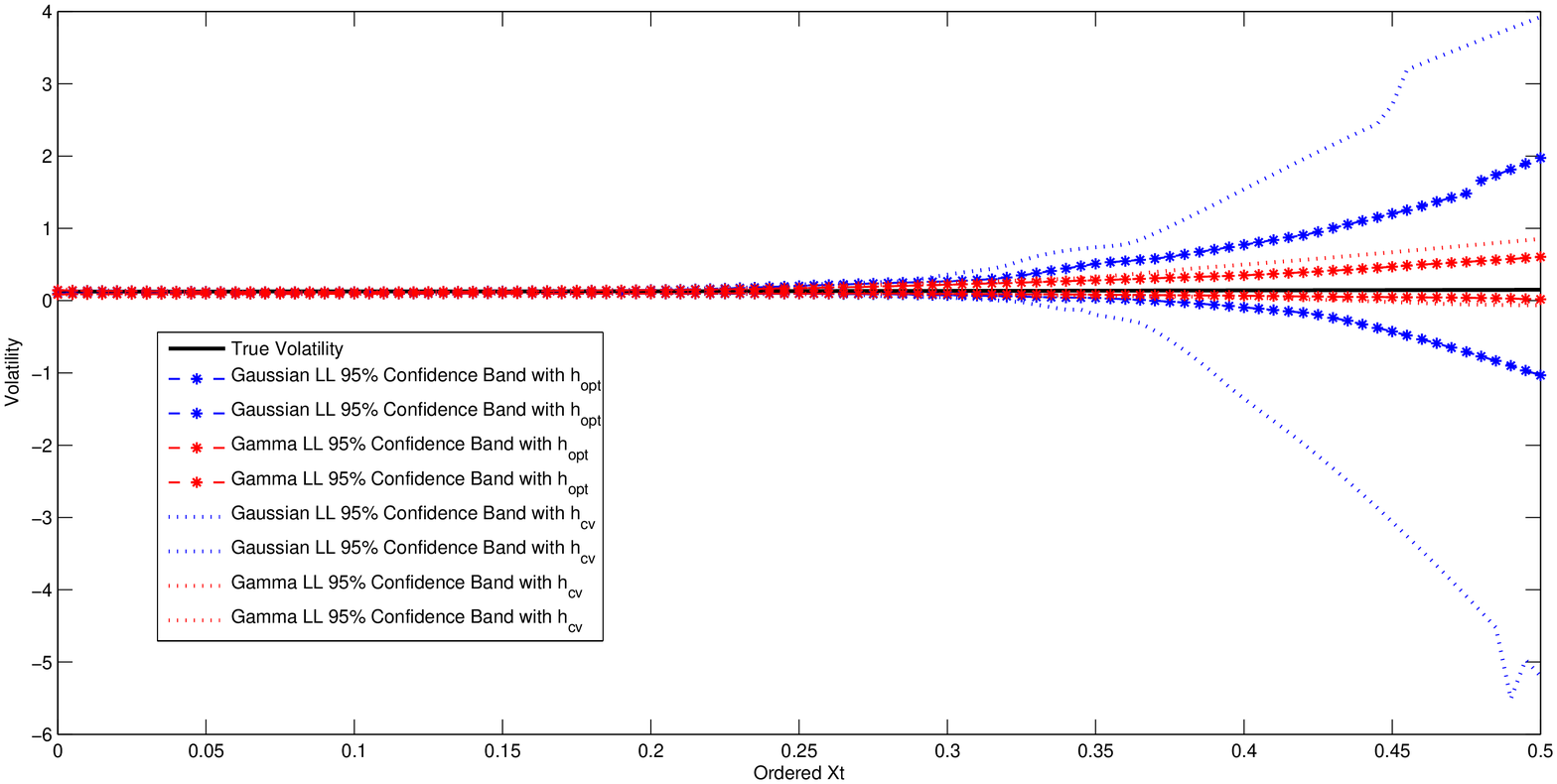}
\caption{ 95\% Monte Carlo confidence intervals for $M(x) = 0.1 +
0.1*x^{2} + 2*0.036^{2}$ based on Gaussian and Gamma kernels with $T
= 10,~n = 1000,~h_{opt} = 1.7*\hat{S}*T^{-2/5} = 0.0441,~h_{cv} =
0.025$}
\end{figure}

Figures 5 and 6 represent the local linear estimator and $95\%$
Monte Carlo confidence intervals constructed with Gamma asymmetric
kernels and Gaussian symmetric kernels for $\mu(x)$ from a sample
with $T = 10$ and $n = 1000$ under two rules of thumb for $h_{n},$
which show the local linear estimator constructed with Gamma
asymmetric kernels performs a little better and the $95\%$ Monte
Carlo confidence intervals constructed with Gamma asymmetric kernels
are shorter than that constructed with Gaussian symmetric kernels
which reveals smaller variability, especially at the sparse design
point. In addition, the local linear estimator constructed with
$h_{n} = c * \hat{S} * T^{-2/5}$ performs a little better and the
$95\%$ Monte Carlo confidence intervals constructed with $h_{n} = c
* \hat{S} * T^{-2/5}$ are shorter than that constructed with
$h_{cv}.$ In the subsequent numeral calculations such as MSE,
coverage rate and lengths of confidence band, we will calculate the
estimators with $h_{n} = c * \hat{S} * T^{-2/5}.$ Additionally, from
Figure 5 we can observe that the local linear estimator constructed
with Gaussian symmetric kernels exhibits a higher downward bias than
that constructed with Gamma asymmetric kernels which is practically
unbiased, which coincides with that in Gospodinov and Hirukawa
\cite{gh}.

Similar results for estimation of $M(x)$ can be observed from Figure
4, 7 and 8. Inconsistently, the optimal bandwidths $h_{n}$ on
estimating $M(x)$ by means of minimizing MSE or $k-$block CV are
$h_{opt} = 0.0441$ with $c_{opt} = 1.7$ which coincides with that in
Xu and Phillips \cite{xp} and $h_{cv} = 0.025.$ Compared with the
optimal smoothing parameter for $\mu(x),$ it is observed that a
narrower one for $M(x)$, as documented in Chapman and Pearson
\cite{cp}. From Figure 7, we can observe that the local linear
estimator constructed with Gaussian symmetric kernels exhibits a
higher upward bias than that constructed with Gamma asymmetric
kernels, which may be caused by the discretization bias similarly as
the microstructure noise in empirical market.

\begin{remark}
\label{r4.1} As proposed in Hirukawa and Sakudo \cite{hs}, the
``rule of thumb'' smoothing bandwidth should be under consideration
for each cases such as various $T$ and $n,$ which should be
determined by the specific financial data. As depicted in Figure 9,
contrary to $h_{opt}$ with $c_{opt} = 2.8$ for a sample with $T =
10$ and $n = 1000,$ the estimated parameters $c_{opt}$ in $h_{opt}$
based on a sample with $T = 50$ and $n = 5000$ or $T = 100$ and $n =
5000$ are 3 or 3.7. The $c_{opt}$ in $h_{opt}$ get larger as the
time span $T$ expands larger, which may be due to the fact that
$h_{opt} = c \cdot \hat{S} \cdot T^{-2/5}$ is inversely proportional
to $T$ and the fact that the smaller $h_{n}$, the larger bias from
Figure 3, 4 and 9. In order to effectively compare the local linear
estimator constructed with Gamma asymmetric kernel with that using
Gaussian symmetric kernel in terms of MSE, we should calculate the
$c_{opt}$ in $h_{opt}$ for various $T$ and $n.$ Here we omit these
calculations for $c_{opt}.$
\begin{figure}[!htb]
\centering
  \subfigure[ $C_{opt} = 3$ for $T = 50,~n = 5000$ with $h_{opt} = 0.0441$ (0.0412 if $C = 2.8$) ]{
  \label{fig:subfig} 
    \includegraphics[width=0.479\textwidth]{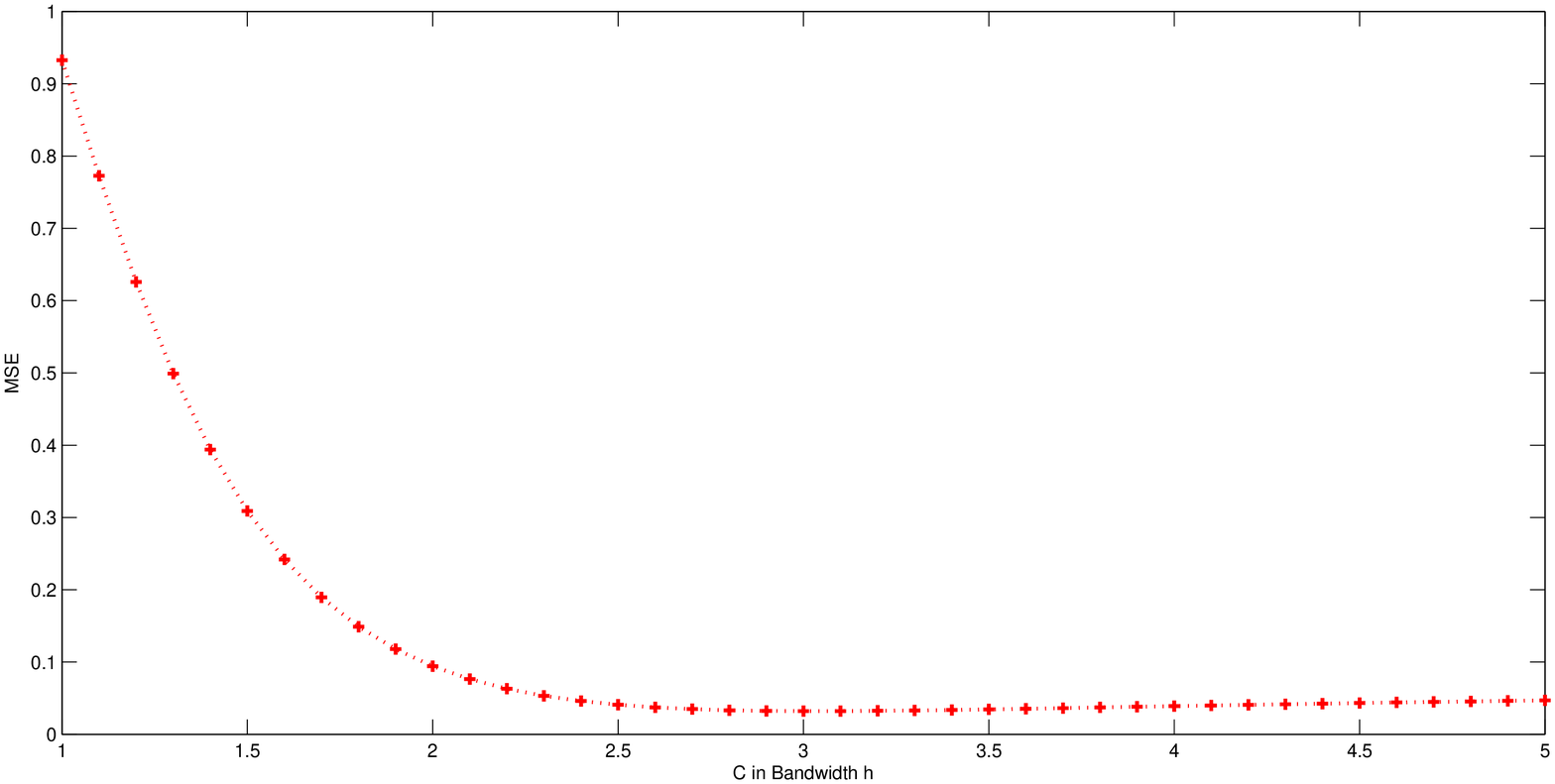}}
 \hspace{0.01in}
  \subfigure[ $C_{opt} = 3.7$ for $T = 100,~n = 5000$ with $h_{opt} = 0.0435$ (0.0329 if $C = 2.8$) ]{
  \label{fig:subfig} 
    \includegraphics[width=0.479\textwidth]{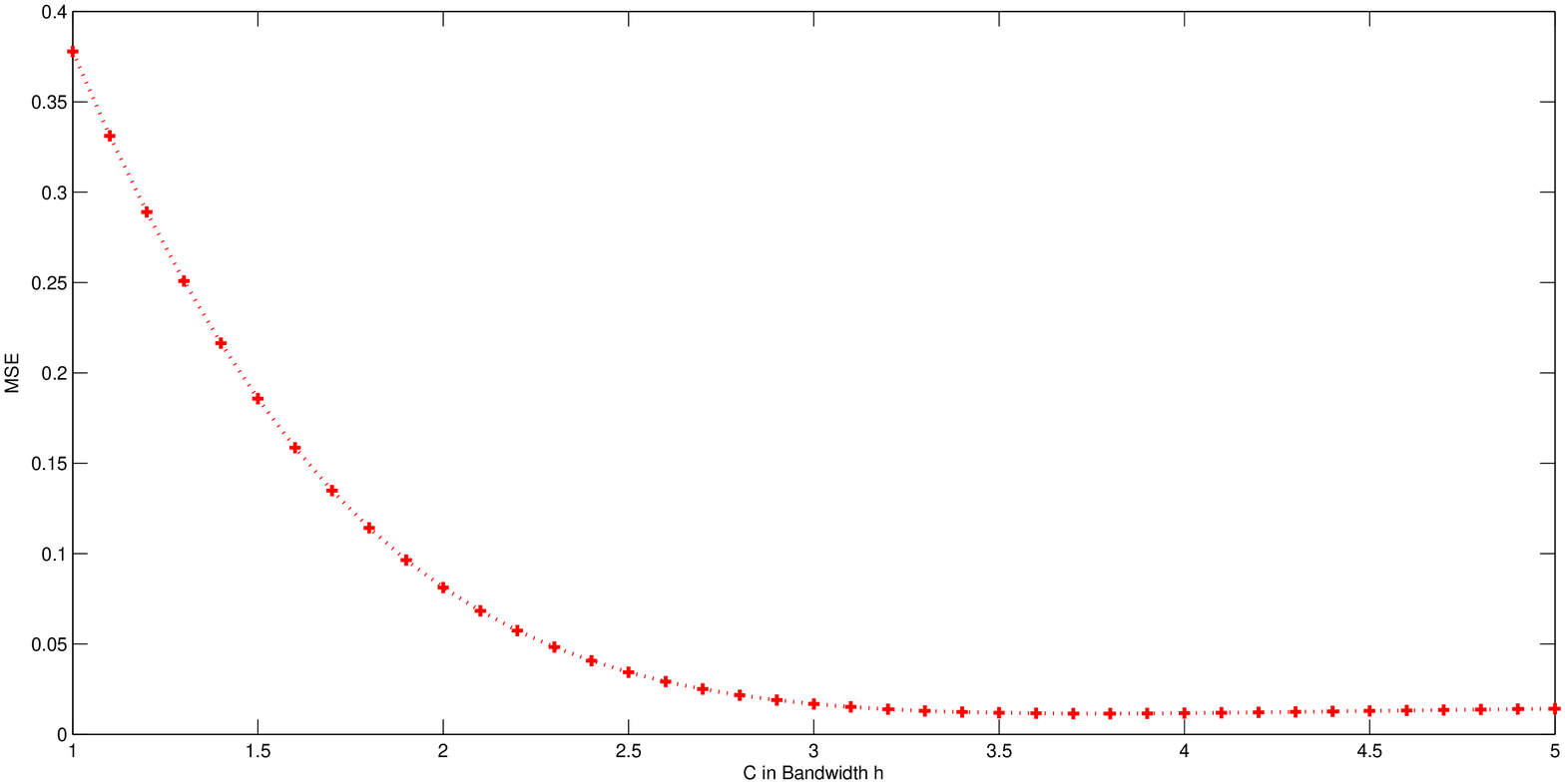}}
  \caption{ Curve for $MSE(h)$ versus $C$ in $h = C*\hat{S}*T^{-2/5}$ under $n = 5000$ and various $T$ of $\mu(x) = 1 - 10*x$ }
  \label{fig:subfig} 
\end{figure}
\end{remark}

We will assess the performance of the local linear estimators
constructed with Gamma asymmetric kernels and those constructed with
Gaussian symmetric kernels for both drift and conditional variance
functions via the Mean Square Errors (MSE)
\begin{equation}
\label{mc5} MSE = \frac{1}{m}\sum_{k=0}^{m}\left\{\hat{\mu}(x_{k}) -
\mu(x_{k})\right\}^{2}, \end{equation} where $\hat{\mu}(x)$ is the
estimator of $\mu(x)$ and $\{x_{k}\}_{1}^{m}$ are chosen uniformly
to cover the range of sample path of $X_{t}.$ Table 2 gives the
results on the MSE of local linear estimator constructed with Gamma
asymmetric kernels (MSE-LL(Gamma)) and local linear estimator
constructed with Gaussian symmetric kernels (MSE-LL(Gaussian)) for
the drift function $\mu(x)$ with jump size $Z_{n} \sim
\mathscr{N}(0,0.036^{2})$ over 500 replicates. Table 3 reports the
results on MSE-LL(Gamma) and MSE-LL(Gaussian) for $M(x)$.

From Table 2 and 3, we can make the following remarks.
\begin{itemize}
\item { The local linear estimator for $\mu(x)$ and $M(x)$ constructed with Gamma
asymmetric kernels performs a little better than that constructed
with Gaussian symmetric kernels in terms of MSE;}
\item { For the same time interval $T$, as the sample sizes $n$ tends larger, the
performances of the estimators for $\mu(x)$ or $M(x)$ are improved
due to the fact that more information for estimation procedure is
sampled as $\Delta_{n} \rightarrow 0;$}
\item { For the same sample sizes $n$, as the time interval $T$ expands
larger, the performance of the estimator for $\mu(x)$ is improved
due to the fact that the more information about drift coefficient is
obtained as the time span get larger, however, the performance of
the estimator for $M(x)$ gets worse, especially when $T = 100$, due
to the fact that more jumps happens in larger time interval $T$ in
steps 3 of Algorithm 1;}
\item { To some extent, the previous remark confirms
that the drift and conditional variance functions cannot be
identified in a fixed time span, which corresponds to the results of
Theorem \ref{thm2}.}
\end{itemize}

\begin{table*}[!htb]
\centering \caption{Simulation results on MSE-LL(Gaussian),
MSE-LL(Gamma) for three lengths of time interval (T) and three
sample sizes for $\mu(x) = 1 - 10 x$ with arrival intensity $\lambda
\cdot T = 20$, jump size $Z_{n} \sim \mathscr{N}(0,0.036^{2})$ and
$h_{opt}$ over 500 replicates.}
\begin{tabular}{crrrrrc}
\hline T & \multicolumn{1}{c}{$~$} & \multicolumn{1}{c}{$n = 500$} &
\multicolumn{1}{c}{$n = 1000$} & \multicolumn{1}{c}{$n = 5000$}\\
\hline
10             & MSE-LL(Gaussian)     & 0.5818         & 0.5253        & 0.0833    \\
               & MSE-LL(Gamma)        & 0.3010         & 0.1511        & 0.0424    \\
50             & MSE-LL(Gaussian)     & 0.5075         & 0.7499        & 0.1596    \\
               & MSE-LL(Gamma)        & 0.1373         & 0.0874        & 0.0154    \\
100            & MSE-LL(Gaussian)     & 0.0364         & 0.0351        & 0.0078    \\
               & MSE-LL(Gamma)        & 0.0187         & 0.0093        & 0.0022    \\
\hline
\end{tabular}
\end{table*}

\begin{table*}[!htb]
\centering \caption{Simulation results on MSE-LL(Gaussian),
MSE-LL(Gamma) for three lengths of time interval (T) and three
sample sizes for $M(x) = 0.1 + 0.1*x^{2} + \lambda*0.036^{2}$ with
arrival intensity $\lambda \cdot T = 20$, jump size $Z_{n} \sim
\mathscr{N}(0,0.036^{2})$ and $h_{opt}$ over 500 replicates.}
\begin{tabular}{crrrrrc}
\hline T & \multicolumn{1}{c}{$~$} & \multicolumn{1}{c}{$n = 500$} &
\multicolumn{1}{c}{$n = 1000$} & \multicolumn{1}{c}{$n = 5000$}\\
\hline
10             & MSE-LL(Gaussian)     & 0.1268         & 0.0337        & 0.0019    \\
               & MSE-LL(Gamma)        & 0.0084         & 0.0047        & $3.7014 \times 10^{-4}$    \\
50             & MSE-LL(Gaussian)     & 0.1153         & 0.0169        & 0.0079    \\
               & MSE-LL(Gamma)        & 0.0079         & $5.7058 \times 10^{-4}$        & $3.1990 \times 10^{-4}$    \\
100            & MSE-LL(Gaussian)     & 10.7405         & 0.0649        & 0.0089    \\
               & MSE-LL(Gamma)        & 27.9961         & 0.0025        & 0.0011    \\
\hline
\end{tabular}
\end{table*}

Figure 10 and 11 give the QQ plots for the local linear estimators
of the drift function $\mu(x)$ and conditional variance function
$M(x)$ constructed with Gamma asymmetric kernels and those
constructed with Gaussian symmetric kernels with $T = 50$ and
$\Delta_{n} = 0.01$. This reveals the normality of the local linear
estimators of the drift function $\mu(x)$ and conditional variance
function $M(x)$ constructed with Gamma asymmetric kernels, which
confirms the results in Theorems \ref{thm2}.

\begin{figure}[!htb]
\centering
  \subfigure[ QQ plot using Gaussian kernels  ]{
  \label{fig:subfig} 
    \includegraphics[width=0.479\textwidth]{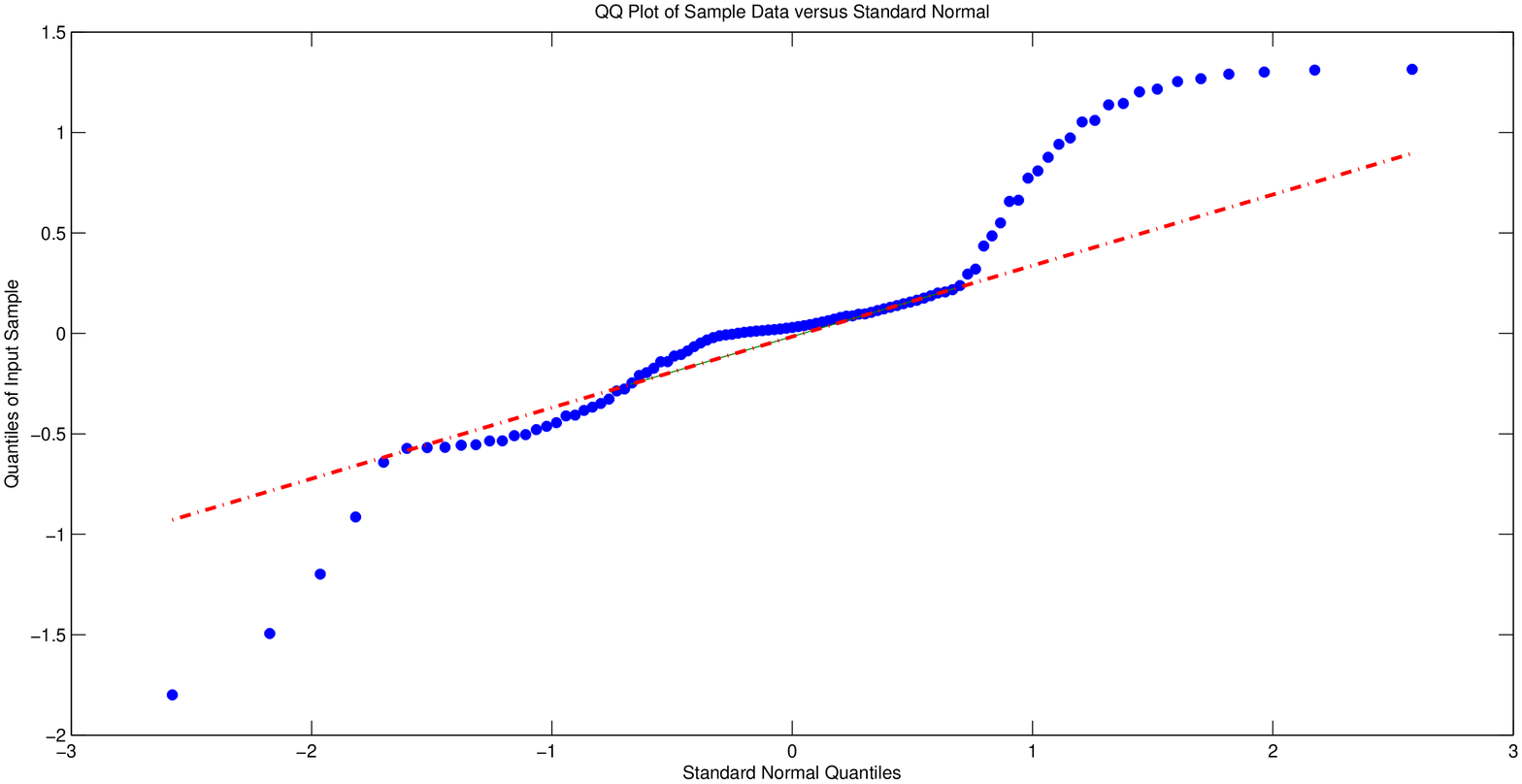}}
 \hspace{0.01in}
  \subfigure[ QQ plot using Gamma kernels ]{
  \label{fig:subfig} 
    \includegraphics[width=0.479\textwidth]{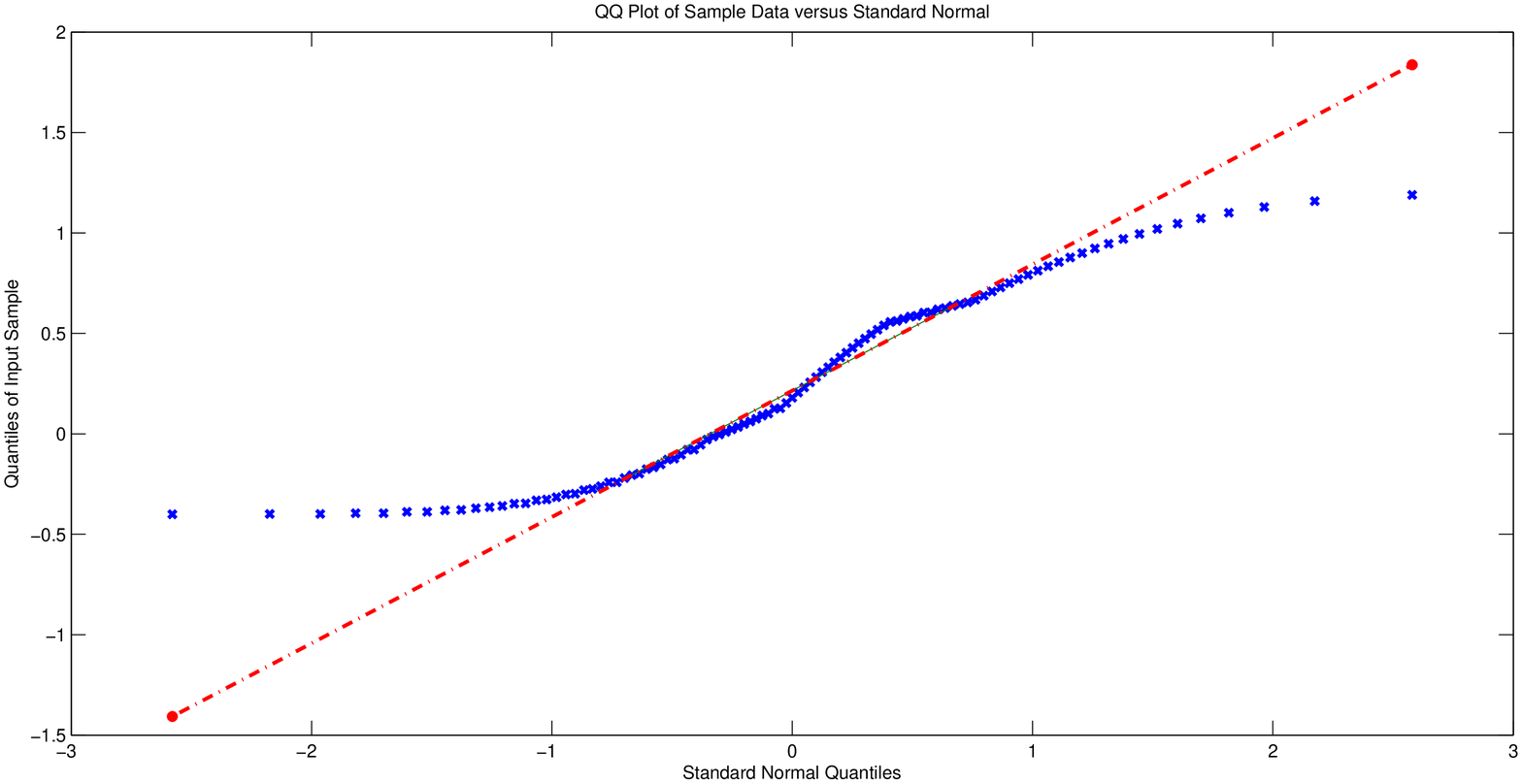}}
  \caption{ QQ plot of local linear estimators for $\mu(x) = 1 -
10*x$ using Gaussian and Gamma kernels with $T = 50,~n =
5000,~h_{opt} = 3*\hat{S}*T^{-2/5} = 0.0441$ }
  \label{fig:subfig} 
\end{figure}

\begin{figure}[!htb]
\centering
  \subfigure[ QQ plot using Gaussian kernels  ]{
  \label{fig:subfig} 
    \includegraphics[width=0.479\textwidth]{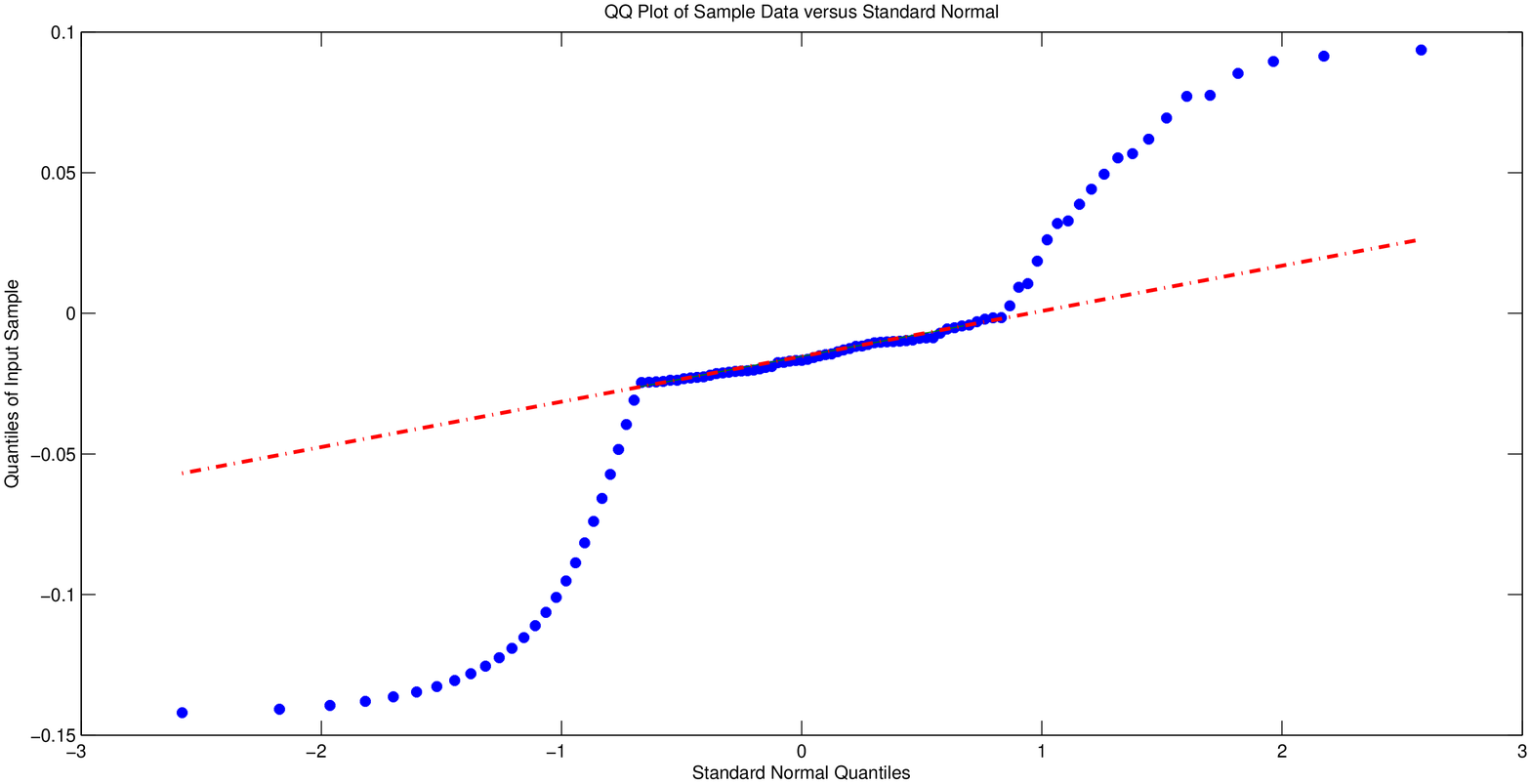}}
 \hspace{0.01in}
  \subfigure[ QQ plot using Gamma kernels ]{
  \label{fig:subfig} 
    \includegraphics[width=0.479\textwidth]{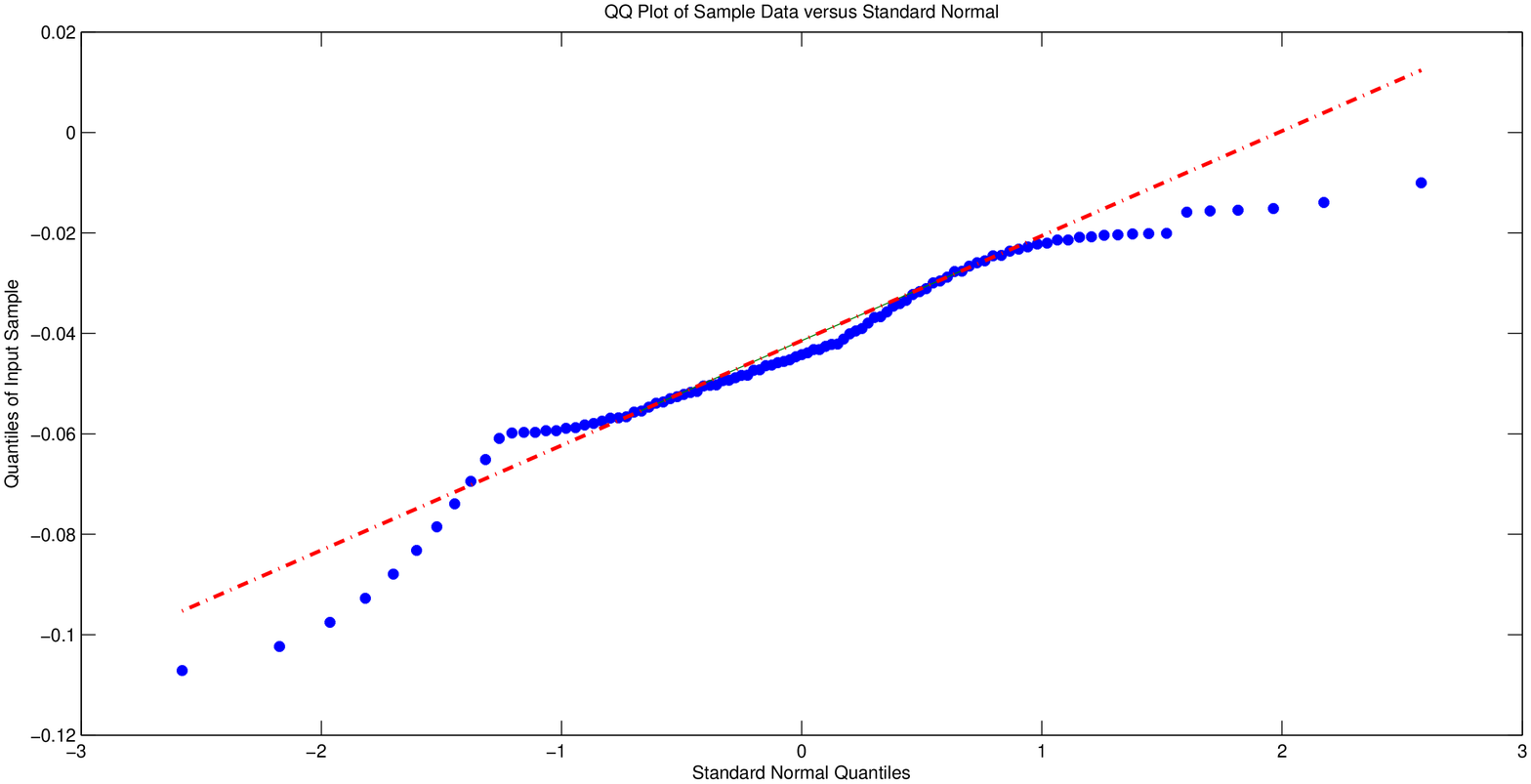}}
  \caption{ QQ plot of local linear estimators for $M(x) = 0.1 +
0.1*x^{2} + 20/50*0.036^{2}$ using Gaussian and Gamma kernels with
$T = 50,~n = 5000,~h_{opt} = 1.9*\hat{S}*T^{-2/5} = 0.0294$ }
  \label{fig:subfig} 
\end{figure}

The computational results about comparison of coverage rate and
length of confidence band for drift $\mu(x)$ are summarized in Table
4 - 6 and conditional variance $M(x)$ in Table 9 for various
intensity and amplitude of jumps with $T = 50,~n = 5000.$ The
confidence intervals are constructed as those in Remark \ref{r3.5}.
For ``boundary $x = 0.05$'' or ``interior $x = 0.15$'', we also
calculate the adjusted lengths of confidence band for $\mu(x)$ which
was argued for in Xu \cite{xkl}. The adjusted lengths of the
calibrated confidence intervals constructed with the actual critical
values are presented in Table 7 - 8. The actual critical values
marked as $2.5\%$ or $97.5\%$ quantile in Table 7 - 8 are adapt to
actual simulation data and no longer 1.96 as that in the standard
normal distribution table for asymptotic normality such that the
coverage rates are adjusted to $95\%.$ Note that the ratio in tables
denotes the ratio of the length of confidence band constructed with
Gaussian symmetric kernel (CB-GSK) to that constructed with Gamma
asymmetric kernel (CB-GAK). From these six tables, we can obtain the
following findings.
\begin{itemize}
\item { From Table 4, for the drift function, according to the coverage rates in percent,
CB-GSK or CB-GAK are more under-covered as $h_{n}$ expands larger,
especially at the ``boundary $x = 0.05$'' or the sparse ``interior
$x = 0.30$''. When $x = 0.15,$ CB-GSK or CB-GAK has favorable
lengths for the bandwidths $h_{n}.$ The reasons behind the
phenomenon may be the following three aspects. Firstly, estimations
for asymptotic variance at design points $x = 0.05, 0.30$ based on
the formula in Theorem \ref{thm2} are more slightly below the sample
variance of the randomly generated data, which may be the main
reason. Secondly, the process visits the design points $x = 0.05,
0.30$ less frequently such that there are relatively less sample
points near their neighborhood to estimate the unknown quantity. The
first two reasons coincide with those observed in Xu \cite{xkl}.
Thirdly, as shown in Remark \ref{r3.5}, the bias in the normal
confidence interval is estimated by taking the second derivative of
the local linear estimators,which may give rise to estimation bias.
Also, one can see that at the ``boundary point x = 0.05'', the
coverage rate with GSK is a little better than that with GAK, which
may be due to the facts: CB-GSK can allocate weight to the
observations less than zero while CB-GAK can't.}

\item { From Table 4 - 6, the absolute mean of bias and variance for estimator
constructed with Gamma asymmetric kernels are less than that
constructed with Gaussian symmetric kernels, which indicates that
compared with that constructed with Gaussian symmetric kernels,
estimator constructed with Gamma asymmetric kernels is practically
unbiased and exhibits smaller variability for either the boundary
point or the spare design point. This coincides with the discussion
of coverage rate in Remark \ref{r3.6}, as documented in Gospodinov
and Hirukawa \cite{gh}.}

\item { From Table 4 - 6, it is observed that for the same design point $x$,
the ratio is more than 1 under the condition that the confidence
rate is approximate, that is CB-GSK is longer than CB-GAK.
Meanwhile, CB-GSK or CB-GAK expands larger and the ratio becomes
bigger as $x$ increases, especially at the sparse design point $x =
0.30.$ One can also find that the ratio becomes smaller as the
smoothing parameter $h_{n}$ gets larger, especially at the boundary
design point $x = 0.05,$ which coincides with the theoretical
results discussed in Remark \ref{r3.6}. Note that the smaller the
smoothing parameter $h_{n}$, the larger the bias. Hence the choice
of smoothing parameter $h_{n}$ is not recommended to be too large or
too small and should depend on the need for higher coverage rate or
lower bias. As the intensity and amplitude of jump increases, CB-GSK
or CB-GAK expands larger to improve coverage rate since it needs
more information to cover more and larger jump similarly as that
mentioned in Bandi and Nguyen \cite{bn}.}

\item { Here we discuss the adjusted lengths of confidence intervals constructed with various
kernels since the two findings above indicate that the reasonably
correct coverage rate depends on the choice of smoothing bandwidth
$h_{n}$ and the consistently estimated variance. From Table 7 for
``boundary $x = 0.05$'', after confidence rates are adjusted to
$95\%,$ the adjusted lengths of confidence intervals constructed
with Gamma asymmetric kernel (ACI-GAK) are shorter than that
constructed with Gaussian symmetric kernel (ACI-GSK) with the
smoothing parameter $h_{n} = 0.01, 0.02, 0.03,$ however, ACI-GAK are
longer than ACI-GSK with $h_{n} = 0.0441, 0.05,$ which coincides
with the conclusion in Remark \ref{r3.6} that the closer to the
boundary point or the larger bandwidth for fixed ``boundary $x$'',
the shorter the length of confidence interval based on Gaussian
symmetric kernel. In addition, the lengths of ACI-GAK are more
robust as the smoothing parameter $h_{n}$ changes. One can also
observe that the absolute critical values for ACI-GAK or ACI-GSK are
almost larger than 1.96 for adjusted lengths, which provides a
reference point for empirical analysis when constructing the
confidence intervals. In Table 8, for ``interior $x = 0.15$'',
ACI-GAK are shorter than ACI-GSK with all $h_{n}$ mentioned and the
absolute critical values for ACI-GAK or ACI-GSK are almost larger
than 1.96 for adjusted lengths. Constructed with the Gamma
asymmetric kernel, the adjustment for $\mu(0.15)$ is milder relative
to $\mu(0.05)$ for any given bandwidth or intensity and amplitude of
jumps, which is because that there are plenty of sample points used
to estimate $\mu(x)$ near the design point $x = 0.15$ such that the
coverage rate for $\mu(0.15)$ is close to $95\%$ in Table 4 - 6. }

\item {For simplicity, similar observations for the conditional variance $M(x)$ only at ``interior $x = 0.15$''
are shown in Table 9. Similarly as the drift function, the ratio is
more than 1 under the condition that the confidence rate is
approximate, that is CB-GAK is shorter than CB-GSK. Meanwhile, the
ratio becomes smaller as the smoothing parameter $h_{n}$ gets larger
and the intensity and amplitude of jump increases. Compared with the
mean of bias of the estimator for $\mu(0.15)$ in Table 4 - 6, the
mean of bias for $M(0.15)$ constructed with Gamma asymmetric kernels
is larger than that constructed with Gaussian symmetric kernels. As
mentioned in Theorem \ref{thm2} and Remark \ref{r3.5}, for the
estimators constructed with Gamma asymmetric kernels at ``interior
$x = 0.15$'', the bias for ${\mu}(x)$ in model (\ref{mc1}) is $h_{n}
B_{\hat{\mu}_{n}(x)} = h_{n} \frac{x}{2}\mu^{''}(x) \equiv 0,$ which
is equal to the bias of the estimator constructed with Gaussian
symmetric kernels, that is also $0.$ However, the bias for ${M}(x)$
is $h_{n} B_{\hat{M}_{n}(x)} = h_{n} \frac{x}{2}M^{''}(x) = h_{n}
\cdot 0.1 \cdot x =  O(h_{n})$ which is far greater than the bias of
the estimator constructed with Gaussian symmetric kernels, that is
$O(h_{n}^{2}).$ Furthermore, this observation does not contradict
with the result on MSE considered previously because the variance of
the estimator for $M(0.15)$ constructed with Gamma asymmetric
kernels is less than that constructed with Gaussian symmetric
kernels and the variance dominates the MSE.}

\end{itemize}

\begin{sidewaystable}[!htb]
  \centering
  \caption{Estimation for $\mu(x)$ with arrival intensity $\lambda \cdot T = 20$, jump size $Z_{n} \sim \mathscr{N}(0,0.036^{2})$ over 500 replicates.}
    \begin{tabular}{cccccccccccccc}
    \toprule
    \multirow{2}[4]{*}{Bandwidth} & \multicolumn{2}{c}{Coverage Rate} & \multicolumn{2}{c}{Mean of Bias} & \multicolumn{2}{c}{Variance} & \multicolumn{4}{c}{Estimation for Variance} & \multicolumn{3}{c}{Length of Confidence Band} \\
\cmidrule{2-14}          & Sym   & Asym  & Sym   & Asym  & Sym   & Asym  & Sym   & Std   & Asym  & Std   & Sym   & Asym  & Ratio \\
    \midrule
    \multicolumn{14}{c}{x = 0.05} \\
    \midrule
    $h_{1} = 0.01$ & 94.6  & 92.4  & 0.0045 & 0.0028 & 0.0144 & 0.0067 & 0.0136 & 0.0011 & 0.0057 & 0.0004 & 0.457 & 0.2954 & 1.5471  \\
    \midrule
    $h_{2} = 0.02$ & 94.6  & 89.6  & 0.0032 & 0.0011 & 0.007 & 0.0055 & 0.0069 & 0.0005 & 0.0039 & 0.0002 & 0.3262 & 0.2456 & 1.3282  \\
    \midrule
    $h_{3} = 0.03$ & 94.2  & 88    & 0.0022 & 0.0003 & 0.0049 & 0.005 & 0.0048 & 0.0003 & 0.0032 & 0.0002 & 0.2701 & 0.2217 & 1.2183  \\
    \midrule
    $h_{4} = 0.05$ & 93.6  & 84    & 0.0018 & -0.0002 & 0.0037 & 0.0046 & 0.0031 & 0.0002 & 0.0025 & 0.0001 & 0.2178 & 0.1976 & 1.1022  \\
    \midrule
    $h_{opt} = 0.0441$ & 93.8  & 84.6  & 0.0018 & -0.0001 & 0.0038 & 0.0047 & 0.0033 & 0.0002 & 0.0026 & 0.0001 & 0.2253 & 0.2012 & 1.1198  \\
    \midrule
    \multicolumn{14}{c}{x = 0.15} \\
    \midrule
    $h_{1} = 0.01$ & 94.4  & 94.2  & -0.0038 & -0.004 & 0.0139 & 0.0042 & 0.0142 & 0.0013 & 0.0042 & 0.0003 & 0.4659 & 0.2538 & 1.8357  \\
    \midrule
    $h_{2} = 0.02$ & 95.6  & 94    & -0.005 & -0.0036 & 0.0072 & 0.0035 & 0.0072 & 0.0005 & 0.0033 & 0.0002 & 0.332 & 0.225 & 1.4756  \\
    \midrule
    $h_{3} = 0.03$ & 93.2  & 94    & -0.0045 & -0.0035 & 0.0051 & 0.0032 & 0.0049 & 0.0003 & 0.0029 & 0.0001 & 0.2745 & 0.212 & 1.2948  \\
    \midrule
    $h_{4} = 0.05$ & 93.6  & 93.4  & -0.0034 & -0.0033 & 0.0036 & 0.003 & 0.0032 & 0.0002 & 0.0026 & 0.0001 & 0.2208 & 0.1993 & 1.1079  \\
    \midrule
    $h_{opt} = 0.0441$ & 93.6  & 93.2  & -0.0036 & -0.0033 & 0.0038 & 0.003 & 0.0034 & 0.0002 & 0.0026 & 0.0001 & 0.2285 & 0.2012 & 1.1357  \\
    \midrule
    \multicolumn{14}{c}{x = 0.30} \\
    \midrule
    $h_{1} = 0.01$ & 94.4  & 95    & -0.1856 & -0.0913 & 0.863 & 0.14  & 0.7196 & 0.5996 & 0.125 & 0.0538 & 3.1606 & 1.3587 & 2.3262  \\
    \midrule
    $h_{2} = 0.02$ & 89.6  & 89.2  & -0.1702 & -0.0581 & 0.4349 & 0.0853 & 0.2866 & 0.1626 & 0.058 & 0.0196 & 2.0382 & 0.9319 & 2.1871  \\
    \midrule
    $h_{3} = 0.03$ & 86.2  & 87.2  & -0.138 & -0.0421 & 0.2695 & 0.0644 & 0.1443 & 0.0607 & 0.0369 & 0.0109 & 1.4629 & 0.7458 & 1.9615  \\
    \midrule
    $h_{4} = 0.05$ & 76.6  & 81    & -0.0806 & -0.0277 & 0.125 & 0.0474 & 0.0477 & 0.0119 & 0.0214 & 0.0052 & 0.8498 & 0.5691 & 1.4932  \\
    \midrule
    $h_{opt} = 0.0441$ & 79.6  & 81.8  & -0.0929 & -0.0299 & 0.1486 & 0.05  & 0.0596 & 0.0196 & 0.0235 & 0.0064 & 0.946 & 0.5962 & 1.5867  \\
    \bottomrule
    \end{tabular}%
  \label{tab:addlabel}%
\end{sidewaystable}
\clearpage

\begin{sidewaystable}[!htb]
  \centering
  \caption{Estimation for $\mu(x)$ with arrival intensity $\lambda \cdot T = 50$, jump size $Z_{n} \sim \mathscr{N}(0,0.036^{2})$ over 500 replicates.}
    \begin{tabular}{cccccccccccccc}
    \toprule
    \multirow{2}[4]{*}{Bandwidth} & \multicolumn{2}{c}{Coverage Rate} & \multicolumn{2}{c}{Mean of Bias} & \multicolumn{2}{c}{Variance} & \multicolumn{4}{c}{Estimation for Variance} & \multicolumn{3}{c}{Length of Confidence Band} \\
\cmidrule{2-14}          & Sym   & Asym  & Sym   & Asym  & Sym   & Asym  & Sym   & Std   & Asym  & Std   & Sym   & Asym  & Ratio \\
    \midrule
    \multicolumn{14}{c}{x = 0.05} \\
    \midrule
    $h_{1} = 0.01$ & 95.2  & 92    & 0.0027 & 0.0005 & 0.0138 & 0.0069 & 0.0137 & 0.0013 & 0.0057 & 0.0004 & 0.4586 & 0.2965 & 1.5467  \\
    \midrule
    $h_{2} = 0.02$ & 95    & 90    & 0.0019 & 0.0006 & 0.0068 & 0.0056 & 0.007 & 0.0005 & 0.004 & 0.0002 & 0.3271 & 0.2465 & 1.3270  \\
    \midrule
    $h_{3} = 0.03$ & 94.6  & 88.8  & 0.0018 & 0.0007 & 0.0049 & 0.0051 & 0.0048 & 0.0003 & 0.0032 & 0.0002 & 0.2709 & 0.2226 & 1.2170  \\
    \midrule
    $h_{4} = 0.05$ & 92.6  & 86.4  & 0.002 & 0.0009 & 0.0036 & 0.0047 & 0.0031 & 0.0002 & 0.0026 & 0.0001 & 0.2185 & 0.1984 & 1.1013  \\
    \midrule
    $h_{opt} = 0.0441$ & 93.8  & 86.8  & 0.0019 & 0.0009 & 0.0038 & 0.0047 & 0.0033 & 0.0002 & 0.0027 & 0.0001 & 0.2259 & 0.202 & 1.1183  \\
    \midrule
    \multicolumn{14}{c}{x = 0.15} \\
    \midrule
    $h_{1} = 0.01$ & 95.4  & 94.2  & -0.0029 & -0.0091 & 0.0134 & 0.0045 & 0.0143 & 0.0013 & 0.0043 & 0.0003 & 0.469 & 0.2556 & 1.8349  \\
    \midrule
    $h_{2} = 0.02$ & 95.8  & 93.4  & -0.0071 & -0.0094 & 0.007 & 0.0039 & 0.0073 & 0.0005 & 0.0033 & 0.0002 & 0.3343 & 0.2265 & 1.4759  \\
    \midrule
    $h_{3} = 0.03$ & 94.4  & 92.6  & -0.0085 & -0.0095 & 0.0052 & 0.0036 & 0.005 & 0.0003 & 0.003 & 0.0001 & 0.2764 & 0.2135 & 1.2946  \\
    \midrule
    $h_{4} = 0.05$ & 92.2  & 91.2  & -0.009 & -0.0096 & 0.004 & 0.0034 & 0.0032 & 0.0002 & 0.0026 & 0.0001 & 0.2222 & 0.2006 & 1.1077  \\
    \midrule
    $h_{opt} = 0.0441$ & 92.6  & 91.2  & -0.009 & -0.0096 & 0.0042 & 0.0035 & 0.0034 & 0.0002 & 0.0027 & 0.0001 & 0.2299 & 0.2025 & 1.1353  \\
    \midrule
    \multicolumn{14}{c}{x = 0.30} \\
    \midrule
    $h_{1} = 0.01$ & 93.4  & 94    & -0.1792 & -0.0942 & 0.8427 & 0.1411 & 0.7276 & 0.6509 & 0.1244 & 0.0672 & 3.1579 & 1.3501 & 2.3390  \\
    \midrule
    $h_{2} = 0.02$ & 92.2  & 87.8  & -0.1605 & -0.06 & 0.4074 & 0.0893 & 0.2826 & 0.1555 & 0.0576 & 0.0217 & 2.0249 & 0.9271 & 2.1841  \\
    \midrule
    $h_{3} = 0.03$ & 89.2  & 83    & -0.1365 & -0.0431 & 0.257 & 0.0679 & 0.1429 & 0.058 & 0.0367 & 0.0115 & 1.4564 & 0.7426 & 1.9612  \\
    \midrule
    $h_{4} = 0.05$ & 76    & 77.8  & -0.083 & -0.0275 & 0.1272 & 0.0495 & 0.0476 & 0.0114 & 0.0213 & 0.0054 & 0.8491 & 0.5674 & 1.4965  \\
    \midrule
    $h_{opt} = 0.0441$ & 79.2  & 78.6  & -0.0945 & -0.0298 & 0.1483 & 0.0522 & 0.0588 & 0.0187 & 0.0233 & 0.0066 & 0.9396 & 0.5929 & 1.5848  \\
    \bottomrule
    \end{tabular}%
  \label{tab:addlabel}%
\end{sidewaystable}
\clearpage

\begin{sidewaystable}[!htb]
  \centering
  \caption{Estimation for $\mu(x)$ with arrival intensity $\lambda \cdot T = 20$, jump size $Z_{n} \sim \mathscr{N}(0,0.1^{2})$ over 500 replicates.}
    \begin{tabular}{cccccccccccccc}
    \toprule
    \multirow{2}[4]{*}{Bandwidth} & \multicolumn{2}{c}{Coverage Rate} & \multicolumn{2}{c}{Mean of Bias} & \multicolumn{2}{c}{Variance} & \multicolumn{4}{c}{Estimation for Variance} & \multicolumn{3}{c}{Length of Confidence Band} \\
\cmidrule{2-14}          & Sym   & Asym  & Sym   & Asym  & Sym   & Asym  & Sym   & Std   & Asym  & Std   & Sym   & Asym  & Ratio \\
    \midrule
    \multicolumn{14}{c}{x = 0.05} \\
    \midrule
    $h_{1} = 0.01$ & 94.2  & 93.6  & -0.0114 & -0.0074 & 0.0142 & 0.0071 & 0.0142  & 0.0014  & 0.0059  & 0.0004  & 0.4671 & 0.3019 & 1.5472  \\
    \midrule
    $h_{2} = 0.02$ & 95.6  & 88.6  & -0.0068 & -0.0056 & 0.0075 & 0.006 & 0.0072  & 0.0006  & 0.0041  & 0.0003  & 0.3333 & 0.2509 & 1.3284  \\
    \midrule
    $h_{3} = 0.03$ & 94    & 86    & -0.0045 & -0.0046 & 0.0055 & 0.0055 & 0.0050  & 0.0003  & 0.0033  & 0.0002  & 0.2758 & 0.2264 & 1.2182  \\
    \midrule
    $h_{4} = 0.05$ & 91.2  & 83.2  & -0.0019 & -0.0035 & 0.004 & 0.0051 & 0.0032  & 0.0002  & 0.0026  & 0.0001  & 0.2222 & 0.2017 & 1.1016  \\
    \midrule
    $h_{opt} = 0.0441$ & 91.2  & 83    & -0.0022 & -0.0036 & 0.0042 & 0.0051 & 0.0034  & 0.0002  & 0.0027  & 0.0001  & 0.2283 & 0.2046 & 1.1158  \\
    \midrule
    \multicolumn{14}{c}{x = 0.15} \\
    \midrule
    $h_{1} = 0.01$ & 95    & 94.8  & -0.0073 & -0.005 & 0.0142 & 0.0043 & 0.0147 & 0.0015 & 0.0044 & 0.0003 & 0.4753 & 0.2587 & 1.8373  \\
    \midrule
    $h_{2} = 0.02$ & 94.4  & 94    & -0.0062 & -0.0045 & 0.0073 & 0.0035 & 0.0075 & 0.0006 & 0.0034 & 0.0002 & 0.3387 & 0.2292 & 1.4777  \\
    \midrule
    $h_{3} = 0.03$ & 93.8  & 94    & -0.0055 & -0.0043 & 0.0051 & 0.0033 & 0.0051 & 0.0003 & 0.003 & 0.0002 & 0.2799 & 0.2159 & 1.2964  \\
    \midrule
    $h_{4} = 0.05$ & 93.4  & 93.2  & -0.0043 & -0.004 & 0.0037 & 0.0031 & 0.0033 & 0.0002 & 0.0027 & 0.0001 & 0.225 & 0.203 & 1.1084  \\
    \midrule
    $h_{opt} = 0.0441$ & 93.4  & 93.2  & -0.0045 & -0.0041 & 0.0038 & 0.0031 & 0.0035 & 0.0002 & 0.0027 & 0.0001 & 0.2314 & 0.2045 & 1.1315  \\
    \midrule
    \multicolumn{14}{c}{x = 0.30} \\
    \midrule
    $h_{1} = 0.01$ & 95    & 94.4  & -0.1758 & -0.0828 & 0.764 & 0.1299 & 0.6683 & 0.5664 & 0.1073 & 0.0468 & 3.0489 & 1.2594 & 2.4209  \\
    \midrule
    $h_{2} = 0.02$ & 90.8  & 88.4  & -0.1562 & -0.0542 & 0.3728 & 0.0807 & 0.2704 & 0.1552 & 0.051 & 0.0176 & 1.9791 & 0.8742 & 2.2639  \\
    \midrule
    $h_{3} = 0.03$ & 87.6  & 84.6  & -0.1256 & -0.0402 & 0.2393 & 0.0616 & 0.1392 & 0.0619 & 0.033 & 0.0099 & 1.4352 & 0.7045 & 2.0372  \\
    \midrule
    $h_{4} = 0.05$ & 80.8  & 79.8  & -0.0734 & -0.0275 & 0.1176 & 0.0461 & 0.0475 & 0.0131 & 0.0194 & 0.0049 & 0.8471 & 0.5425 & 1.5615  \\
    \midrule
    $h_{opt} = 0.0441$ & 82.6  & 79.8  & -0.0823 & -0.0291 & 0.1355 & 0.0481 & 0.0569 & 0.0199 & 0.021 & 0.0058 & 0.9234 & 0.5628 & 1.6407  \\
    \bottomrule
    \end{tabular}%
  \label{tab:addlabel}%
\end{sidewaystable}
\clearpage

\begin{sidewaystable}[!htb]
  \centering
  \caption{Adjusted length of confidence band for $\mu(0.05)$ over 500 replicates.}
    \begin{tabular}{cccccccccc}
    \toprule
    \multirow{2}[4]{*}{Bandwidth} & \multicolumn{2}{c}{Coverage Rate} & \multicolumn{2}{c}{Sym Quantile} & \multicolumn{2}{c}{Asym Quantile} & \multicolumn{3}{c}{Length of Confidence Band} \\
\cmidrule{2-10}          & Sym   & Asym  & 2.50\% & 97.50\% & 2.50\% & 97.50\% & Sym   & Asym  & Ratio \\
    \midrule
    \multicolumn{10}{c}{Arrival intensity $\lambda \cdot T = 20$, Jump size $Z_{n} \sim \mathscr{N}(0,0.036^{2})$} \\
    \midrule
    $h_{1} = 0.01$ & 95.2  & 95.2  & -2.1465  & 1.8766 & -2.2883  & 2.2854 & 0.4703 & 0.3453 & 1.362 \\
    \midrule
    $h_{2} = 0.02$ & 95.2  & 95.2  & -2.0760  & 2.0444 & -2.3716  & 2.4169 & 0.3435 & 0.3002 & 1.1442 \\
    \midrule
    $h_{3} = 0.03$ & 95.2  & 95.2  & -2.1586  & 2.1179 & -2.5243  & 2.5791 & 0.295 & 0.2888 & 1.0215 \\
    \midrule
    $h_{4} = 0.05$ & 95.2  & 95.2  & -2.1865  & 2.2459 & -2.7060  & 2.7412 & 0.2464 & 0.2746 & 0.8973 \\
    \midrule
    $h_{opt} = 0.0441$ & 95.2  & 95    & -2.1899  & 2.1468 & -2.7222  & 2.6919 & 0.2493 & 0.278 & 0.8968 \\
    \midrule
          & \multicolumn{2}{c}{mean} & -2.1515 & 2.0863 & -2.5225 & 2.5429 & 0.3209 & 0.2974 & 1.0791 \\
\cmidrule{2-10}          &       &       &       &       &       &       &       &       &  \\
    \midrule
    \multicolumn{10}{c}{Arrival intensity $\lambda \cdot T = 50$, Jump size $Z_{n} \sim \mathscr{N}(0,0.036^{2})$} \\
    \midrule
    $h_{1} = 0.01$ & 95.2  & 95.2  & -1.9490  & 2.0748 & -2.2026  & 2.2388 & 0.4721 & 0.3368 & 1.4017 \\
    \midrule
    $h_{2} = 0.02$ & 95.2  & 95.2  & -2.0882  & 2.074 & -2.4857  & 2.5502 & 0.3484 & 0.3173 & 1.098 \\
    \midrule
    $h_{3} = 0.03$ & 95.2  & 95.2  & -2.0592  & 2.1506 & -2.6234  & 2.5998 & 0.2917 & 0.2883 & 1.0118 \\
    \midrule
    $h_{4} = 0.05$ & 95.2  & 95.2  & -2.3158  & 2.0811 & -2.7272  & 2.7163 & 0.2455 & 0.2758 & 0.8901 \\
    \midrule
    $h_{opt} = 0.0441$ & 95.2  & 94.8  & -2.2977  & 2.0508 & -2.7281  & 2.6843 & 0.2508 & 0.279 & 0.8989 \\
    \midrule
          & \multicolumn{2}{c}{mean} & -2.142 & 2.0863 & -2.5534 & 2.5579 & 0.3217 & 0.2994 & 1.0743 \\
\cmidrule{2-10}          &       &       &       &       &       &       &       &       &  \\
    \midrule
    \multicolumn{10}{c}{Arrival intensity $\lambda \cdot T = 20$, Jump size $Z_{n} \sim \mathscr{N}(0,0.1^{2})$} \\
    \midrule
    $h_{1} = 0.01$ & 95.2  & 95.2  & -2.0544  & 1.8202 & -2.1361  & 2.1092 & 0.4599 & 0.3258 & 1.4116 \\
    \midrule
    $h_{2} = 0.02$ & 95.2  & 95.2  & -1.8318  & 1.8576 & -2.4108  & 2.2472 & 0.3124 & 0.2973 & 1.0508 \\
    \midrule
    $h_{3} = 0.03$ & 95.2  & 95.2  & -2.0650  & 1.9779 & -2.5856  & 2.4245 & 0.2834 & 0.2811 & 1.0081 \\
    \midrule
    $h_{4} = 0.05$ & 95.2  & 95.2  & -2.1245  & 2.086 & -2.8735  & 2.5503 & 0.2381 & 0.2786 & 0.8546 \\
    \midrule
    $h_{opt} = 0.0441$ & 95.2  & 95.2  & -2.1531  & 2.0896 & -2.8347  & 2.5017 & 0.2468 & 0.2783 & 0.8868 \\
    \midrule
          & \multicolumn{2}{c}{mean} & -2.0458 & 1.9663 & -2.5681 & 2.3666 & 0.3081 & 0.2922 & 1.0544 \\
\cmidrule{2-10}    \end{tabular}%
  \label{tab:addlabel}%
\end{sidewaystable}
\clearpage

\begin{sidewaystable}[!htb]
  \centering
  \caption{Adjusted length of confidence band for $\mu(0.15)$ over 500 replicates.}
    \begin{tabular}{cccccccccc}
    \toprule
    \multirow{2}[4]{*}{Bandwidth} & \multicolumn{2}{c}{Coverage Rate} & \multicolumn{2}{c}{Sym Quantile} & \multicolumn{2}{c}{Asym Quantile} & \multicolumn{3}{c}{Length of Confidence Band} \\
\cmidrule{2-10}          & Sym   & Asym  & 2.50\% & 97.50\% & 2.50\% & 97.50\% & Sym   & Asym  & Ratio \\
    \midrule
    \multicolumn{10}{c}{Arrival intensity $\lambda \cdot T = 20$, Jump size $Z_{n} \sim \mathscr{N}(0,0.036^{2})$} \\
    \midrule
    $h_{1} = 0.01$ & 95.2  & 95.2  & -2.0012  & 2.0126 & -1.9856  & 2.1802 & 0.4761 & 0.2695 & 1.7666 \\
    \midrule
    $h_{2} = 0.02$ & 95.2  & 95.2  & -1.9823  & 2.0154 & -1.9710  & 2.1613 & 0.3381 & 0.237 & 1.4266 \\
    \midrule
    $h_{3} = 0.03$ & 95.2  & 95.2  & -2.0890  & 2.1366 & -1.9350  & 2.3658 & 0.2956 & 0.2325 & 1.2714 \\
    \midrule
    $h_{4} = 0.05$ & 95.2  & 95.2  & -2.0606  & 2.2294 & -1.9769  & 2.3788 & 0.2415 & 0.2214 & 1.0908 \\
    \midrule
    $h_{opt} = 0.0441$ & 95.2  & 94.8  & -2.0697  & 2.2194 & -1.9547  & 2.3841 & 0.2498 & 0.2226 & 1.1222 \\
    \midrule
          & \multicolumn{2}{c}{mean} & -2.0406 & 2.1227 & -1.9646 & 2.294 & 0.3202 & 0.2366 & 1.3534 \\
\cmidrule{2-10}          &       &       &       &       &       &       &       &       &  \\
    \midrule
    \multicolumn{10}{c}{Arrival intensity $\lambda \cdot T = 50$, Jump size $Z_{n} \sim \mathscr{N}(0,0.036^{2})$} \\
    \midrule
    $h_{1} = 0.01$ & 95.2  & 95.2  & -2.1808  & 2.0987 & -1.9418  & 2.1174 & 0.51  & 0.2638 & 1.9333 \\
    \midrule
    $h_{2} = 0.02$ & 95.2  & 95.2  & -1.9529  & 2.3189 & -2.1120  & 2.0546 & 0.363 & 0.2401 & 1.5119 \\
    \midrule
    $h_{3} = 0.03$ & 95.2  & 95.2  & -1.9840  & 2.3023 & -2.1067  & 2.0558 & 0.3013 & 0.2261 & 1.3326 \\
    \midrule
    $h_{4} = 0.05$ & 95.2  & 95.2  & -2.1573  & 2.1565 & -2.2084  & 2.101 & 0.244 & 0.2201 & 1.1086 \\
    \midrule
    $h_{opt} = 0.0441$ & 95.2  & 95    & -2.1034  & 2.1938 & -2.2020  & 2.1022 & 0.2511 & 0.2218 & 1.1321 \\
    \midrule
          & \multicolumn{2}{c}{mean} & -2.0757 & 2.214 & -2.1142 & 2.0862 & 0.3339 & 0.2344 & 1.4245 \\
\cmidrule{2-10}          &       &       &       &       &       &       &       &       &  \\
    \midrule
    \multicolumn{10}{c}{Arrival intensity $\lambda \cdot T = 20$, Jump size $Z_{n} \sim \mathscr{N}(0,0.1^{2})$} \\
    \midrule
    $h_{1} = 0.01$ & 95.2  & 95.2  & -2.0582  & 2.0251 & -2.0616  & 2.2115 & 0.494 & 0.2821 & 1.7512 \\
    \midrule
    $h_{2} = 0.02$ & 95.2  & 95.2  & -2.0097  & 2.1434 & -2.1454  & 2.2018 & 0.3585 & 0.2544 & 1.4092 \\
    \midrule
    $h_{3} = 0.03$ & 95.2  & 95.2  & -2.0967  & 2.1874 & -2.1593  & 2.2739 & 0.3059 & 0.2445 & 1.2511 \\
    \midrule
    $h_{4} = 0.05$ & 95.2  & 95.2  & -2.1265  & 2.2835 & -2.1661  & 2.3895 & 0.2534 & 0.2363 & 1.0724 \\
    \midrule
    $h_{opt} = 0.0441$ & 95.2  & 95.2  & -2.0901  & 2.3064 & -2.1598  & 2.3709 & 0.2596 & 0.2367 & 1.0967 \\
    \midrule
          & \multicolumn{2}{c}{mean} & -2.0762 & 2.1892 & -2.1384 & 2.2895 & 0.3343 & 0.2508 & 1.3329 \\
\cmidrule{2-10}    \end{tabular}%
  \label{tab:addlabel}%
\end{sidewaystable}
\clearpage

\begin{sidewaystable}[!htb]
  \centering
  \caption{ Estimation for $M(0.15)$ over 500 replicates.}
    \begin{tabular}{cccccccccccccc}
    \toprule
    \multirow{2}[4]{*}{Bandwidth} & \multicolumn{2}{c}{Coverage Rate} & \multicolumn{2}{c}{Mean of Bias} & \multicolumn{2}{c}{Variance ($\times 10^{-4}$)} & \multicolumn{4}{c}{Estimation for Variance ($\times 10^{-4}$)} & \multicolumn{3}{c}{Length of Confidence Band} \\
\cmidrule{2-14}          & Sym   & Asym  & Sym   & Asym  & Sym   & Asym  & Sym   & Std   & Asym  & Std   & Sym   & Asym  & Ratio \\
    \midrule
    \multicolumn{14}{c}{Arrival intensity $\lambda \cdot T = 20$, Jump size $Z_{n} \sim \mathscr{N}(0,0.036^{2})$} \\
    \midrule
    $h_{1} = 0.01$ & 94.8  & 94.8  & 0.004 & 0.0049 & 0.3221 & 0.0991 & 0.4581 & 0.0679 & 0.1369 & 0.0125 & 0.0265 & 0.0145 & 1.8276  \\
    \midrule
    $h_{2} = 0.02$ & 92    & 90.4  & 0.0043 & 0.0053 & 0.1669 & 0.0821 & 0.233 & 0.0268 & 0.108 & 0.0087 & 0.0189 & 0.0129 & 1.4651  \\
    \midrule
    $h_{3} = 0.04$ & 76.4  & 79.2  & 0.0048 & 0.0058 & 0.0985 & 0.0731 & 0.1242 & 0.0108 & 0.0895 & 0.0065 & 0.0138 & 0.0117 & 1.1795  \\
    \midrule
    $h_{4} = 0.05$ & 68.6  & 68.6  & 0.0051 & 0.0059 & 0.0867 & 0.0714 & 0.1037 & 0.0082 & 0.0852 & 0.006 & 0.0126 & 0.0114 & 1.1053  \\
    \midrule
    $h_{opt} = 0.0294$ & 85    & 86.2  & 0.0045 & 0.0056 & 0.1237 & 0.0766 & 0.1647 & 0.0162 & 0.097 & 0.0073 & 0.0159 & 0.0122 & 1.3033  \\
    \midrule
    \multicolumn{14}{c}{} \\
    \midrule
    \multicolumn{14}{c}{Arrival intensity $\lambda \cdot T = 50$, Jump size $Z_{n} \sim \mathscr{N}(0,0.036^{2})$} \\
    \midrule
    $h_{1} = 0.01$ & 94.6  & 93    & 0.0034 & 0.0045 & 0.329 & 0.1114 & 0.4626 & 0.0692 & 0.1392 & 0.0132 & 0.0266 & 0.0146 & 1.8219  \\
    \midrule
    $h_{2} = 0.02$ & 91    & 89.8  & 0.0038 & 0.005 & 0.1807 & 0.0911 & 0.2361 & 0.0276 & 0.11  & 0.0094 & 0.019 & 0.013 & 1.4615  \\
    \midrule
    $h_{3} = 0.04$ & 79.8  & 82.6  & 0.0045 & 0.0056 & 0.109 & 0.0803 & 0.1263 & 0.0116 & 0.0913 & 0.0071 & 0.0139 & 0.0118 & 1.1780  \\
    \midrule
    $h_{4} = 0.05$ & 72.2  & 72.6  & 0.0048 & 0.0057 & 0.096 & 0.0784 & 0.1056 & 0.0089 & 0.087 & 0.0066 & 0.0127 & 0.0116 & 1.0948  \\
    \midrule
    $h_{opt} = 0.0294$ & 87    & 85.4  & 0.0041 & 0.0053 & 0.1356 & 0.0843 & 0.1662 & 0.0168 & 0.0987 & 0.0079 & 0.016 & 0.0123 & 1.3008  \\
    \midrule
    \multicolumn{14}{c}{} \\
    \midrule
    \multicolumn{14}{c}{Arrival intensity $\lambda \cdot T = 20$, Jump size $Z_{n} \sim \mathscr{N}(0,0.1^{2})$} \\
    \midrule
    $h_{1} = 0.01$ & 95    & 93.8  & 0.0038 & 0.0046 & 0.5407 & 0.1829 & 0.6642 & 0.3883 & 0.1961 & 0.0698 & 0.0311 & 0.0171 & 1.8187  \\
    \midrule
    $h_{2} = 0.02$ & 93    & 91.2  & 0.004 & 0.0052 & 0.3149 & 0.1456 & 0.3372 & 0.1509 & 0.154 & 0.0487 & 0.0224 & 0.0152 & 1.4737  \\
    \midrule
    $h_{3} = 0.04$ & 85.2  & 85.2  & 0.0046 & 0.0057 & 0.179 & 0.1276 & 0.1778 & 0.0614 & 0.1272 & 0.0371 & 0.0163 & 0.0139 & 1.1727  \\
    \midrule
    $h_{4} = 0.05$ & 78.2  & 79.2  & 0.0049 & 0.0058 & 0.1544 & 0.1246 & 0.1481 & 0.0471 & 0.1211 & 0.0347 & 0.0149 & 0.0135 & 1.1037  \\
    \midrule
    $h_{opt} = 0.0294$ & 88.6  & 87.8  & 0.0043 & 0.0055 & 0.2285 & 0.1339 & 0.2323 & 0.0889 & 0.1373 & 0.041 & 0.0186 & 0.0144 & 1.2917  \\
    \bottomrule
    \end{tabular}%
  \label{tab:addlabel}%
\end{sidewaystable}
\clearpage

\section{ Empirical Analysis }

In this section, we apply the second-order jump-diffusion to model
the return of stock index in Shanghai Stock Exchange between July
2014 and Dec 2014 from China under five-minute high frequency data,
and then apply the local linear estimators to estimate the unknown
coefficients in model (\ref{1.2}) based on Gamma asymmetric kernels
and Gaussian symmetric kernels. The real-world financial market data
set analyzed consists of 6048 observations.

We assume that
\begin{equation}
\label{ea1} \left\{
\begin{array}{ll}
d \log Y_{t} = X_{t}dt,\\
dX_{t} = \mu(X_{t-})dt + \sigma(X_{t-})dW_{t}+
\int_{\mathscr{E}}c(X_{t-} , z) r(\omega, dt, dz),
\end{array}
\right.\end{equation} where $\log Y_{t}$ is the log integrated
process for stock index and $X_{t}$ is the latent process for the
log-returns. According to (\ref{2.2}), we can get the proxy of the
latent process
\begin{equation}
\label{ea2} \widetilde{X}_{i\Delta_n}=\frac{\log Y_{i\Delta_n} -
\log Y_{(i-1)\Delta_n}}{\Delta_n}.
\end{equation}

The plots of the stock index and its proxy (\ref{ea2}) from China,
i.e. Shanghai Composite Index in high frequency data are shown in
Figure 12.

\begin{figure}[!htb]
\centering
  \subfigure[ Shanghai Composite Index (2014) ]{
  \label{fig:subfig} 
    \includegraphics[width=0.45\textwidth]{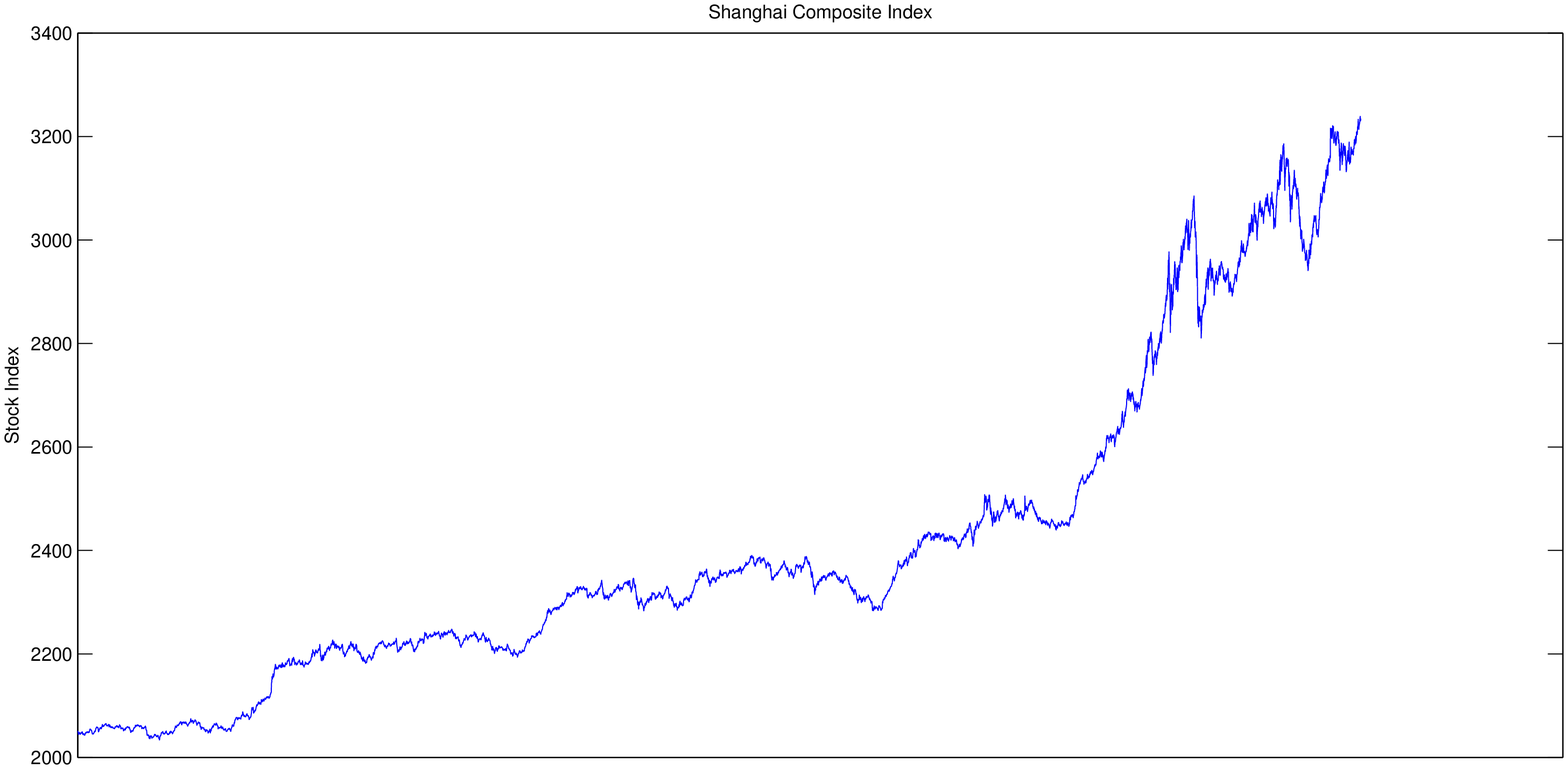}}
 \hspace{0.25in}
  \subfigure[ Proxy of Shanghai Composite Index (2014) ]{
  \label{fig:subfig} 
    \includegraphics[width=0.45\textwidth]{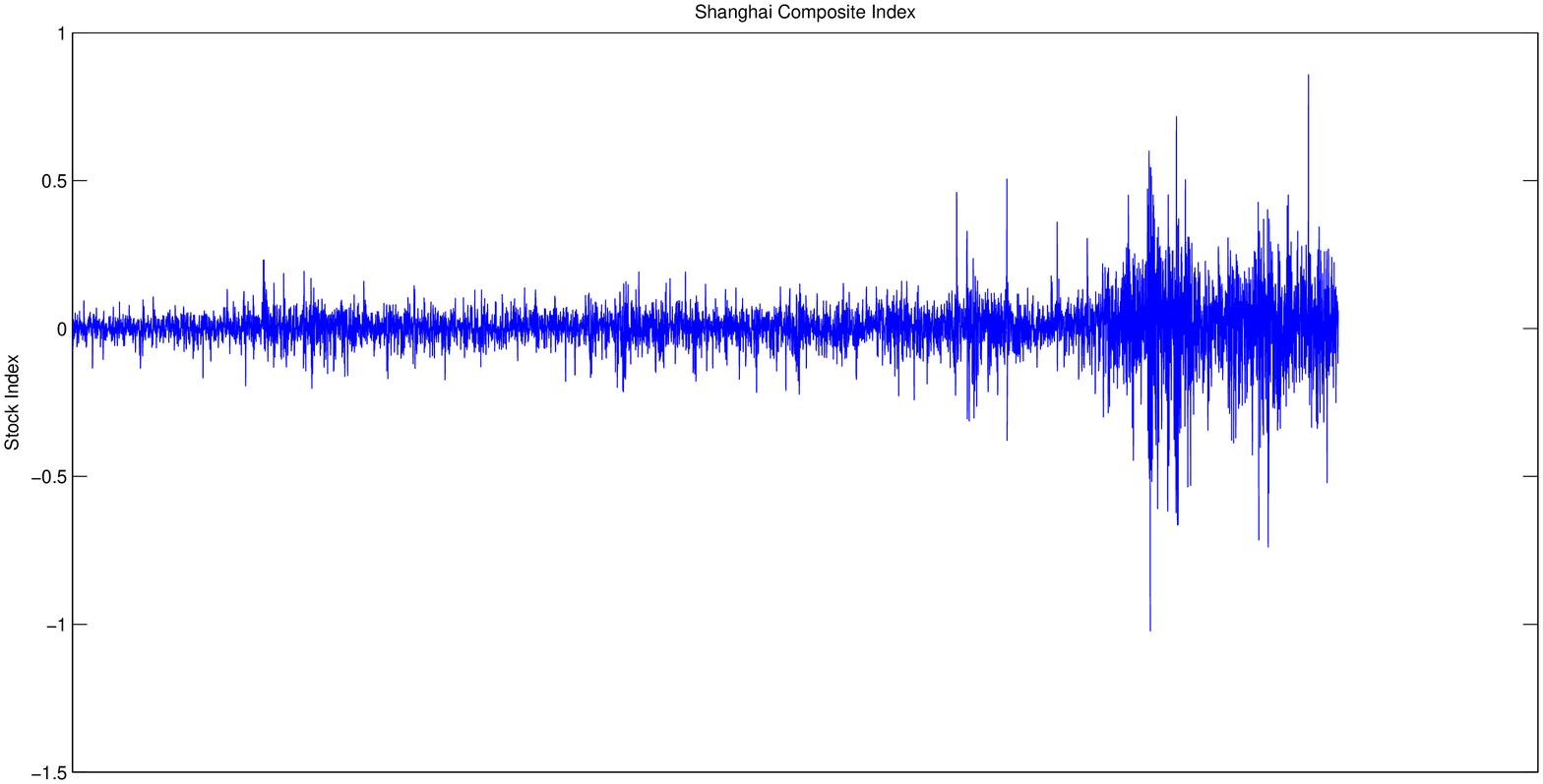}}
  \caption{ Time Series and Proxy of Shanghai Composite Index (2014) from July 01, 2014 to Dec 31, 2014 }
  \label{fig:subfig} 
\end{figure}

First, we test the existence of jumps for the proxy $X_{t}$ through
the test statistic proposed in Barndorff-Nielsen and Shephard
\cite{bs2} (denoted by BS Statistic). For five-minute high frequency
data, the value of BS Statistic is -3.9955, which exceeds [-1.96,
1.96], so there exists jumps in high frequency data at the 5\%
significance level, which confirms the validity of model (\ref{1.2})
not model (\ref{nm}) for the return of stock index by the
second-order process. Based on the Augmented Dickey-Fuller test
statistic, we can easily get that the null hypothesis of
non-stationarity is accepted at the 5\% significance level for the
stock index $Y_{t}$, but is rejected for the proxy of $X_{t}$, which
confirms the assumption of stationary by differencing.

Here we use two alternative smoothing parameters $h_{cv}$ which is
selected by the $k-$block cross-validation method, and $h_{T} = c
\cdot \hat{S} \cdot T^{-\frac{2}{5}}$ with $c = 4$ for drift and $c
= 2$ for volatility (chosen only for illustration). Figure 13
depicts the curves of CV($h_{cv}$) versus $h_{cv}$ for Shanghai
Composite Index showing that CV($h_{cv}$) is minimized at $h_{cv} =
0.095$ for drift estimator, $h_{cv} = 0.040$ for volatility
estimator, respectively.

\begin{figure}[!htb]
\centering
  \subfigure[ $CV(h_{n})$ for drift estimator ]{
  \label{fig:subfig} 
    \includegraphics[width=0.45\textwidth]{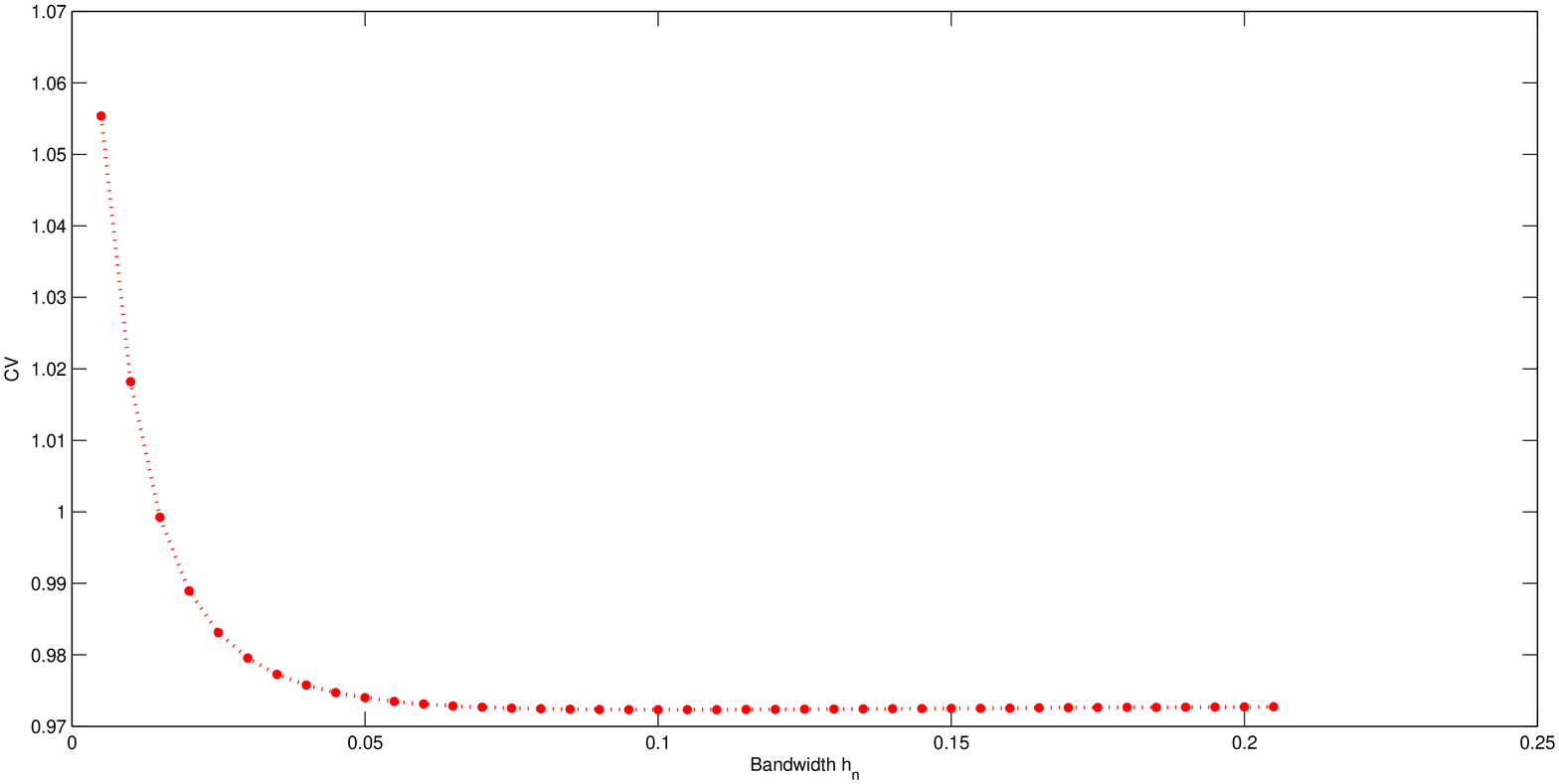}}
 \hspace{0.25in}
  \subfigure[ $CV(h_{n})$ for volatility estimator ]{
  \label{fig:subfig} 
    \includegraphics[width=0.45\textwidth]{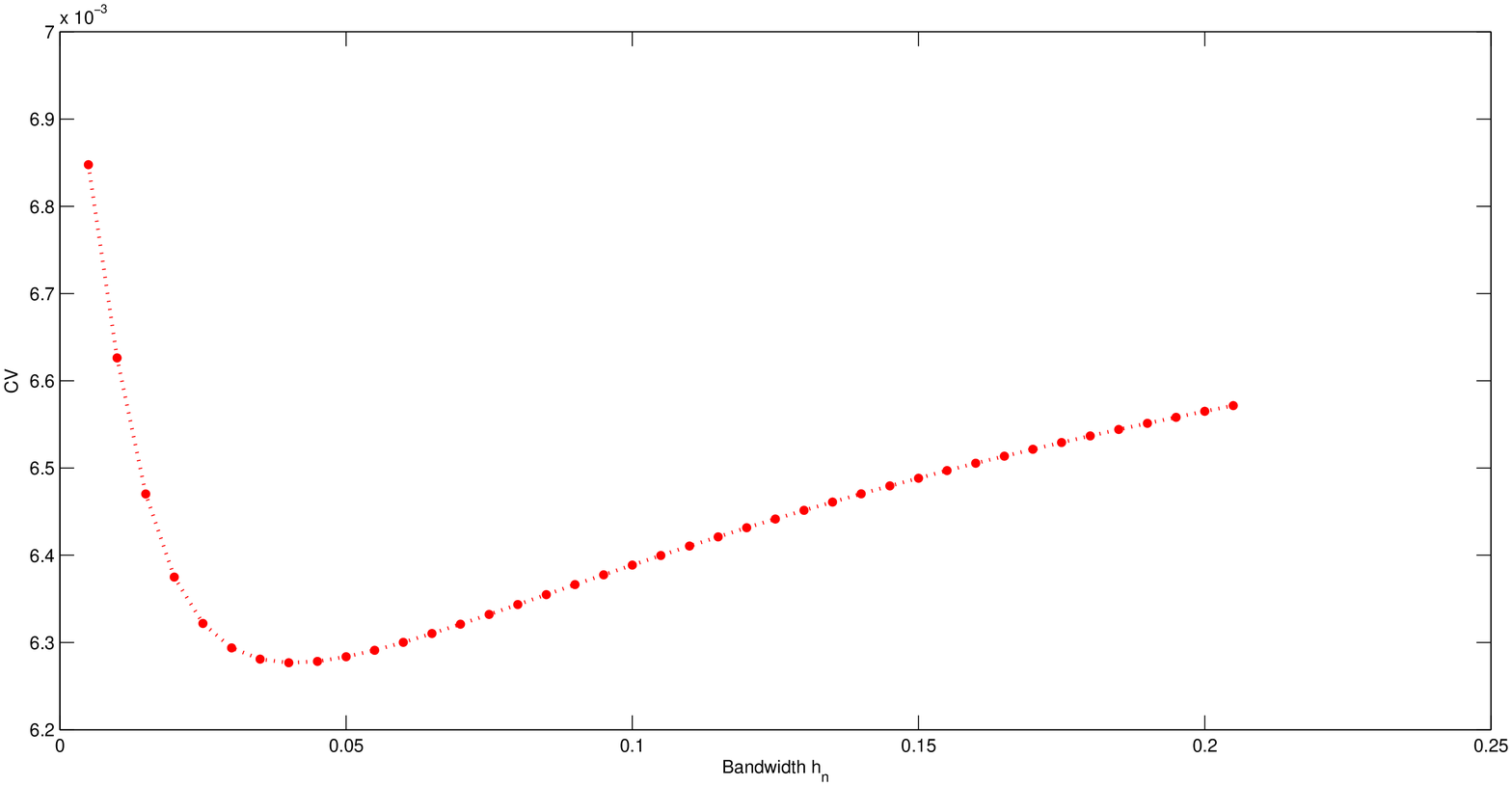}}
  \caption{ Curve of $CV(h_{n})$ versus $h_{n}$ for Shanghai Composite Index (2014) }
  \label{fig:subfig} 
\end{figure}

Then, we will employ the local linear estimators based on Gamma
asymmetric kernels (\ref{2.7}) and (\ref{2.8}) to estimate the
unknown coefficients under (\ref{ea2}) and $\Delta_{n} =
\frac{1}{48}$ for five-minute data (t = 1 meaning one day) with
various bandwidth such as $h_{cv}$ and $h_{T}.$ The estimation
curves for unknown qualities in five-minute high frequency data are
displayed in Figure 14.

It is observed that the linear shape with negative coefficient for
drift estimator in FIG 14 (a) \& (b) for various bandwidths which
indicates that the higher log-return increments correspond to the
lower drift in the latent process (this fact coincides with the
economic phenomenon of mean reversion), which reveals a negative
correlation. It is also shown the quadratic form with positive
coefficient for volatility estimator with a minimum at 0.3 in FIG 14
(c) \& (d) which reveals that the higher absolute value of
log-return increments correspond to the higher volatility in the
latent process (this fact coincides with the economic phenomenon of
volatility smile). These findings are consistent with those in
Nicolau \cite{ni2}.

\begin{figure}[!htb]
\centering
  \subfigure[ Drift estimator with $h_{cv}$ ]{
  \label{fig:subfig} 
    \includegraphics[width=0.45\textwidth]{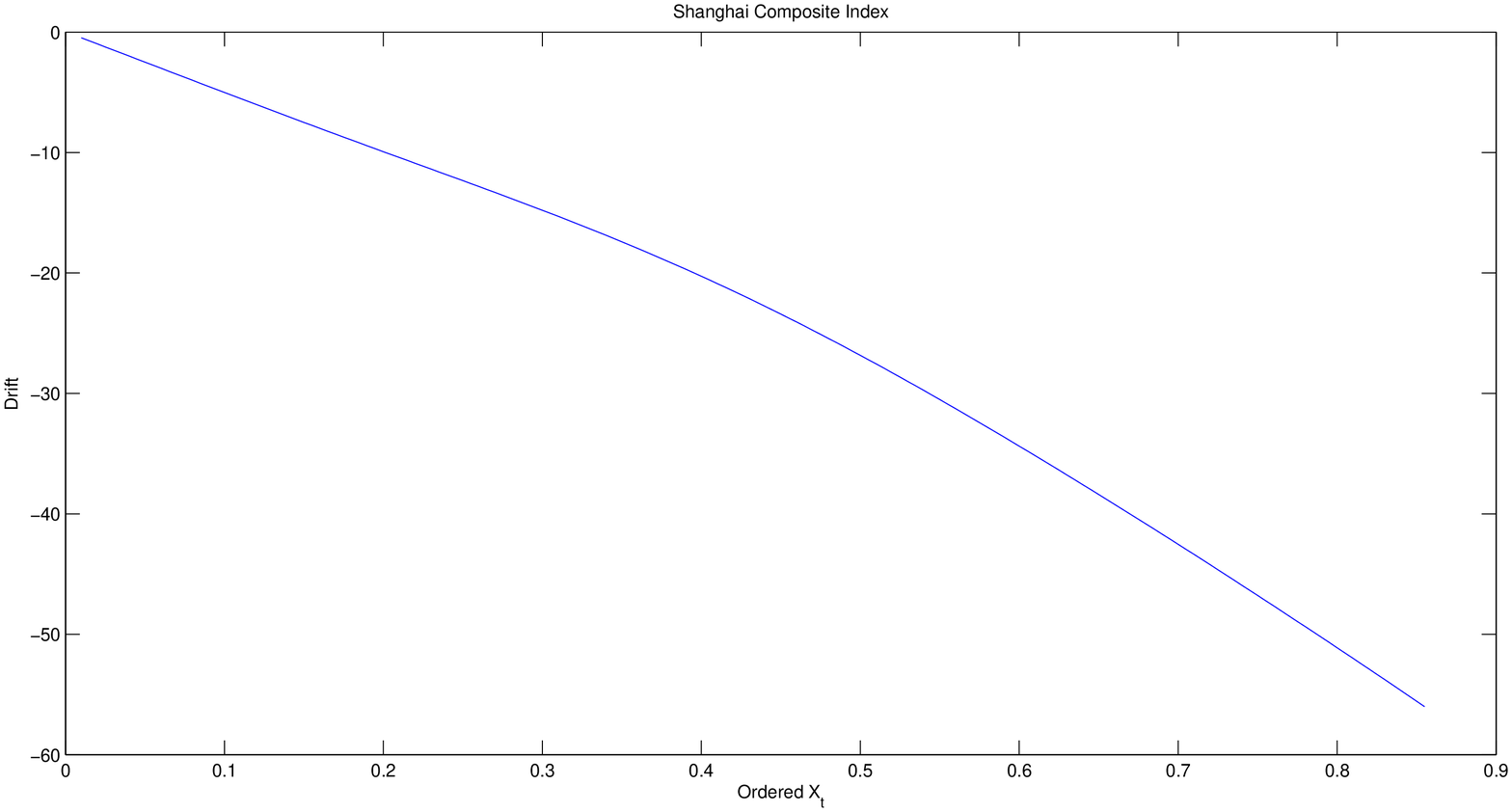}}
 \hspace{0.25in}
  \subfigure[ Drift estimator with $h_{T}$ ]{
  \label{fig:subfig} 
    \includegraphics[width=0.45\textwidth]{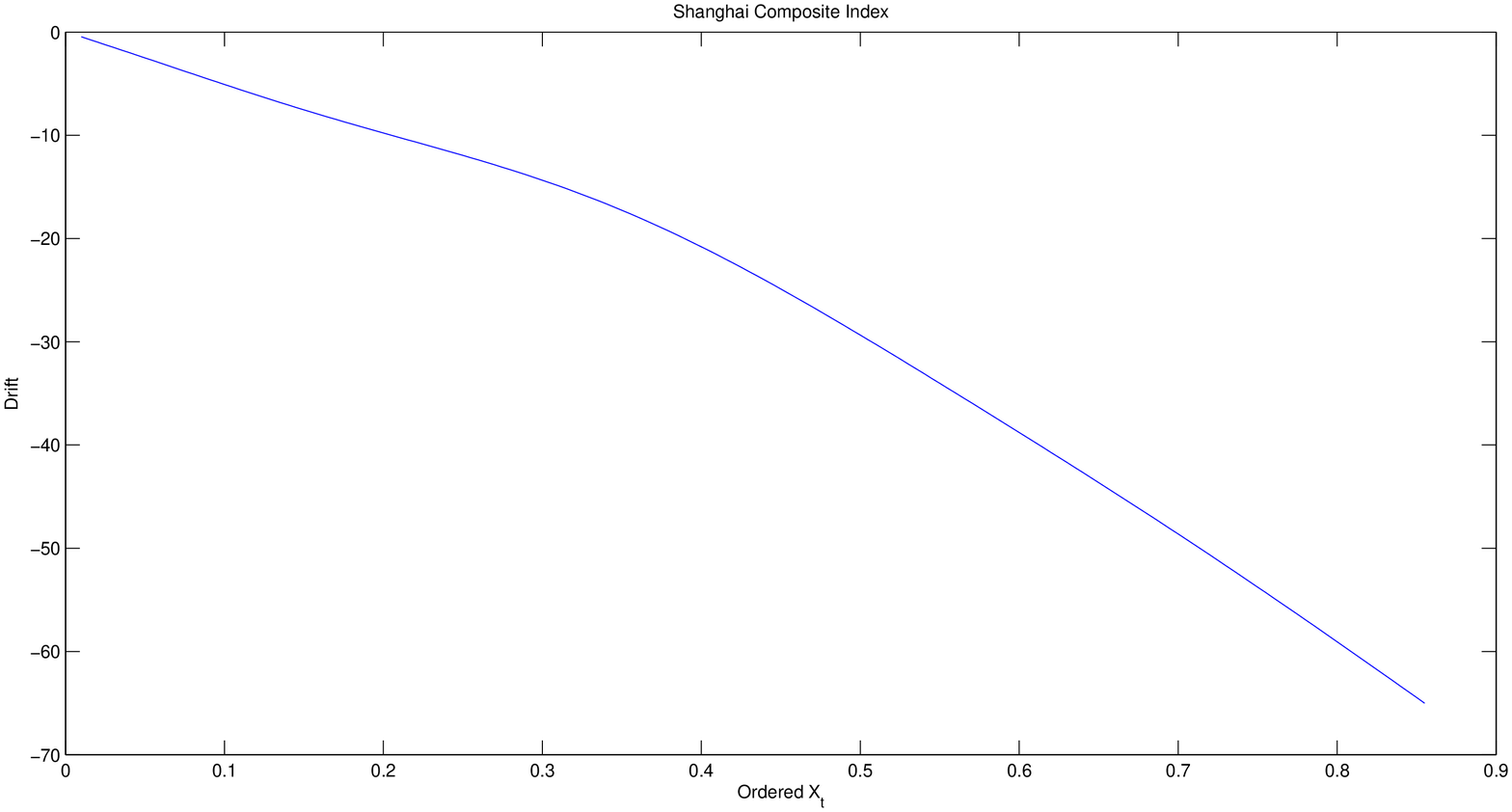}}
  \subfigure[ Volatility estimator with $h_{cv}$ ]{
  \label{fig:subfig} 
    \includegraphics[width=0.45\textwidth]{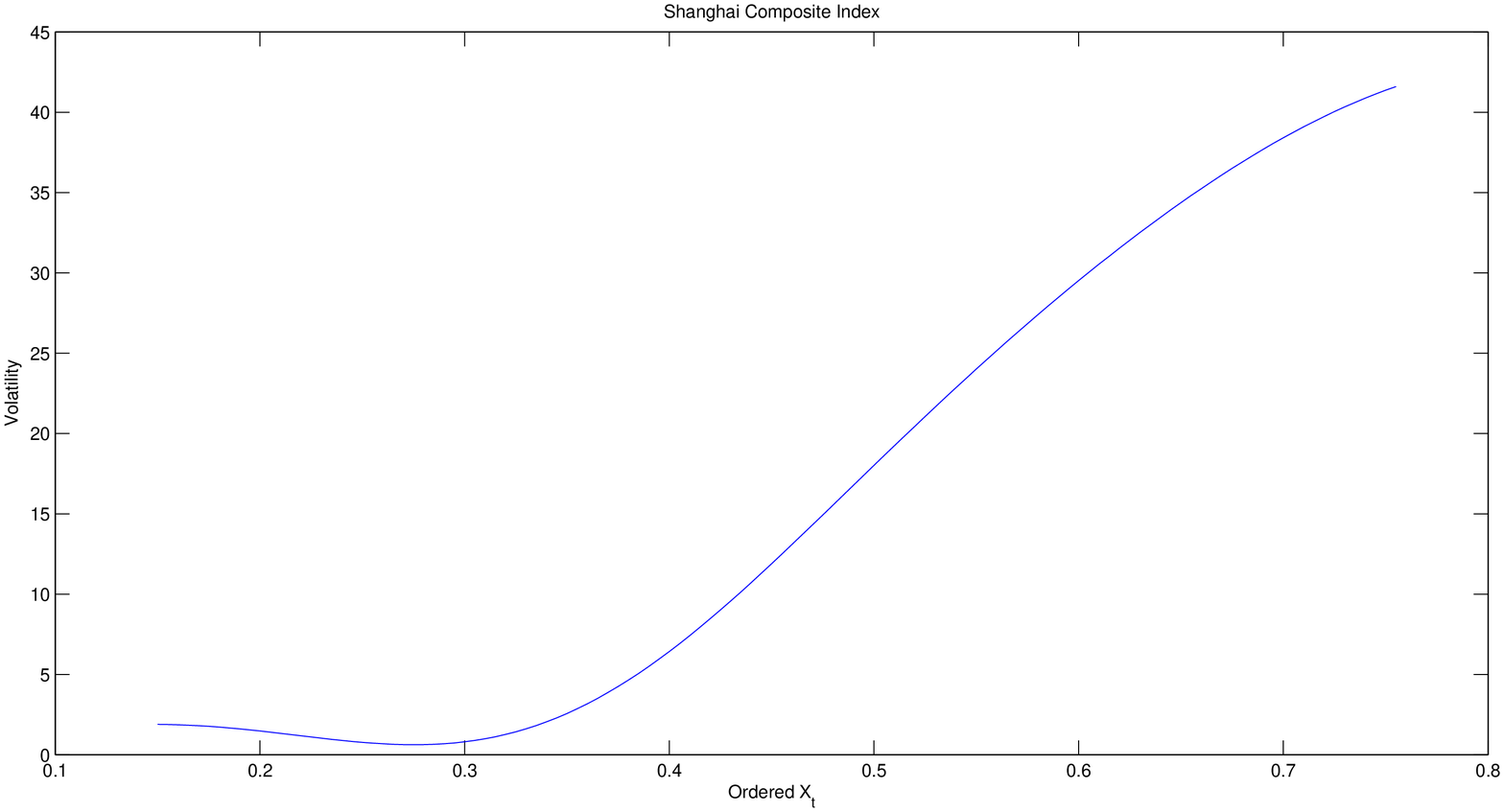}}
 \hspace{0.25in}
  \subfigure[ Volatility estimator with $h_{T}$ ]{
  \label{fig:subfig} 
    \includegraphics[width=0.45\textwidth]{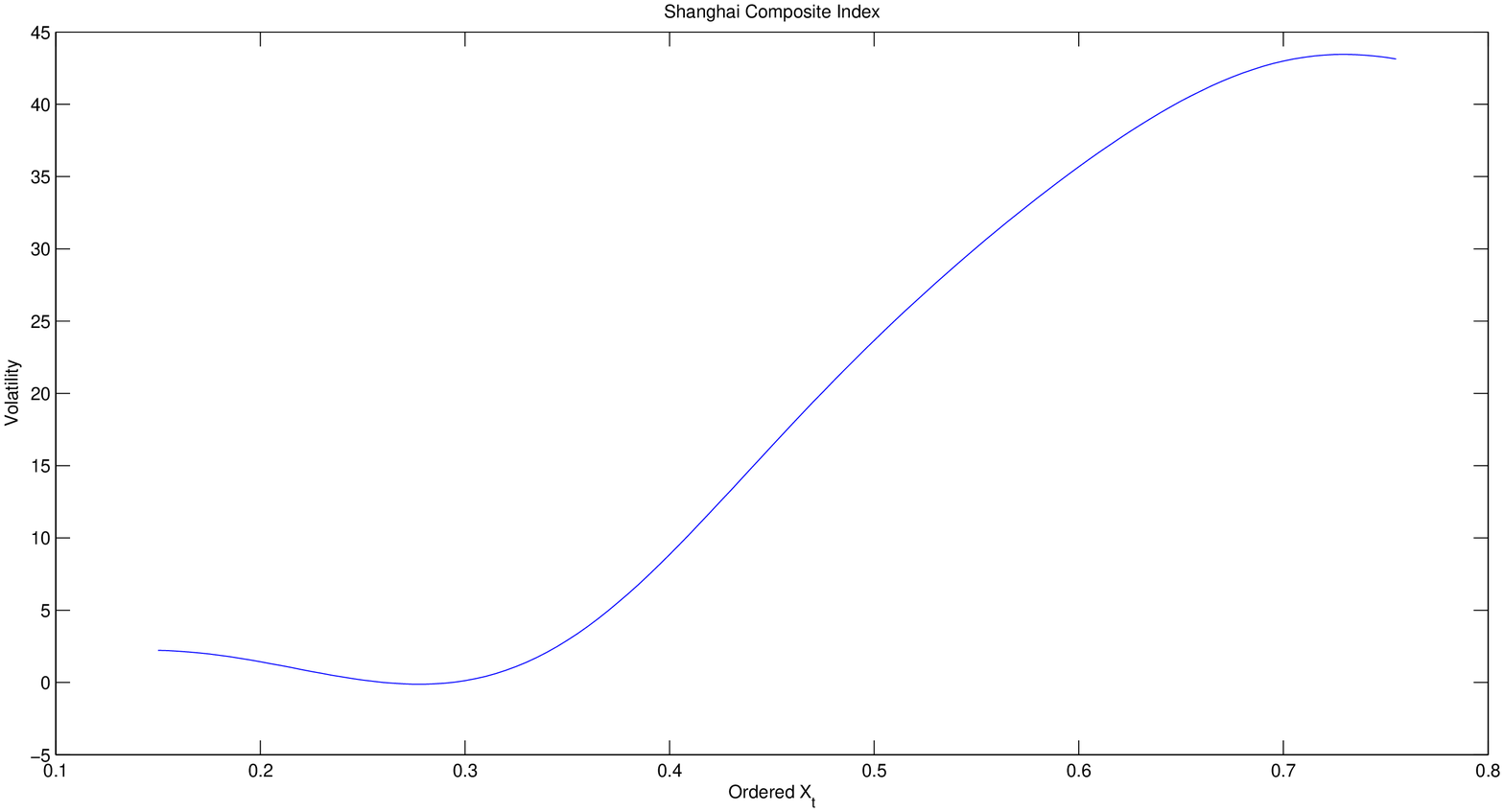}}
  \caption{ Local linear estimators of the drift and volatility coefficients for
  Shanghai Composite Index (2014) based on Gamma kernels using various bandwidths }
  \label{fig:subfig} 
\end{figure}

Finally, we will construct 95\% normal confidence intervals for the
unknown coefficients based on Gamma asymmetric kernels and Gaussian
symmetric kernels under (\ref{ea2}) and $\Delta_{n} = \frac{1}{48}$
for five-minute data (t = 1 meaning one day) with various bandwidth
such as $h_{cv}$ and $h_{T}.$ The 95\% normal confidence bands for
the drift and volatility functions are demonstrated in Figure 15.
All the quantities are computed at 120 equally spaced nonnegative
ordered $\widetilde{X}_{i\Delta_n}$ from 0.01 to 0.601. For a more
intuitive comparison of the lengths of normal confidence intervals
of the drift and volatility coefficients based on Gaussian kernels
and Gamma kernels for Shanghai Composite Index (2014) using various
bandwidths, here the ratios of the length of confidence band
constructed with Gaussian symmetric kernel (CB-GSK) to that
constructed with Gamma asymmetric kernel (CB-GAK) are shown in
Figure 16. Note that the blue dotted lines in Figure 16 represent
the ratio value of one.

From Figures 15 and 16, we can observe the following findings.
\begin{itemize}
\item { As for the quantities close to zero, the ratios are less than one,
which coincides with the discussion in Remark \ref{r3.6} that the
closer to the boundary point or the larger bandwidth for fixed
``boundary $x$'', the shorter the length of confidence interval
based on Gaussian symmetric kernel. }
\item {As the points increase,
especially at the sparse points, CB-GSK tends to be longer than
CB-GAK and the ratios gradually become larger and greater than 1,
which effectively verifies the efficiency gains and resistance to
sparse points of local linear smoothing using Gamma asymmetric
kernel through the real high frequency financial data. }
\item { Note that some values for the ratios in Figure 16 are zero, which is due to the
fact that the local linear estimators based on Gaussian symmetric
kernel for conditional variance and conditional fourth moment are
negative. Fortunately, the local linear estimators based on Gamma
asymmetric kernel for conditional variance and conditional fourth
moment are positive. }
\end{itemize}

\begin{figure}[!htb]
\centering
  \subfigure[ Drift confidence band with $h_{cv}$ ]{
  \label{fig:subfig} 
    \includegraphics[width=0.45\textwidth]{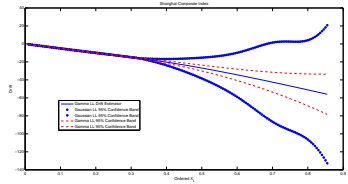}}
 \hspace{0.25in}
  \subfigure[ Drift confidence band with $h_{T}$ ]{
  \label{fig:subfig} 
    \includegraphics[width=0.45\textwidth]{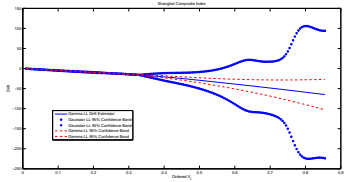}}
  \subfigure[ Volatility confidence band with $h_{cv}$ ]{
  \label{fig:subfig} 
    \includegraphics[width=0.45\textwidth]{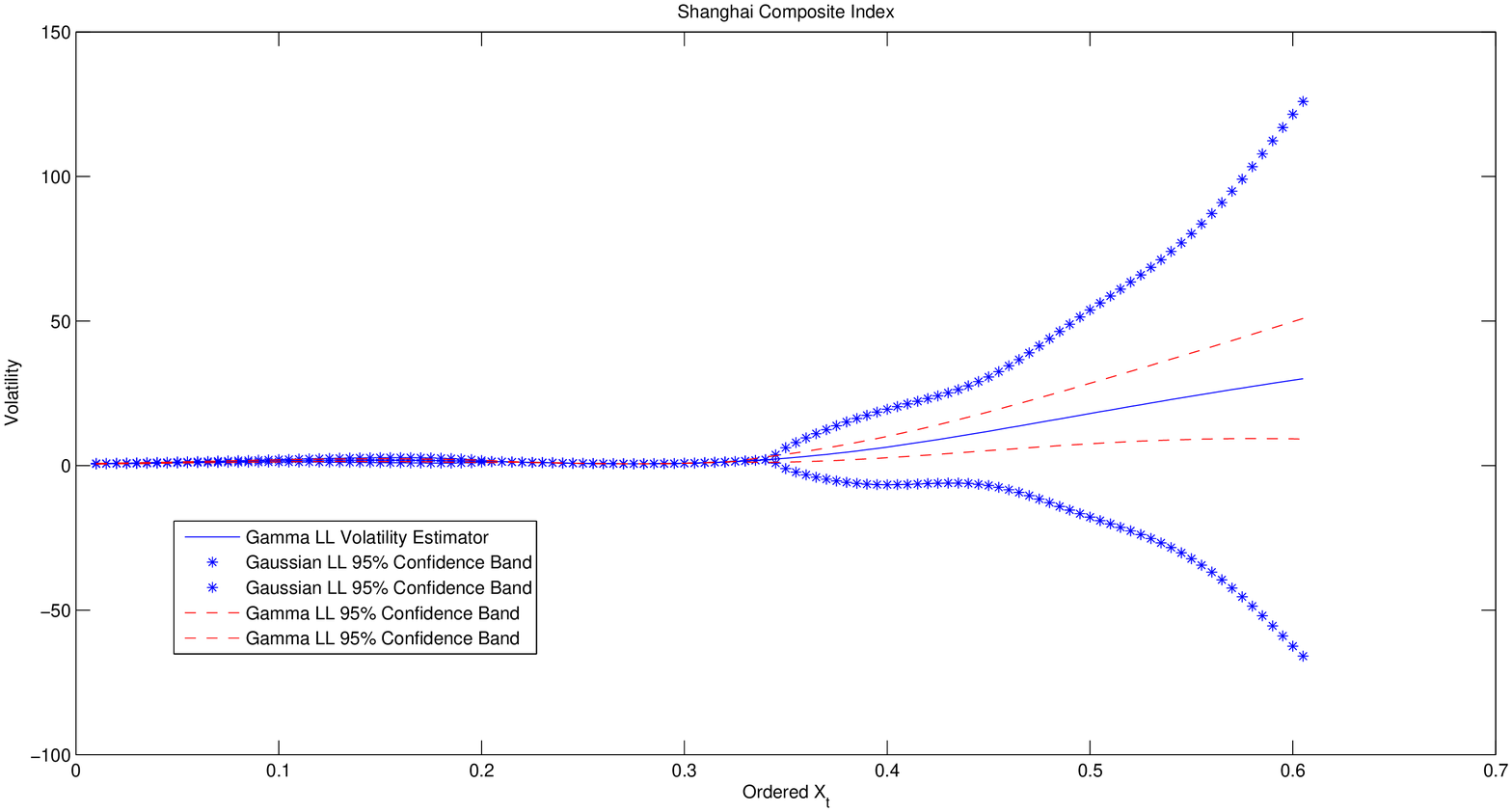}}
 \hspace{0.25in}
  \subfigure[ Volatility confidence band with $h_{T}$ ]{
  \label{fig:subfig} 
    \includegraphics[width=0.45\textwidth]{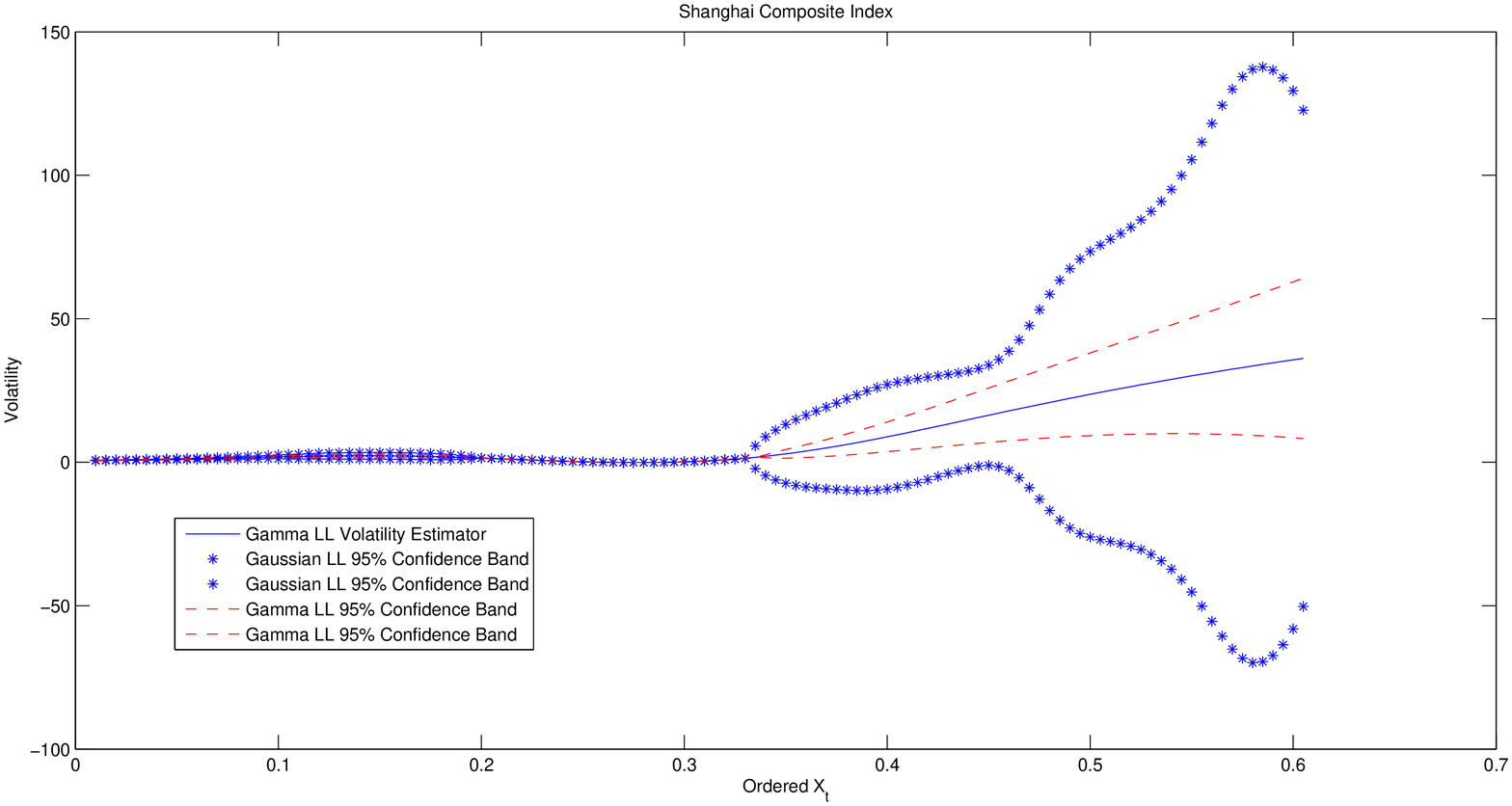}}
  \caption{ 95\% normal confidence intervals of the drift and volatility coefficients for
  Shanghai Composite Index (2014) based on Gamma kernels and Gaussian kernels using various bandwidths }
  \label{fig:subfig} 
\end{figure}

\begin{figure}[!htb]
\centering
  \subfigure[ Length Ratios for Drift with $h_{cv}$ ]{
  \label{fig:subfig} 
    \includegraphics[width=0.45\textwidth]{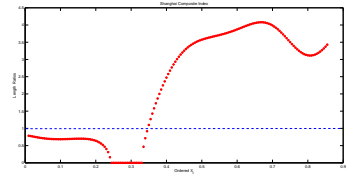}}
 \hspace{0.25in}
  \subfigure[ Length Ratios for Drift with $h_{T}$ ]{
  \label{fig:subfig} 
    \includegraphics[width=0.45\textwidth]{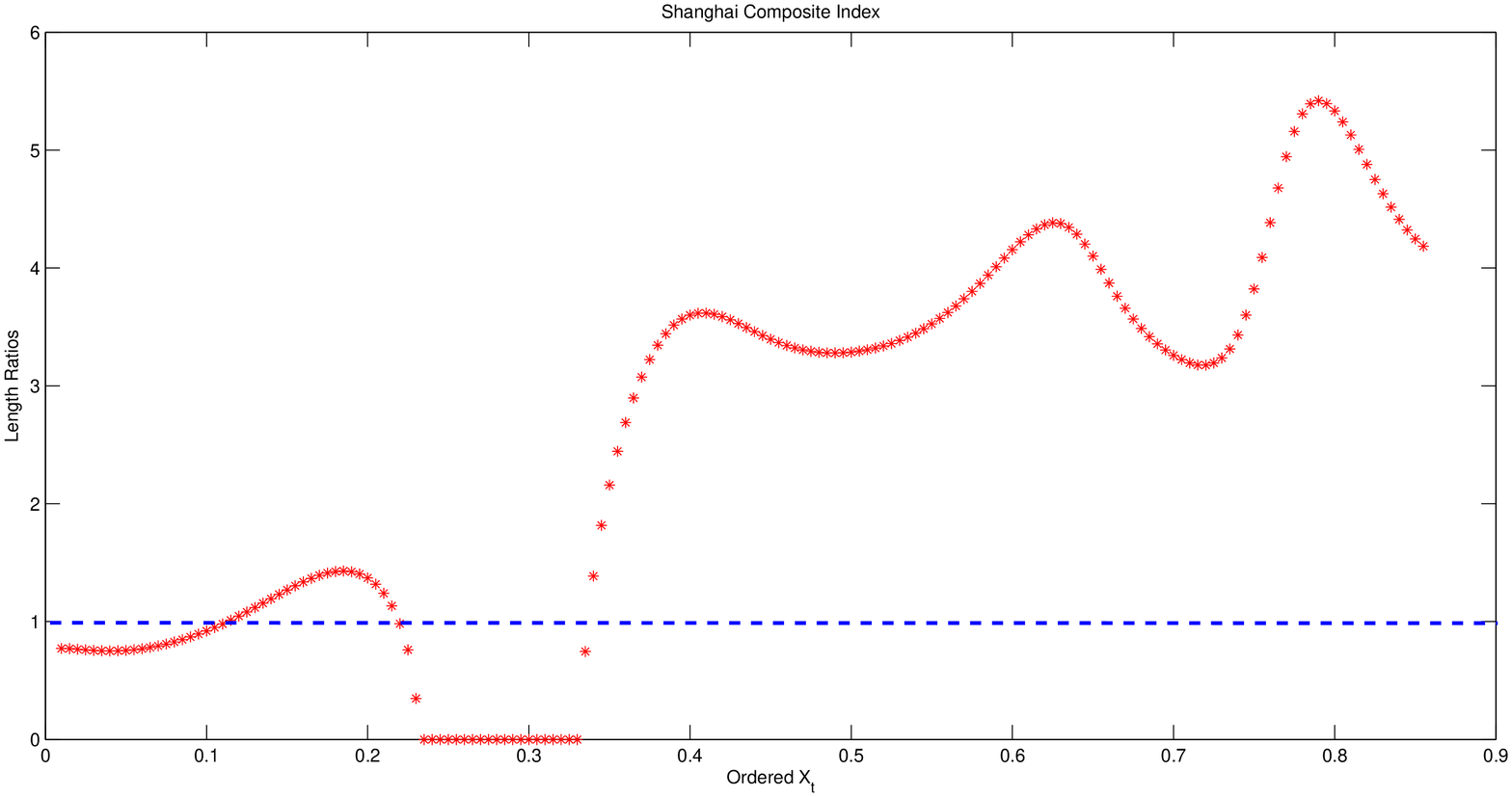}}
  \subfigure[ Length Ratios for Volatility with $h_{cv}$ ]{
  \label{fig:subfig} 
    \includegraphics[width=0.45\textwidth]{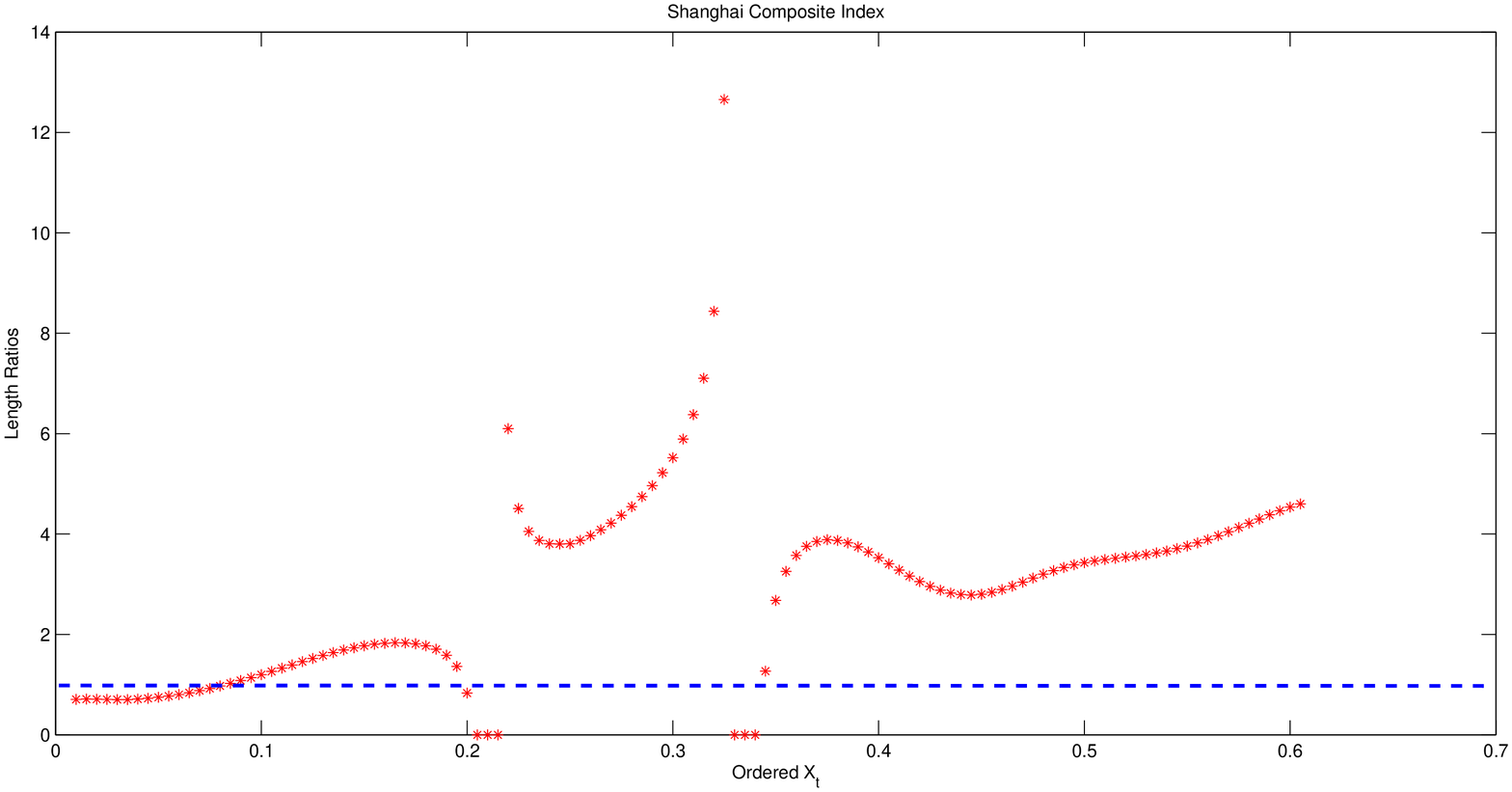}}
 \hspace{0.25in}
  \subfigure[ Length Ratios for Volatility with $h_{T}$ ]{
  \label{fig:subfig} 
    \includegraphics[width=0.45\textwidth]{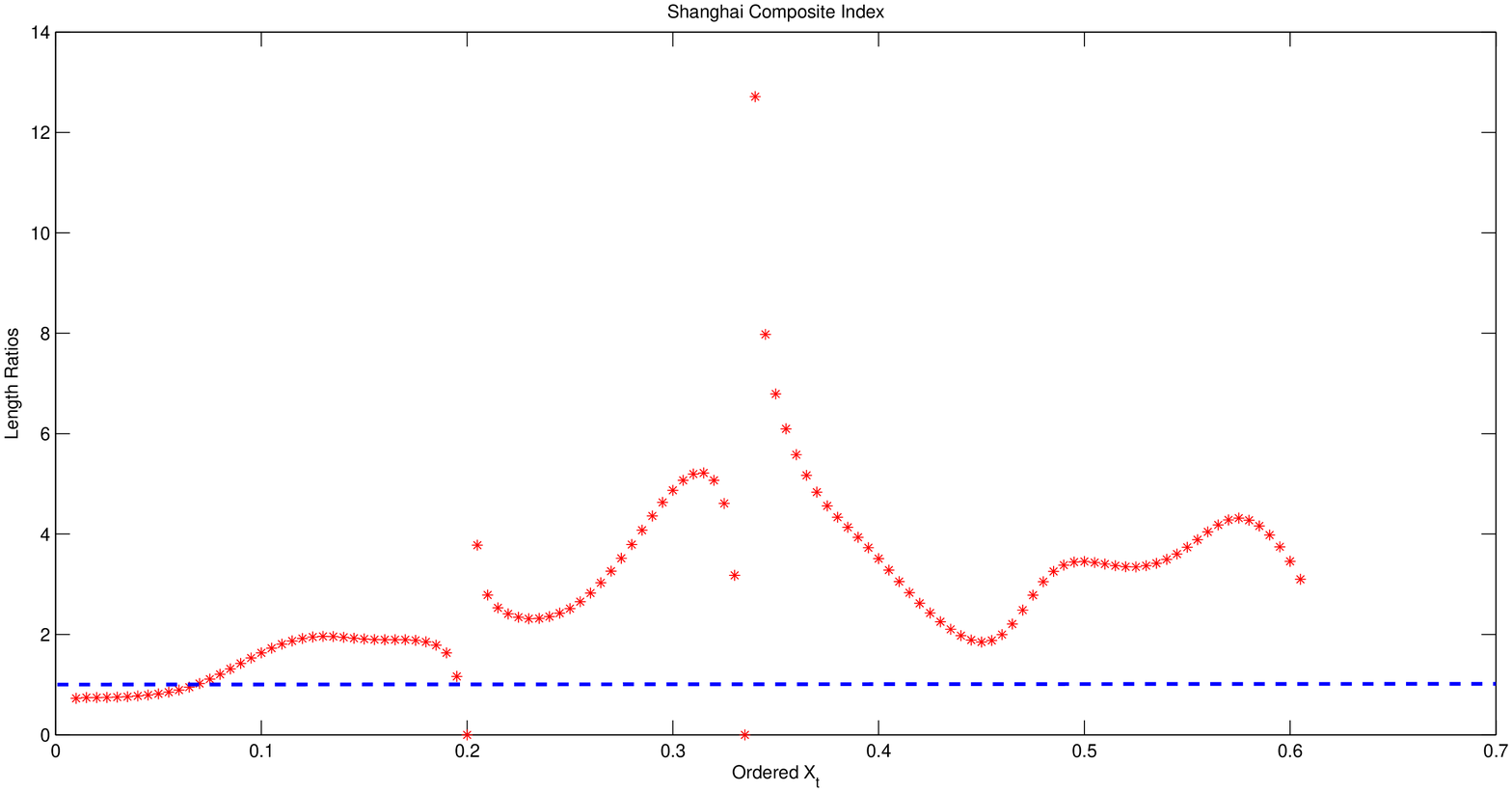}}
  \caption{ The ratios between the lengths of normal confidence intervals of the drift and volatility
  coefficients based on Gaussian kernels and Gamma kernels for
  Shanghai Composite Index (2014) using various bandwidths }
  \label{fig:subfig} 
\end{figure}

\section{ Conclusion }

In this paper, the local linear estimators based on Gamma asymmetric
kernels for the unknown drift and conditional variance in
second-order jump-diffusion model. Besides the standard properties
of the local linear estimation constructed with Gaussian symmetric
kernels such as simple bias representation and boundary bias
correction, the local linear smoothing using Gamma asymmetric
kernels possesses some extra advantages such as variable bandwidth,
variance reduction and resistance to sparse design, which is
validated through finite sample simulation study. Theoretically,
under appropriate regularity conditions, we prove that the
estimators constructed with Gamma asymmetric kernels possess the
consistency and asymptotic normality for large sample and verify the
advantages such as bias reduction, robustness and shorter length of
confidence band through simulation experiments for finite sample.

Empirically, the estimators are illustrated through stock index in
China under five-minute high-frequency data and possess some
advantages mentioned above. This means the second-order
jump-diffusion model may be an alternative model to describe the
dynamic of some financial data, especially to explain integrated
economic phenomena that the current observation in empirical finance
usually behaves as the cumulation of all past perturbations.

\bibliographystyle{amsplain}

\section{ Proofs }
\subsection{ Procedure for Assumption \ref{a3}}

Notice that the expectation with respect to the distribution of
$\xi_{n,i}$ depends on the stationary densities of $X_{n,i}$ and
$\widetilde{X}_{n,i}$ because $\xi_{n,i}$ is a convex linear
combination of $X_{n,i}$ and $\widetilde{X}_{n,i}.$

For the case (i): $ E\big[~|h K^{'}(X_{(i-1)\Delta_{n}})|\big] <
\infty.$ For $K_{G(x/h + 1, h)} (u),$ its first-order derivative has
the form of $K^{'}_{G(x/h + 1, h)} (u) =
\Big(\frac{x}{h}\Big)\frac{u^{x/h - 1} \exp(-u/h)}{h^{x/h + 1}
\Gamma(x/h + 1)} - \Big(\frac{1}{h}\Big)\frac{u^{x/h} \exp
(-u/h)}{h^{x/h + 1} \Gamma(x/h + 1)} := \frac{1}{h} K_{1} (u) +
\frac{1}{h} K_{2} (u).$ Then using the well-known properties of the
$\Gamma$ function, the mean of Gamma distribution and the derivative
of the function $gp(x) := g(x) \cdot p(x)$ for stationary process
$X_{t}$, we have
\begin{eqnarray*}
& ~ & E\big[~|h K^{'}(X_{(i-1)\Delta_{n}}) g(X_{(i-1)\Delta_{n}})|
\big]\\ & = & E\big[~|K_{1} (X_{(i-1)\Delta_{n}})
g(X_{(i-1)\Delta_{n}})| \big] + E\big[~|K_{2} (X_{(i-1)\Delta_{n}})
g(X_{(i-1)\Delta_{n}})|
\big]\\
& = & x \frac{h^{x/h} \Gamma(x/h)}{h^{x/h + 1} \Gamma(x/h + 1)}
\int_{0}^{\infty}{\frac{y^{x/h - 1} \exp(- y/h)}{h^{x/h}
\Gamma(x/h)} gp(y)} dy\\ & ~ & + \int_{0}^{\infty}{\frac{y^{x/h}
\exp
(-y/h)}{h^{x/h + 1} \Gamma(x/h + 1)} gp(y)}dy\\
& = & \int_{0}^{\infty}{\frac{y^{x/h - 1} \exp(- y/h)}{h^{x/h}
\Gamma(x/h)} gp(y)} dy + \int_{0}^{\infty}{\frac{y^{x/h} \exp
(-y/h)}{h^{x/h + 1} \Gamma(x/h + 1)} gp(y)}dy\\
& = & E[gp(\xi_{1})] + E[gp(\xi_{2})]\\
& = & E[gp(E(\xi_{1}) + \xi_{1} - E(\xi_{1}))] + E[gp(E(\xi_{2}) +
\xi_{2} - E(\xi_{2}))]\\
& = & 2gp(x) + O(h) < \infty,
\end{eqnarray*}
where $\xi_{1} \stackrel{\mathcal {D}} = G (x/h, h)$, $\xi_{2}
\stackrel{\mathcal {D}} = G (x/h + 1, h)$ and $G$ denotes the Gamma
distribution.

For the case (ii):

 $E\big[~|h^{2} K^{'2}(X_{(i-1)\Delta_{n}})
g(X_{(i-1)\Delta_{n}})|\big] \leq 2
E\big[~|K^{2}_{1}(X_{(i-1)\Delta_{n}}) g(X_{(i-1)\Delta_{n}})|\big]
+ 2 E\big[~|K^{2}_{2}(X_{(i-1)\Delta_{n}})
g(X_{(i-1)\Delta_{n}})|\big].$

\noindent Now we only deal with the first part (the second part can
be dealt with in the similar way). Note that $K_{1}(u) =
\frac{u^{x/h - 1} \exp(- u/h)}{h^{x/h} \Gamma(x/h)}$ can be
considered as a density function for a random variable $\xi_{1}
\stackrel{\mathcal {D}} = G (x/h, h).$ By the property of the
$\Gamma$ function, we have with $\eta_{x} \stackrel{\mathcal {D}} =
G (2 x/h - 1, h),$
\begin{eqnarray*}
& ~ & E\big[~|K^{2}_{1}(X_{(i-1)\Delta_{n}})
g(X_{(i-1)\Delta_{n}})|\big]\\
& = & B_{b}(x) E[g p(\eta_{x})]\\
& = & B_{b}(x) E[g p(E(\xi_{1}) + \xi_{1} - E(\xi_{1}))] = B_{b}(x)
g p(x)\\
& \approx & g p(x) \Big\{
\begin{array}{ll}
\frac{1}{2 \sqrt{\pi}}h^{-1/2}x^{-1/2}~~~if~x/b \rightarrow \infty~(``interior~x'');\\
h^{-1}\frac{\Gamma(2\kappa - 1)}{2^{2\kappa -
1}\Gamma^{2}(\kappa)}~~~if~x/b \rightarrow \kappa~(``boundary~x''),
\end{array}
\end{eqnarray*}
where $B_{b}(x) = \frac{h^{- 1}\Gamma(2x/h - 1)}{2^{2x/b - 1}
\Gamma^{2}(x/h + 1)}$ and the last equation follows from Chen
(\cite{ch2}, P474). Hence, the results of $\lim_{h \rightarrow 0}
h^{1/2} E\big[~|h^{2} K^{'2}(\xi_{n,i}) g(\xi_{n,i})|\big] < \infty$
for ``interior x'' and $\lim_{h \rightarrow 0} h E\big[~|h^{2}
K^{'2}(\xi_{n,i}) g(\xi_{n,i})|\big] < \infty$ for ``boundary x''
hold.

\subsection{ Some Technical Lemmas with Proofs}

We lay out some notations. For $x = (x_{1}, \cdot \cdot \cdot ,
x_{d})$, $\partial_{x_{j}} := \frac{\partial}{\partial{x_{j}}}$,
$\partial^{2}_{x_{j}} := \frac{\partial^{2}}{\partial{x_{j}^{2}}}$,
$\partial^{2}_{x_{i}x_{j}} :=
\frac{\partial^{2}}{\partial{x_{i}}\partial{x_{j}}}$, $\partial{x}
:= (\partial_{x_{1}} , \cdot \cdot \cdot , \partial{x_{d}})^{\ast}$,
and $\partial^{2}_{x} = \big(\partial^{2}_{x_{i}x_{j}}\big)_{1\leq
i,j \leq d}$, where $\ast$ stands for the transpose.

\begin{lemma}
\label{l1} (Shimizu and Yoshida \cite{sy}) Let $Z$ be a
$d$-dimensional solution-process to the stochastic differential
equation $$Z_{t} = Z_{0} + \int_{0}^{t}\mu(Z_{s-})ds +
\int_{0}^{t}\sigma(Z_{s-})dW_{s} +
\int_{0}^{t}\int_{\mathscr{E}}c(Z_{s-} , z)r(\omega, dt, dz),$$
where $Z_{0}$ is a random variable, $\mathscr{E} = \mathbb{R}^{d}
\setminus\{0\}$, $\mu(x), c(x , z)$ are $d$-dimensional vectors
defined on $\mathbb{R}^{d}, \mathbb{R}^{d}\times\mathscr{E}$
respectively, $\sigma(x)$ is a $d \times d$ diagnonal matrix defined
on $\mathbb{R}^{d}$, and $W_{t}$ is a $d$-dimensional vector of
independent Brownian motions.

Let $g$ be a $C^{2(l+1)}$-class function whose derivatives up to
2$(l+1)$th are of polynomial growth. Assume that the coefficients
$\mu(x), \sigma(x),$ and $c(x ,z)$ are $C^{2l}$-class function whose
derivatives with respective to $x$ up to 2$l$th are of polynomial
growth. Under Assumption 5, the following expansion holds
\begin{equation}
\label{6.1} E[g(Z_{t})|\mathscr{F}_{s}] =
\sum_{j=0}^{l}L^{j}g(Z_{s})\frac{\Delta_{n}^{j}}{j!} + R,
\end{equation}
for $t > s$ and $\Delta_{n} = t - s$, where  $R =
\int_{0}^{\Delta_{n}}\int_{0}^{u_{1}} \dots \int_{0}^{u_{l}}
E[L^{l+1}g(Z_{s+u_{l+1}})|\mathscr{F}_{s}]du_{1} \dots du_{l+1}$ is
a stochastic function of order $\Delta_{n}^{l+1}, Lg(x) =
\partial_{x}^{\ast}g(x)\mu(x) + \frac{1}{2}
tr[\partial_{x}^{2}g(x)\sigma(x)\sigma^{\ast}(x)] +
\int_{\mathscr{E}}\{g(x + c(x , z)) - g(x) -
\partial_{x}^{\ast}g(x)c(x , z)\}f(z)dz.$
\end{lemma}

\begin{remark}
\label{r6.1} Consider a particularly important model:

$$\left\{
\begin{array}{ll}dY_{t} = X_{t-}dt,\\
dX_{t} = \mu(X_{t-})dt + \sigma(X_{t-})dW_{t} +
\int_{\mathscr{E}}c(X_{t-} , z)r(w , dt , dz).\end{array} \right.
$$

\noindent As $d$ = 2, we have
\begin{equation}
\label{6.2}
\begin{array}{ll}
Lg(x , y) = x(\partial g/\partial y) + \mu(x)(\partial g/\partial x)
+ \frac{1}{2}\sigma^{2}(x)(\partial^{2} g/\partial
x^{2})\\
~~~~~~~~~~~~~+ \int_{\mathscr{E}}\{g(x + c(x , z) , y) - g(x , y) -
\frac{\partial g}{\partial x} \cdot c(x , z) \}f(z)dz.
\end{array}
\end{equation}
Based on the second-order infinitesimal operator (\ref{6.2}), we can
calculate many mathematical expectations involving
$\widetilde{X}_{i\Delta_{n}}$, for instance (\ref{2.3}) and
(\ref{2.4}) which provide the basis for estimators (\ref{2.5}) and
(\ref{2.6}).
\end{remark}

\begin{lemma}
\label{l2} Under Assumption \ref{a1}, \ref{a2} and \ref{a5}, let
$$\mu_{n}^{\ast} (x) = \frac{\sum_{i=1}^{n}
w_{i-1}^{\ast}(\frac{X_i-X_{i-1}}{\Delta_{n}})}{\sum_{i=1}^{n}
w_{i-1}^{\ast}}$$ and
$$ M_{n}^{\ast} (x) = \frac{\sum_{i=1}^{n}
w_{i-1}^{\ast}\frac{(X_i-X_{i-1})^2}{\Delta_{n}}}{\sum_{i=1}^{n}
w_{i-1}^{\ast}}.$$ where
\begin{eqnarray*}
w_i^\ast & = & K_{G(x/h + 1, h)}(X_{i})(\sum_{j=1}^n K_{G(x/h + 1,
h)}(X_{j-1})(X_{j-1} - x)^2\\ & ~ & - (X_{i} - x) \sum_{j=1}^{n}
K_{G(x/h + 1, h)}(X_{j - 1})(X_{j-1} - x))
\end{eqnarray*} then
$$\mu_{n}^{\ast}(x) \stackrel{p}{\rightarrow} \mu(x),~~~M_{n}^{\ast}(x) \stackrel{p}{\rightarrow} M(x).$$ Furthermore, for ``interior x'', if
$h = O((n\Delta_{n})^{-2/5}),$ then
\begin{align*}
& \sqrt{n\Delta_{n} h^{1/2}}\big(\mu_{n}^{\ast}(x) - \mu(x) - h
B_{\mu_{n}^{\ast}(x)}\big) \stackrel{d}{\rightarrow} N\Big( 0 ,
\frac{M(x)}{2 \sqrt{\pi} x^{1/2} p(x)}\Big),\\
& \sqrt{n\Delta_{n} h^{1/2}}\big(M_{n}^{\ast}(x) - M(x) - h
B_{M_{n}^{\ast}(x)}\big) \stackrel{d}{\rightarrow} N\Big( 0 ,
\frac{\int_{\mathscr{E}}c^{4}(x , z)f(z)dz}{2 \sqrt{\pi} x^{1/2}
p(x)}\Big),
\end{align*}
for ``boundary x'', if $h = O((n\Delta_{n})^{-1/5}),$ then
\begin{align*}
& \sqrt{n\Delta_{n} h}\big(\mu_{n}^{\ast}(x) - \mu(x) - h^{2}
B^{'}_{\mu_{n}^{\ast}(x)}\big) \stackrel{d}{\rightarrow} N\Big( 0 ,
\frac{M(x) \Gamma(2\kappa + 1)}{2^{2\kappa + 1} \Gamma^{2}(\kappa + 1) p(x)}\Big),\\
& \sqrt{n\Delta_{n} h}\big(M_{n}^{\ast}(x) - M(x) - h^{2}
B^{'}_{M_{n}^{\ast}(x)}\big) \stackrel{d}{\rightarrow} N\Big( 0 ,
\frac{\int_{\mathscr{E}}c^{4}(x , z)f(z)dz \Gamma(2\kappa +
1)}{2^{2\kappa + 1} \Gamma^{2}(\kappa + 1) p(x)}\Big),
\end{align*}
where
$B_{\mu_{n}^{\ast}(x)},~B^{'}_{\mu_{n}^{\ast}(x)},~B_{M_{n}^{\ast}(x)},~B^{'}_{M_{n}^{\ast}(x)}$
denotes the bias of the estimators of
$\mu_{n}^{\ast}(x),~M_{n}^{\ast}(x),$ respectively, that is
\begin{align*}
& B_{\mu_{n}^{\ast}(x)} = \frac{x}{2}\mu^{''}(x),~~~~~B^{'}_{\mu_{n}^{\ast}(x)} = \frac{1}{2}(2 + \kappa)\mu^{''}(x)\\
& B_{M_{n}^{\ast}(x)} =
\frac{x}{2}M^{''}(x),~~~~~B^{'}_{M_{n}^{\ast}(x)} = \frac{1}{2}(2 +
\kappa)M^{''}(x).
\end{align*}
\end{lemma}

\begin{remark}
This lemma considered the asymptotic properties of the local linear
estimation for stationary jump-diffusion model (\ref{1.1}) using
Gamma asymmetric kernels, which is different from that in Hanif
\cite{hm1}.

After carefully sketching the paper of Hanif \cite{hm1}, we found
that the part $S_{n, k}$ of the weight $\omega_{i}^{LL}(x, b)$ (that
is $w_i^\ast$ here) in (2.8) or (2.9) in Hanif \cite{hm1} should be
$S_{n, k} = \sum_{i = 1}^{n} K_{G(x/b + 1, b)}(X_{i\Delta_{n, T}})
\cdot (X_{i\Delta_{n, T}} - x)^{k}$ , not $S_{n, k} = \sum_{i =
1}^{n} K_{G(x/b + 1, b)}(X_{i\Delta_{n, T}}) \cdot (X_{i\Delta_{n,
T}})^{k}.$ So in the detailed proof of Lemma 4, Theorem 1 and
Theorem 2 in Hanif \cite{hm1}, we should consider $K_{G(x/b + 1,
b)}(X_{i\Delta_{n, T}}) \cdot (X_{i\Delta_{n, T}} - x)^{k}$ and
$K_{G(x/b + 1, b)}(X_{s-}) \cdot (X_{s-} - x)^{k}$, not $K_{G(x/b +
1, b)}(X_{i\Delta_{n, T}}) \cdot (X_{i\Delta_{n, T}})^{k}$ or
$K_{G(x/b + 1, b)}(X_{s-}) \cdot (X_{s-})^{k}.$ According to the
similar approach as Chen \cite{ch3}, we will give a modified proof
to the stationary results of Lemma 4, Theorem 1 and Theorem 2 in
Hanif \cite{hm1}. Hence, the central limit theorems of
$\mu_{n}^{\ast}(x)$ and $M_n^{\ast}(x)$ are different from those in
Hanif \cite{hm1} for the stationary case.
\end{remark}
\begin{proof}

For convenience, we still use the same notations as that in Hanif
\cite{hm1}. The part $S_{n, k}$ of the weight $\omega_{i}^{LL}(x,
b)$ in (2.8) or (2.9) in Hanif \cite{hm1} should be $S_{n, k} =
\sum_{i = 1}^{n} K_{G(x/b + 1, b)}(X_{i\Delta_{n, T}}) \cdot
(X_{i\Delta_{n, T}} - x)^{k}$ , not $S_{n, k} = \sum_{i = 1}^{n}
K_{G(x/b + 1, b)}(X_{i\Delta_{n, T}}) \cdot (X_{i\Delta_{n,
T}})^{k}.$ So in the detailed proof of Lemma 4, Theorem 1 and
Theorem 2 in Hanif \cite{hm1}, we should consider $K_{G(x/b + 1,
b)}(X_{i\Delta_{n, T}}) \cdot (X_{i\Delta_{n, T}} - x)^{k}$ and
$K_{G(x/b + 1, b)}(X_{s-}) \cdot (X_{s-} - x)^{k}$, not $K_{G(x/b +
1, b)}(X_{i\Delta_{n, T}}) \cdot (X_{i\Delta_{n, T}})^{k}$ or
$K_{G(x/b + 1, b)}(X_{s-}) \cdot (X_{s-})^{k}.$

The key point of the detailed proof for stationary case of Lemma 4
in Hanif \cite{hm1} is
\begin{eqnarray*}
& ~ & \frac{1}{T} \int_{0}^{T} K_{G(x/b + 1, b)}(X_{s-})(X_{s-} -
x)^{k}\frac{d[X]^{c}_{s}}{\sigma^{2}(X_{s-})}\\
& = & \frac{1}{T} \int_{0}^{\infty} K_{G(x/b + 1, b)}(a)(a -
x)^{k}\frac{L_{X}(T, a)}{\sigma^{2}(a)}da\\
& = & \int_{0}^{\infty} K_{G(x/b + 1, b)}(a)(a - x)^{k}
\frac{\bar{L}_{X}(T, a)}{T} da\\
& = & \int_{0}^{\infty} K_{G(x/b + 1, b)}(a)(a - x)^{k} p(a) da\\
& = & E[(\xi - x)^{k} p(\xi)] := \alpha_{k}(x),
\end{eqnarray*}
where $k = 0, 1, 2$, $\xi \stackrel{\mathcal {D}} = G (x/h + 1, h)$
and $G$ denotes the Gamma distribution.

According to the result (A.1) and (A.2) in Chen (\cite{ch3}, P321),
it can be shown that
$$\alpha_{k}(x) = \sum_{j = 0}^{2 - k}p^{(j)}(x)E(\xi - x)^{j + k}/j! + o_{p}\{E(\xi - x)^{2}\}.$$
As $\xi$ is the $G (x/h + 1, h)$ random variable, $E(\xi - x) =
h_{n}$, $E(\xi - x)^{2} = x h_{n} + 2 h_{n}^{2}$ and $E(\xi - x)^{l}
= O(h_{n}^{2})$ for $3 \leq l.$ Thus, we can deduce
\begin{align}
& \label{6.3} \alpha_{0}(x) = p(x) + p^{(1)}(x)h_{n} +
\frac{p^{(2)}(x)}{2}[x
h_{n} + 2h_{n}^{2}] + o_{p}(h_{n}^{2}),\\
& \label{6.4} \alpha_{1}(x) = p(x)h_{n} + p^{(1)}(x)[x h_{n} +
2h_{n}^{2}] +
o_{p}(h_{n}^{2}),\\
& \label{6.5} \alpha_{2}(x) = p(x)[x h_{n} + 2h_{n}^{2}] +
o_{p}(h_{n}^{2}).
\end{align}

\noindent $\bf{Firstly},$ we calculate the bias $A_{22}$ for
$\hat{M}_{LL}^{1}(x, b) - M^{1}(x)$ in Hanif (\cite{hm1}, P966). We
write $A_{22}$ as (4.6) in Hanif (\cite{hm1}, P966)
\begin{eqnarray*}
A_{22} & = & \frac{\sum_{i=1}^{n}
w_{i-1}^{\ast}(\frac{X_i-X_{i-1}}{\Delta_{n}} -
\mu(x))}{\sum_{i=1}^{n} w_{i-1}^{\ast}}\\
& = & \frac{\sum_{i = 1}^{n}[S_{n, 2} - S_{n, 1} \cdot
(X_{i\Delta_{n, T}} - x)]K_{G(x/b + 1, b)}(X_{i\Delta_{n,
T}})\Delta_{n, T}(\mu(X_{i\Delta_{n, T}}) - \mu(x)) +
o_{a.s.}(1)}{\Delta_{n, T}(S_{n, 0} \cdot S_{n, 2} - S^{2}_{n,
1})}\\
& = & \frac{\frac{1}{n^{2}}\sum_{i = 1}^{n}[S_{n, 2} - S_{n, 1}
\cdot (X_{i\Delta_{n, T}} - x)]K_{G(x/b + 1, b)}(X_{i\Delta_{n,
T}})(\mu(X_{i\Delta_{n, T}}) - \mu(x)) +
o_{a.s.}(1)}{\frac{1}{n^{2}}(S_{n, 0} \cdot S_{n, 2} - S^{2}_{n,
1})},
\end{eqnarray*}
where $S_{n, k} = \sum_{i = 1}^{n} K_{G(x/b + 1, b)}(X_{i\Delta_{n,
T}}) \cdot (X_{i\Delta_{n, T}} - x)^{k}.$

\noindent Substituting results (\ref{6.3}) - (\ref{6.5}) to the
denominator of $A_{22}$, we may derive
$$A_{22}^{Den} = \alpha_{0}(x) \cdot \alpha_{2}(x) - \alpha^{2}_{1}(x) = p^{2}(x) [x h_{n} + 2h_{n}^{2}] + o(h^{2}_{n}).$$
Taylor expanding $\mu(X_{i\Delta_{n, T}})$ at $x$ for the numerator
of $A_{22}$,
\begin{eqnarray*}
A_{22}^{Num} & = & \frac{1}{n^{2}}\sum_{i = 1}^{n}[S_{n, 2} - S_{n,
1} \cdot (X_{i\Delta_{n, T}} - x)]K_{G(x/b + 1, b)}(X_{i\Delta_{n,
T}})(\mu(X_{i\Delta_{n, T}}) - \mu(x))\\
& = & \frac{1}{n^{2}}\sum_{i = 1}^{n}[S_{n, 2} - S_{n, 1} \cdot
(X_{i\Delta_{n, T}} - x)]K_{G(x/b + 1, b)}(X_{i\Delta_{n,
T}})\Big(\mu^{'}(x)(X_{i\Delta_{n, T}} - x)\\
& ~ & + \frac{1}{2}\mu^{''}(X_{i\Delta_{n, T}})(X_{i\Delta_{n, T}} -
x)^{2} + \frac{1}{6}\mu^{'''}(\zeta_{n, i})(X_{i\Delta_{n, T}} -
x)^{3}\Big)\\
& = & \frac{1}{n^{2}}\sum_{i = 1}^{n}[S_{n, 2} - S_{n, 1} \cdot
(X_{i\Delta_{n, T}} - x)]K_{G(x/b + 1, b)}(X_{i\Delta_{n, T}})\Big(
\frac{1}{2}\mu^{''}(X_{i\Delta_{n, T}})(X_{i\Delta_{n, T}}
- x)^{2}\\
& ~ & + \frac{1}{6}\mu^{'''}(\zeta_{n, i})(X_{i\Delta_{n, T}} -
x)^{3}\Big)
\end{eqnarray*}
by virtue of the fact that
\begin{eqnarray*}
& ~ & \sum_{i = 1}^{n} \omega_{i}^\ast \times ({X}_{i} - x)\\
& = & \sum_{i = 1}^{n} K_{G(x/h + 1, h)}\big({X}_{i}\big)({X}_{i} -
x) \times \sum_{j=1}^n K_{G(x/h + 1, h)}\big({X}_{j-1}\big)
({X}_{j-1} - x)^2\\
& ~ & - \sum_{i=1}^n K_{G(x/h + 1, h)}\big({X}_{i}\big) ({X}_{i} -
x)^2 \times \sum_{j=1}^{n} K_{G(x/h + 1,
h)}\big({X}_{j-1}\big) ({X}_{j-1} - x)\\
& = & 0,
\end{eqnarray*}
where $\zeta_{n, i} = \theta x + (1 - \theta) X_{i\Delta_{n, T}}.$

With $K_{G(x/b + 1, b)}(\cdot) g(\cdot) (\cdot - x)^{k}$ instead of
$K_{G(x/b + 1, b)}(\cdot) (\cdot - x)^{k}$ in Lemma 4 of Hanif
\cite{hm1}, we can similarly deduce
$$\frac{1}{n}\sum_{i = 1}^{n} K_{G(x/b + 1, b)}(X_{i\Delta_{n, T}}) g(X_{i\Delta_{n, T}}) (X_{i\Delta_{n, T}} - x)^{k}
\stackrel{P} \rightarrow E[g(\xi)(\xi - x)^{k} p(\xi)] :=
\beta_{k}(x),$$ where $k = 0, 1, 2$, $\xi \stackrel{\mathcal {D}} =
G (x/h + 1, h)$, $G$ denotes the Gamma distribution and $g(\cdot) =
\frac{1}{2} \mu^{''}(\cdot).$

According to the result (A.1) and (A.3) in Chen (\cite{ch3}, P321),
it can be shown that with $r(x) = g(x) \cdot p(x)$
$$\beta_{k}(x) = \sum_{j = 0}^{2 - k}r^{(j)}(x)E(\xi - x)^{j + k}/j! + o_{p}\{E(\xi - x)^{2}\}.$$
As $\xi$ is the $G (x/h + 1, h)$ random variable, $E(\xi - x) =
h_{n}$, $E(\xi - x)^{2} = x h_{n} + 2 h_{n}^{2}$ and $E(\xi - x)^{l}
= O(h_{n}^{2})$ for $3 \leq l.$ Thus, we can deduce
\begin{align}
& \label{6.6} \beta_{0}(x) = r(x) + r^{(1)}(x)h_{n} +
\frac{r^{(2)}(x)}{2}[x
h_{n} + 2h_{n}^{2}] + o_{p}(h_{n}^{2}),\\
& \label{6.7} \beta_{1}(x) = r(x)h_{n} + r^{(1)}(x)[x h_{n} +
2h_{n}^{2}] +
o_{p}(h_{n}^{2}),\\
& \label{6.8} \beta_{2}(x) = r(x)[x h_{n} + 2h_{n}^{2}] + o_{p}(h_{n}^{2}),\\
& \label{6.9} \beta_{3}(x) = O(h_{n}^{2}).
\end{align}
\noindent Substituting results (\ref{6.6}) - (\ref{6.9}) to the
numerator of $A_{22}$, we may derive
$$A_{22}^{Num} = \alpha_{2}(x) \cdot \beta_{2}(x) - \alpha_{1}(x) \cdot \beta_{3}(x) =  \frac{1}{2} \mu^{''}(x) \cdot p^{2}(x) [x
h_{n} + 2h_{n}^{2}]^{2} + o(h^{2}_{n}).$$ So the bias for
$\hat{M}_{LL}^{1}(x, b) - M^{1}(x)$ in Hanif (\cite{hm1}, P966) is
\begin{eqnarray*}
A_{22} & = & \frac{A_{22}^{Num}}{A_{22}^{Den}}\\
& = & \frac{\frac{1}{2} \mu^{''}(x) \cdot p^{2}(x) [x h_{n} +
2h_{n}^{2}]^{2} + o(h^{2}_{n})}{p^{2}(x) [x h_{n} + 2h_{n}^{2}] +
o(h^{2}_{n})}\\
& = & \frac{1}{2} \mu^{''}(x)[x h_{n} + 2h_{n}^{2}] + o(h^{2}_{n})\\
& = & \Big\{ \begin{array}{ll} \frac{x}{2} \mu^{''}(x) h_{n} + o(h_{n})~~~if~x/b \rightarrow \infty~(``interior~x'');\\
h_{n}^{2}[\frac{1}{2} \mu^{''}(x) (2 + \kappa)]~~~if~x/b \rightarrow
\kappa~(``boundary~x'')
\end{array}\\
\end{eqnarray*}
\noindent $\bf{Secondly},$ we calculate two parts $[B_{22}, B_{22}]$
and $[C_{22}, C_{22}]$ related to the variance of the asymptotic
normality for $\hat{M}_{LL}^{1}(x, b) - M^{1}(x)$ in Hanif
(\cite{hm1}, P966).
\begin{eqnarray*}
& ~ & [B_{22}, B_{22}]\\ & = & \frac{\sum_{i = 1}^{n}[\Delta_{n,
T}S_{n, 2} - \Delta_{n, T}S_{n, 1}\cdot (X_{i\Delta_{n, T}} -
x)]^{2}K^{2}_{G(x/b + 1, b)}(X_{i\Delta_{n, T}})\int_{i\Delta_{n,
T}}^{(i+1)\Delta_{n, T}}\sigma^{2}(X_{s-})ds}{\left(\Delta^{2}_{n,
T}(S_{n, 0} \cdot
S_{n, 2} - S^{2}_{n, 1})\right)^{2}}\\
& = & \frac{\sum_{i = 1}^{n}[\Delta_{n, T}S_{n, 2} - \Delta_{n,
T}S_{n, 1}\cdot (X_{i\Delta_{n, T}} - x)]^{2}K^{2}_{G(x/b + 1,
b)}(X_{i\Delta_{n, T}}) \Delta_{n, T} \sigma^{2}(X_{i\Delta_{n, T}})
+ o_{a.s.}(1)}{\left(\Delta^{2}_{n, T}(S_{n, 0} \cdot S_{n, 2} -
S^{2}_{n, 1})\right)^{2}}.
\end{eqnarray*}
Due to $A_{22}^{Den}$, we have
$$[B_{22}, B_{22}]^{Den} = p^{4}(x) [x h_{n} + 2h_{n}^{2}]^{2} + o(h^{4}_{n}).$$
As for $[B_{22}, B_{22}]^{Num}$, we have
\begin{eqnarray*}
& ~ &[B_{22}, B_{22}]^{Num}\\
 & = & \sum_{i = 1}^{n}[\Delta^{2}_{n,
T}S^{2}_{n, 2} - 2 \Delta^{2}_{n, T} S_{n, 2} S_{n, 1}\cdot
(X_{i\Delta_{n, T}} - x)\\ & ~ & + \Delta^{2}_{n, T} S^{2}_{n,
1}\cdot (X_{i\Delta_{n, T}} - x)^{2}] K^{2}_{G(x/b + 1,
b)}(X_{i\Delta_{n,
T}}) \Delta_{n, T} \sigma^{2}(X_{i\Delta_{n, T}})\\
& = & [B_{22}, B_{22}]^{Num}_{1} + [B_{22}, B_{22}]^{Num}_{2} +
[B_{22}, B_{22}]^{Num}_{3}.
\end{eqnarray*}
According to (\ref{6.3}) - (\ref{6.9}) with $g(\cdot) =
\sigma^{2}(\cdot)$ and $A_{h_{n}}(x)$ in Chen (\cite{ch2}, P474), it
can be shown that $[B_{22}, B_{22}]^{Num}_{1}$ is larger than the
others (which has the lowest infinitesimal order). Under the similar
calculus as $\beta_{0}(x)$, the dominant one has the following
expression
$$[B_{22},
B_{22}]^{Num}_{1} = A_{h_{n}}(x) \sigma^{2}(x) p^{3}(x) [x h_{n} +
2h_{n}^{2}]^{2} + o(h^{4}_{n}).$$ Hence, we get $$[B_{22}, B_{22}] =
\frac{A_{h_{n}}(x) \sigma^{2}(x) p^{3}(x) [x h_{n} + 2h_{n}^{2}]^{2}
+ o(h^{4}_{n})}{p^{4}(x) [x h_{n} + 2h_{n}^{2}]^{2} + o(h^{4}_{n})}
= A_{h_{n}}(x) \frac{\sigma^{2}(x)}{p(x)}.$$ Similarly, we can prove
$$[C_{22}, C_{22}] =
A_{h_{n}}(x) \frac{\int_{\mathscr{E}}c^2(x,z)f(z)dz}{p(x)}.$$ So the
variance of the asymptotic normality for $\hat{M}_{LL}^{1}(x, b) -
M^{1}(x)$ in Hanif (\cite{hm1}, P966) should be $A_{h_{n}}(x)
\frac{M(x)}{p(x)}.$

The modified proof of Theorem 2 in Hanif \cite{hm1} is similar to
that of Theorem 1, so we omit it. One can refer to Song \cite{syp2}
for similar procedure.
\end{proof}

\subsection{ The proof of Theorem \ref{thm1}}

\begin{proof}

Here we only prove the first result; the second is analogical. By
Lemma \ref{l2}, it suffice to show that :
$$ \hat{\mu}_n(x)-\mu^\ast _n(x) \stackrel{p}{\rightarrow} 0.$$
$\bf{Firstly}$, we prove that
\begin{equation}
\label{6.10} \frac{1}{n^2}\sum_{i=1}^n
w_{i-1}-\frac{1}{n^2}\sum_{i=1}^n w_{i-1}^\ast
\stackrel{p}{\rightarrow} 0.
\end{equation}

\noindent To this end, we should prove that
\begin{equation}
\label{6.11} \frac{1}{n}\sum_{i=1}^{n}K(\tilde{X}_{i-1})-
\frac{1}{n}\sum_{i=1}^{n}K(X_{i-1}) \stackrel{p}{\rightarrow} 0,
\end{equation}
\begin{equation}
\label{6.12}
\frac{1}{n}\sum_{i=1}^{n}K(\tilde{X}_{i-1})(\tilde{X}_{i - 1} - x) -
\frac{1}{n}\sum_{i=1}^{n}K(X_{i-1})(X_{i-1} - x)
\stackrel{p}{\rightarrow} 0,
\end{equation}
\begin{equation}
\label{6.13} \frac{1}{n}\sum_{i=1}^{n}K(\tilde{X}_{i-1})(\tilde{X}_i
- x)^2- \frac{1}{n}\sum_{i=1}^{n}K(X_{i-1})(X_{i-1} - x)^2
\stackrel{p}{\rightarrow} 0.
\end{equation}

\noindent For (\ref{6.11}), let $\varepsilon_{1 , n} = \frac{1}{n}
\sum_{i=1}^{n} K_{G(x/h + 1,
h)}\big(\widetilde{X}_{(i-1)\Delta_{n}}\big) - \frac{1}{n}
\sum_{i=1}^{n} K_{G(x/h + 1, h)}\big(X_{(i-1)\Delta_{n}}\big).$
\begin{eqnarray}
& & \max_{1 \leq i \leq n}\big|\widetilde{X}_{(i-1)\Delta_{n}} -
{X}_{(i-1)\Delta_{n}}\big| \nonumber \\
  & \leq & \max_{1 \leq i \leq n}\frac{1}{\Delta_{n}}\big|\int
_{(i-2)\Delta_{n}}^{(i-1)\Delta_{n}}(X_{s-} -
{X}_{(i-1)\Delta_{n}}) d s\big| \nonumber \\
& \leq & \max_{1 \leq i \leq n} \sup_{(i-2)\Delta_{n} \leq s \leq
(i-1)\Delta_{n}}\big|X_{s-} - {X}_{(i-1)\Delta_{n}}\big| \nonumber \\
& = & \label{6.14} O_{a.s.}(\sqrt {\Delta_{n}\log(1/\Delta_{n})}),
\end{eqnarray}
the last asymptotic equation for the order of magnitude, one can
refer Bandi and Nguyen (\cite{bn}, equations 94 and 95).

\noindent By the mean-value theorem, stationarity, {Assumptions
\ref{a3}, \ref{a5} and (\ref{6.14})}, we obtain
\begin{eqnarray*}
E[|\varepsilon_{1 , n}|] & \leq & E[\frac{1}{n} \sum
_{i=1}^{n}|K^{'}(\xi_{n,i}) \big(\widetilde{X}_{(i-1)\Delta_{n}} -
{X}_{(i-1)\Delta_{n}}\big)|]\\
& = & E[|K^{'}(\xi_{n,2}) \big(\widetilde{X}_{\Delta_{n}} -
{X}_{\Delta_{n}}\big)|]\\
& \leq & \sqrt{\Delta_{n} \log (1/\Delta_{n})}
E[|K^{'}(\xi_{n,2})|]\\
& = & \frac{\sqrt{\Delta_{n} \log (1/\Delta_{n})}}{h} E[|h
K^{'}(\xi_{n,2})|]
 \rightarrow 0,
\end{eqnarray*}
where $\xi_{n,2}$ = $\theta X_{\Delta_{n}} + (1 -
\theta)\widetilde{X}_{\Delta_{n}} ~0 \leq \theta \leq 1.$ Hence,
(\ref{6.11}) follows from Chebyshev's inequality.

For (\ref{6.12}) we should prove that
\begin{equation}
\label{6.15}
\delta_{1,n}=\frac{1}{n}\sum_{i=1}^{n}K(\tilde{X}_{i-1})(\tilde{X}_i
- x) - \frac{1}{n}\sum_{i=1}^{n}K(X_{i-1})(\tilde{X}_{i} - x)
\stackrel{p}{\rightarrow} 0,
\end{equation}
and
\begin{equation}
\label{6.16}
\delta_{2,n}=\frac{1}{n}\sum_{i=1}^{n}K(X_{i-1})(\tilde{X}_i - x) -
\frac{1}{n}\sum_{i=1}^{n}K(X_{i-1})(X_{i-1} - x)
\stackrel{p}{\rightarrow} 0.
\end{equation}

\noindent Under Assumption \ref{a3} and \ref{a5}, we have
\begin{eqnarray*}
\big|E[\delta_{1,n}]\big| & = & \big|E[ \frac{1}{n} \sum_{i=1}^{n}
\big\{K(\tilde{X}_{i-1})-K(X_{i-1})\big\}(\tilde{X}_{i} - x)]\big|\\
& = & \big|E[ \big\{K(\tilde{X}_{i-1})-K(X_{i-1})\big\}
E[\tilde{X}_{i} - x\big|\mathscr{F}_{i-1}]]\big|\\
& = & |E[ K'(\xi_{n,i})
(X_{i - 1}-\tilde{X}_{i-1}) (X_{i-1} - x + O_P(\Delta_n))]|\\
& \leq & \sqrt{\Delta_{n} \log (1/\Delta_{n})} E[ |K'(\xi_{n,i})
(X_{i-1} - x + O_P(\Delta_n))|]\\
& = & \frac{\sqrt{\Delta_{n} \log (1/\Delta_{n})}}{h_{n}} E[
|h_{n} K'(\xi_{n,i})(X_{i-1} - x + O_P(\Delta_n))|]\\
& \rightarrow & 0,
\end{eqnarray*}
by stationarity, the mean-value theorem and Remark \ref{r6.1}. So
$E[\delta_{1,n}]\rightarrow 0.$

If we can prove $Var[\delta_{1,n}] \rightarrow 0$, then (\ref{6.15})
holds. Now we calculate $Var[\delta_{1,n}]$
\begin{eqnarray*}
Var[\delta_{1,n}] & = &
\frac{1}{nh_{n}}Var[\frac{1}{\sqrt{n}}\sum_{i=1}^{n} \sqrt{h_{n}}
K^{'}(\xi_{n,i})(\widetilde{X}_{i-1}
- X_{i-1})(\widetilde{X}_{i} - x)]\\
& =: & \frac{1}{nh_{n}}Var[\frac{1}{\sqrt{n}}\sum_{i=1}^{n}f_{i}].\\
\end{eqnarray*}
where $f_{i} := \sqrt{h_{n}} K^{'}(\xi_{n,i})(\widetilde{X}_{i-1} -
X_{i-1})(\widetilde{X}_{i} - x).$

\noindent By Remark \ref{r6.1} and Assumption \ref{a3} and \ref{a5},
we get
\begin{eqnarray*}
E[f_{i}^{2}] & = & E\big[h_{n} K^{'2}(\xi_{n,i})(\widetilde{X}_{i-1}
- X_{i-1})^{2}E[(\widetilde{X}_{i} - x)^{2} | \mathscr{F}_{i - 1}]\big]\\
& \leq & \frac{\Delta_{n} \log (1/\Delta_{n})}{h_{n}}
E\big[h^{2}_{n} K^{'2}(\xi_{n,i}) ((X_{i} - x)^{2} +
O_{p}(\Delta_{n}))\big]
\\
& = &  \frac{\Delta_{n}\log(1/\Delta_{n})}{h_{n}} \Big \{
\begin{array}{ll}
\frac{1}{h_{n}^{1/2}} \cdot h_{n}^{1/2} E\Big[h_{n}^{2}
K^{'2}(\xi_{n,i}) ((X_{i} - x)^{2} +
O_{p}(\Delta_{n})) \Big]\\~~~\big(if~x/b \rightarrow \infty~:~``interior~x'' \big);\\
\frac{1}{h_{n}} \cdot h_{n} E\Big[h_{n}^{2}K^{'2}(\xi_{n,i}) ((X_{i}
- x)^{2} + O_{p}(\Delta_{n}))\Big]\\~~~\big(if~x/b \rightarrow
\kappa~:~``boundary~x'' \big)
\end{array}\\
& \approx & C \Big\{ \begin{array}{ll} \frac{\Delta_{n}\log(1/\Delta_{n})}{h_{n}^{3/2}}~~~if~x/b \rightarrow \infty~(``interior~x'');\\
\frac{\Delta_{n}\log(1/\Delta_{n})}{h_{n}^{2}}~~~if~x/b \rightarrow
\kappa~(``boundary~x'')
\end{array}\\
& < & \infty
\end{eqnarray*}

We notice that $f_i$ is stationary under Assumption \ref{a2} and
$\rho$-mixing with the same size as process
$\{\tilde{X}_{i\Delta_n}; i=1,2,...\}$ and $\{X_{i\Delta_n};
i=1,2,...\}$. So from Lemma 10.1.c with p=q=2 in Lin and
Bai(\cite{lb}, p. 132), we have
\begin{eqnarray*}
\big|Var[\frac{1}{\sqrt{n}}\sum_{i=1}^n f_i]\big|& = &
\big|\frac{1}{n}[\sum_{i=1}^n Var(f_i)+2
\sum_{j=1}^{n-1}\sum_{i=j+1}^n (Ef_if_j-Ef_iEf_j)]\big|\\
& = & Var(f_i)+\frac{2}{n}\big| \sum_{j=1}^{n-1}\sum_{i=j+1}^n
(Ef_if_j-Ef_iEf_j)]\big|\\
& \leq & Var(f_i)+ \frac{2}{n}\sum_{j=1}^{n-1}\sum_{i=j+1}^n
\big|Ef_if_j-Ef_iEf_j\big]\\
& \leq &
Var(f_i)+\frac{8}{n}\sum_{j=1}^{n-1}\sum_{i=j+1}^{n}\rho((i-j)\Delta_{n})(Ef^{2}_{i})^{\frac{1}{2}}(Ef^{2}_{j})^{\frac{1}{2}}\\
& = &
Var(f_i)+\frac{8}{n}\sum_{j=1}^{n-1}\sum_{i=j+1}^{n}\rho((i-j)\Delta_{n})Ef^{2}_{i}
\end{eqnarray*}
We have proved $Ef_i^2 < \infty$ above, so the first part in the
last equality $Var(f_i) < \infty $. Moreover, under Assumption
\ref{a2}, we have $\sum_{j=i+1}^n \rho((j-i)\Delta_n)
=O(\frac{1}{\Delta^{\alpha}_n})$. So $Var(\delta_{1,n}) =
\frac{1}{nh_n\Delta^{\alpha}_n}\rightarrow 0$ as $nh_n\Delta^{1 +
\alpha}_n\rightarrow \infty.$

Similar to the proof of (\ref{6.15}), we prove (\ref{6.16}) by
verifying $E[\delta_{2,n}]\rightarrow 0$ and
$Var[\delta_{2,n}]\rightarrow 0$. From the stationarity, Remark
\ref{r6.1} and Assumptions \ref{a3} and \ref{a5}, we have
\begin{eqnarray*}
E[\delta_{2,n}] & = & E\big[K(X_{i-1})E\big[\widetilde{X}_{i} - X_{i-1}|\mathscr{F}_{i-1}\big]\big]\\
& = & \Delta_{n} E\big[K(X_{i-1}) \mu({X}_{i-1})\big] +
O(\Delta_{n}^{2}) \rightarrow 0.
\end{eqnarray*}
and
\begin{eqnarray*}
Var[\delta_{2,n}] & = & \frac{\Delta_{n}}{n h_{n}}
Var\big[\frac{1}{\sqrt{n}}\sum_{i=1}^{n}
\sqrt{h_{n}} K(X_{i-1}) \frac{1}{\sqrt{\Delta_{n}}}(\widetilde{X}_{i} - X_{i-1})\big]\\
& =: & \frac{\Delta_{n}}{n h_{n}}
Var[\frac{1}{\sqrt{n}}\sum_{i=1}^{n}g_{i}],
\end{eqnarray*}
where
\begin{eqnarray*}
E\big[g_{i}^{2}\big] & = &
E\big[h_{n} K^{2}(X_{i-1}) E\big[\frac{(\widetilde{X}_{i} - X_{i-1})^{2}}{\Delta_{n}} | \mathscr{F}_{i - 1}\big]\big]\\
& = & \frac{1}{3}E\big[h_{n} K^{2}(X_{i-1})
\big(\sigma^{2}({X}_{i-1}) +
\int_{\mathscr{E}}c^{2}({X}_{i-1}, z)f(z)dz\big)\big]\\
& \approx & \Big\{ \begin{array}{ll} O(h^{1/2}_{n})~~~if~x/b \rightarrow \infty~(``interior~x'');\\
O(1)~~~if~x/b \rightarrow \kappa~(``boundary~x'')
\end{array}\\
& < & \infty.
\end{eqnarray*}
by the Remark \ref{r6.1} and Assumption \ref{a1}, \ref{a3}. Hence,
$Var[\delta_{2,n}]=O(\frac{\Delta^{1-\alpha}_n}{n h_n})\rightarrow
0$ under Assumption \ref{a5}.

The proof of (\ref{6.13}) is similar to that of (\ref{6.12}), so we
omit it.

\noindent $\bf{Secondly}$, we prove
\begin{equation}
\label{6.17} \delta_n:=\frac{1}{n^2}
\sum_{i=1}^nw_{i-1}\frac{\widetilde{X}_{i+1}-\widetilde{X}_{i}}{\Delta_n}
-\frac{1}{n^2}\sum_{i=1}^nw^\ast_{i-1}\frac{X_{i}-X_{i-1}}{\Delta_n}
\stackrel{p}{\rightarrow} 0.
\end{equation} which suffice to prove
\begin{equation}
\label{6.18}
\frac{1}{n^2}\sum_{i=1}^nw^\ast_{i-1}\frac{\widetilde{X}_{i+1}-\widetilde{X}_{i}}{\Delta_n}
-\frac{1}{n^2}\sum_{i=1}^nw^\ast_{i-1}\frac{X_{i}-X_{i-1}}{\Delta_n}
\stackrel{p}{\rightarrow} 0,
\end{equation} and
\begin{equation}
\label{6.19}
\frac{1}{n^2}\sum_{i=1}^nw_{i-1}\frac{\widetilde{X}_{i+1}-\widetilde{X}_{i}}{\Delta_n}
-\frac{1}{n^2}\sum_{i=1}^nw^\ast_{i-1}\frac{\widetilde{X}_{i+1}-\widetilde{X}_{i}}{\Delta_n}
\stackrel{p}{\rightarrow} 0.
\end{equation}

\noindent For (\ref{6.18}), we need only prove
\begin{equation}
\label{6.20} \delta_{3,n}=\frac{1}{n}\sum_{i=1}^nK(X_{i-1})
\Big[\frac{\widetilde{X}_{i+1}-\widetilde{X}_{i}}{\Delta_n}-\frac{X_{i}-X_{i-1}}{\Delta_n}\Big]
\stackrel{p}{\rightarrow}0,
\end{equation} and
\begin{equation}
\label{6.21} \delta_{4,n}=\frac{1}{n}\sum_{i=1}^nK(X_{i-1})(X_{i-1}
- x)
\Big[\frac{\widetilde{X}_{i+1}-\widetilde{X}_{i}}{\Delta_n}-\frac{X_{i}-X_{i-1}}{\Delta_n}\Big]
\stackrel{p}{\rightarrow}0.
\end{equation}

\noindent the proof of (\ref{6.20}) and (\ref{6.21}) are similar, so
we just prove (\ref{6.21})

By Lemma \ref{r6.1}, we can get
\begin{eqnarray*}
E[\varepsilon_{1,n}] & := & E\Big[\Big(\frac{(\widetilde{X}_{i+1} -
\widetilde{X}_{i})}{\Delta_{n}} - \frac{({X}_{i} -
{X}_{i-1})}{\Delta_{n}}\Big)\Big|\mathscr{F}_{i-1}\Big]\\
& = & E\Big\{E\Big[\Big(\frac{(\widetilde{X}_{i+1} -
\widetilde{X}_{i})}{\Delta_{n}} - \frac{({X}_{i} -
{X}_{i-1})}{\Delta_{n}}\Big)\Big|\mathscr{F}_{i}\Big]\Big|\mathscr{F}_{i-1}\Big\}\\
& = & \frac{\Delta_n}{2}\Big\{\mu(X_{i-1})\mu^{'}(X_{i-1})+\frac{1}{2}\sigma^2(X_{i-1})\mu^{''}(X_{i-1})\\
&~& + \int_{\mathscr{E}}\big\{\mu(X_{i-1}+ c(X_{i-1} , z)) -
\mu(X_{i-1}) - \mu^{'}(X_{i-1})\cdot c(X_{i-1},z)\big\}f(z)dz
\Big\},
\end{eqnarray*}
so by stationarity and assumption \ref{a3}, we have
\begin{eqnarray*}
E[\delta_{4,n}(x)] & = & E\Big\{E\Big[K(X_{i-1})(X_{i-1} -
x)\Big(\frac{(\widetilde{X}_{i+1} - \widetilde{X}_{i})}{\Delta_{n}}
- \frac{({X}_{i} -
{X}_{i-1})}{\Delta_{n}}\Big)|\mathscr{F}_{i-1}\Big]\Big\}\\
& = & \frac{\Delta_{n}}{2}E\Big[K(X_{i-1})(X_{i-1} -
x)\Big(\mu(X_{i-1})\mu^{'}(X_{i-1}) +\frac{1}{2}\sigma^{2}(X_{i-1})\mu^{''}(X_{i-1})\\
& ~ &+ \int_{\mathscr{E}}\big\{\mu(X_{i-1} + c(X_{i-1} , z)) -
\mu(X_{i-1}) - \mu^{'}(X_{i-1})\cdot
c(X_{i-1},z)\big\}f(z)dz\Big)\Big]\\
& = & O(\Delta_{n})
\end{eqnarray*}
and
\begin{eqnarray*}
&&Var[\delta_{4,n}(x)] \\
& = & \frac{1}{n h_{n}
\Delta_{n}}Var\Big[\frac{1}{\sqrt{n}}\sum_{i=1}^{n}{h^{1/2}_{n}K(X_{i-1})\sqrt{\Delta_{n}}(X_{i-1}
- x)\Big(\frac{(\widetilde{X}_{i+1} -
\widetilde{X}_{i})}{\Delta_{n}} - \frac{({X}_{i\Delta_{n}} -
{X}_{i-1})}{\Delta_{n}}\Big)}\Big]\\
& =: & \frac{1}{n h_{n}
\Delta_{n}}Var\Big[\frac{1}{\sqrt{n}}\sum_{i=1}^{n}g_{i}\Big].
\end{eqnarray*}
\noindent By the similar analysis as above, we can easily obtains
$Var[\delta_{4,n}(x)]\rightarrow 0$  under Assumption \ref{a2} if
$E[g_i^2]<\infty$. In fact by Assumption \ref{a1}, \ref{a3} and
\ref{a4}, we have
\begin{eqnarray*}
E[g_{i}^{2}] & = & E\Big[h_{n} K^{2}(X_{i-1})\Delta_{n}(X_{i-1} -
x)^2\Big(\frac{(\widetilde{X}_{i+1} -
\widetilde{X}_{i})}{\Delta_{n}} - \frac{({X}_{i} -
{X}_{i-1})}{\Delta_{n}}\Big)^{2}\Big]\\
& = & E\Big\{h_{n} K^2(X_{i-1})(X_{i-1} - x)^2
E\Big[\Delta_n\Big(\frac{\tilde{X}_{i+1}-\tilde{X}_{i}}{\Delta_n}-\frac{X_i-X_{i-1}}{\Delta_n}\Big)^2|\mathscr{F}_{i-1}\Big]\Big\}\\
& = & E\Big\{h_{n} K^2(X_{i-1})(X_{i-1} - x)^2
\times \frac{2}{3}\Big[\sigma^2(X_{i-1})+\int_{\mathscr{E}}c^2(X_{i-1},z)f(z)dz+O_P(\Delta_n)\Big]\Big\}\\
& \approx & \Big\{ \begin{array}{ll} O(h^{1/2}_{n})~~~if~x/b \rightarrow \infty~(``interior~x'');\\
O(1)~~~if~x/b \rightarrow \kappa~(``boundary~x'')
\end{array}\\
& < & \infty.
\end{eqnarray*}

The proof of (4.10) is similar to that of (\ref{6.19}). Combination
(\ref{6.10}) and (\ref{6.17}), the relationship $\mu^\ast
_n(x)-\hat{\mu}_n(x) \stackrel{p}{\rightarrow} 0$ holds,
 so by Lemma \ref{r6.1} we have $\hat{\mu}_n(x) \stackrel{p}{\rightarrow} \mu(x).$\end{proof}

\subsection{The proof of Theorem \ref{thm2}}

\begin{proof}
Here we only prove the result for $\mu(x)$; the other is analogical.

\noindent By Lemma \ref{l2}, for ``interior x'', if $h =
O((n\Delta_{n})^{-2/5}),$ then $$U_n^\ast(x) := \sqrt{n\Delta_{n}
h^{1/2}}\big(\mu_{n}^{\ast}(x) - \mu(x) - h
B_{\mu_{n}^{\ast}(x)}\big) \stackrel{d}{\rightarrow} N\Big( 0 ,
\frac{M(x)}{2 \sqrt{\pi} x^{1/2} p(x)}\Big),$$
 for ``boundary x'',
if $h = O((n\Delta_{n})^{-1/5}),$ then $$ U_n^\ast(x) :=
\sqrt{n\Delta_{n} h}\big(\mu_{n}^{\ast}(x) - \mu(x) - h^{2}
B^{'}_{\mu_{n}^{\ast}(x)}\big) \stackrel{d}{\rightarrow} N\Big( 0 ,
\frac{M(x) \Gamma(2\kappa + 1)}{2^{2\kappa + 1} \Gamma^{2}(\kappa +
1) p(x)}\Big),$$ where
$B_{\mu_{n}^{\ast}(x)},~B^{'}_{\mu_{n}^{\ast}(x)}$ denotes the bias
of the estimators of $\mu_{n}^{\ast}(x),$ respectively, that is
$$ B_{\mu_{n}^{\ast}(x)} = \frac{x}{2}\mu^{''}(x),~~~~~B^{'}_{\mu_{n}^{\ast}(x)} = \frac{1}{2}(2 +
\kappa)\mu^{''}(x).$$ So by the asymptotic equivalence theorem, it
suffices to prove that
$$\hat{U}_n(x)-U_n^\ast(x) = \sqrt{h_{n} n \Delta_n }(\hat{\mu}_{n,T} (x)-\mu_{n,T}^\ast (x))\stackrel{p}{\rightarrow} 0,$$
where $\hat{U}_n(x) := \sqrt{n\Delta_{n} h}\big(\hat{\mu}_{n}(x) -
\mu(x) - h^{2} B_{\hat{\mu}_{n}(x)}\big)$ or $\sqrt{n\Delta_{n}
h}\big(\hat{\mu}_{n}(x) - \mu(x) - h^{2}
B^{'}_{\hat{\mu}_{n}(x)}\big).$

In fact, from the proof of Theorem \ref{thm1} such as (\ref{6.10})
and (\ref{6.17}), we know that
\begin{eqnarray*}
& ~ & \hat{U}_n(x)-U_n^\ast(x)\\
& = &  \sqrt{h_{n} n \Delta_n }(\hat{\mu}_{n,T} (x)-\mu_{n,T}^\ast (x))\\
& = &  \sqrt{h_{n} n \Delta_n
}\left(\frac{\delta_n}{\frac{1}{n^2}\sum_{i=1}^{n}w^\ast_{i-1}}
+\frac{\sum_{i=1}^n w_{i-1}
\left(\frac{\widetilde{X}_{i+1}-\widetilde{X}_i}{\Delta_n}\right)}{\sum_{i=1}^{n}w_{i-1}}
-\frac{\sum_{i=1}^n w_{i-1} \left(\frac{\widetilde{X}_{i+1}-\widetilde{X}_i}{\Delta_n}\right)}{\sum_{i=1}^{n}w^{\ast}_{i-1}}\right)\\
& = & \sqrt{h_{n} n \Delta_n
}\left(\frac{\delta_n}{\frac{1}{n^2}\sum_{i=1}^{n}w^\ast_{i-1}}\right)+o_p(1).
\end{eqnarray*}
Due to the stationary case of Lemma 3 in Hanif \cite{hm1}, we have
\begin{equation}
\label{6.22} \frac{1}{n}\sum_{i=1}^n K(X_{i-1})
\stackrel{p}{\rightarrow} p(x).
\end{equation}
We first write
\begin{eqnarray*}
& ~ & \frac{1}{n}\sum_{i=1}^n K(X_{i-1})(X_{i-1} - x)\\
& = & \frac{1}{n}\sum_{i=1}^n E\big[K(X_{i-1})(X_{i-1} - x)\big] +
\frac{1}{n}\sum_{i=1}^n K(X_{i-1})(X_{i-1} - x) -
\frac{1}{n}\sum_{i=1}^n E\big[K(X_{i-1})(X_{i-1} - x)\big]\\
& := & \frac{1}{n}\sum_{i=1}^n E\big[K(X_{i-1})(X_{i-1} - x)\big] +
\frac{1}{n}\sum_{i=1}^n \eta_{i - 1}.
\end{eqnarray*}
According to the result (A.2) in Chen (\cite{ch3}, P321), it is
shown that
\begin{eqnarray*}
& ~ & \frac{1}{n}\sum_{i=1}^n E\big[K(X_{i-1})(X_{i-1} - x)\big]\\
& = & p(x)E(\xi - x) + p^{'}(x)E(\xi - x)^{2} + o_{p}(E(\xi - x)^{2})\\
& = & p(x)h_{n} + p^{'}(x)h_{n}(x + 2 h_{n}) + o_{p}(h_{n}),
\end{eqnarray*}
where $\xi \stackrel{\mathcal {D}} = G (x/h + 1, h)$ and $G$ denotes
the Gamma distribution.

\noindent With the same procedure as the proof details of Lemma 3.3
in Lin, Song and Yi \cite{lsy}, we can prove
$\frac{1}{n}\sum_{i=1}^n \eta_{i - 1} \stackrel{a.s.} \rightarrow
0.$ Hence, we get
\begin{equation}
\label{6.23} \frac{1}{n}\sum_{i=1}^n K(X_{i-1})(X_{i-1} - x)
\stackrel{p}{\rightarrow} p(x)h_{n} + p^{'}(x)h_{n}(x + 2 h_{n}).
\end{equation}

\noindent In the similar procedure, we can obtain
\begin{equation}
\label{6.24} \frac{1}{n}\sum_{i=1}^n K(X_{i-1})(X_{i-1} - x)^2
\stackrel{p}{\rightarrow} p(x)h_{n}(x + 2 h_{n}).
\end{equation}

\noindent According the results of (\ref{6.22}) - (\ref{6.24}), we
have
\begin{eqnarray*}
& ~ & \frac{1}{n^2}\sum_{i=1}^n w^\ast_{i-1}\\ & = & \frac{1}{n^2}
\sum_{i=1}^n K(X_{i-1})\left(\sum_{j=1}^n K(X_{j-1})(X_{j-1} -
x)^2-(X_{i-1} - x)
\sum_{j=1}^n K(X_{j-1})(X_{j-1} - x)\right)\\
& = & \frac{1}{n^2} \sum_{i=1}^n K(X_{i-1})\sum_{j=1}^n
K(X_{j-1})(X_{j-1} - x)^2
- \frac{1}{n^2}\left(\sum_{i=1}^n K(X_{i-1})(X_{i-1} - x)\right)^2\\
& \stackrel{p}{\rightarrow} & h_{n} p^{2}(x) \cdot (x + 2h_{n}) +
h^{2}_{n}[p(x) + p^{'}(x) \cdot x]^{2}\\
& = & h_{n} p^{2}(x) \cdot x + o_{p}(h_{n}).
\end{eqnarray*}
Hence,
$$\hat{U}_n(x)- U^\ast_n(x)=\sqrt{h_{n} n \Delta_n }O_p\Big(\frac{\Delta_n}{h_{n}}\Big)
\stackrel{p}{\rightarrow}0$$ by assumption \ref{a5}.
\end{proof}

\end{document}